\documentclass[11pt]{article}
\usepackage[mathscr]{eucal}
\DeclareMathAlphabet{\mathpzc}{OT1}{pzc}{m}{it}
\newcommand{\W}{{\textsf{{W}$_\sharp^{2}$}}}
\newcommand{{\M}}{{\textsf{{L}$_\sharp^{2}$}}}

\newcommand{\ts}{{\tilde{\sf s}}}

\newcommand{\bsfW}{\boldsymbol{\mathsf W}}
\newcommand{\bsfU}{\boldsymbol{\mathsf U}}
\newcommand{\bsff}{\boldsymbol{\mathsf f}}
\newcommand{\bsfF}{\boldsymbol{\mathsf F}}

\newcommand{\bsfu}{\boldsymbol{\mathsf u}}
\newcommand{\bsfz}{\boldsymbol{\mathsf z}}
\newcommand{\bsfw}{\boldsymbol{\mathsf w}}
\newcommand{\bsfe}{\boldsymbol{\mathsf e}}
\newcommand{\comp}{{\mathbb C}}
\newcommand{\comps}{{\mbox{\tiny $\comp$}}}
\usepackage{bbm}
\usepackage{dsfont}
\usepackage{color}
\usepackage{mathrsfs}
\usepackage{amssymb}
\usepackage{calligra}
\usepackage{fontenc}
\usepackage{amsbsy}

\newcommand{\bfsf}[1]{{\textbf{\textsf{#1}}}}

\usepackage{graphicx}
\graphicspath{%
    {converted_graphics/}
    {/}
}
\begin{document}

\newcommand{\sfU}{{\sf U}}
\newcommand{\sfF}{{\sf F}}
\newcommand{\sfW}{{\sf W}}
\newcommand{\sfD}{{\sf D}}
\newcommand{\sfL}{{\sf L}}
\newcommand{\sfu}{{\sf u}}
\newcommand{\sfw}{{\sf w}}
\newcommand{\sfz}{{\sf z}}
\newcommand{\sfp}{{\sf p}}
\newcommand{\sfv}{{\sf v}}
\newcommand{\simge}{\ba{cc}\vspace*{-2.4mm}>\\ \sim\ea }
\newcommand{\simle}{\ba{cc}\vspace*{-2.4mm}<\\ \sim\ea }
\newcommand{\Cdot}{\!\cdot\!}
\newcommand{\sq}{{$\sqcap\!\!\!\!\sqcup$}}
\newcommand{\Eu}{{\rm I\,\!\! E}}
\newcommand{\Io}{\Int{\Omega}{}}
\newcommand{\Id}{\Int{\cald}{}}
\newcommand{\Div}{\mbox{\rm div}\,}
\newcommand{\tr}{\mbox{\rm tr}\,}
\newcommand{\grad}{\mbox{\rm grad}\,}
\newcommand{\supp}{\mbox{\rm supp}\,}
\newcommand{\curl}{\mbox{\rm curl}\,}
\newcommand{\Ido}{\Int{\partial\Omega}{}}
\newcommand{\IdS}{\Int{\Sigma}{}}
\newcommand{\Oint}[2]{{\displaystyle \oint_{#1}^{#2}}}
\newcommand{\Int}[2]{{\displaystyle \int_{ #1}^{ #2}}}
\newcommand{\Lim}[1]{{\displaystyle \lim_{ #1}}}
\newcommand{\Limsup}[1]{{\displaystyle \limsup_{\footnotesize #1}}}
\newcommand{\Liminf}[1]{{\displaystyle \liminf_{\footnotesize #1}}}
\newcommand{\Sup}[1]{{\displaystyle \sup_{#1}}}
\newcommand{\Inf}[1]{{\displaystyle \inf_{#1}}}
\newcommand{\Max}[1]{{\displaystyle \max_{#1}}}
\newcommand{\Min}[1]{{\displaystyle \min_{#1}}}
\newcommand{\Sum}[2]{{\displaystyle \sum_{#1}^{#2}}}
\newcommand{\Prod}[2]{{\displaystyle \prod_{#1}^{#2}}}
\newcommand{\BCup}[2]{{\displaystyle \bigcup_{#1}^{#2}}}
\newcommand{\BCap}[2]{{\displaystyle \bigcap_{#1}^{#2}}}
\newcommand{\Frac}[2]{\displaystyle{\frac{\displaystyle{#1}}{\displaystyle{#2}}}}
\newcommand{\norm}[1]{\left\|{#1}\right\|}
\newcommand{\Norm}[1]{\langle\langle{#1}\rangle\rangle_q}
\newcommand{\No}[1]{\langle\!\langle{#1}\rangle\!\rangle}
\newcommand{\NO}[1]{{\langle{#1}\rangle}_{\lambda,q}}
\newcommand{\beea}{\begin{eqnarray}}
\newcommand{\eeea}{\end{eqnarray}}
\newcommand{\ms}{\medskip\smallskip}
\newcommand{\bs}{\bigskip}
\newcommand{\ps}{\par\smallskip}
\newcommand{\bfe}{{\mbox{\boldmath $e$}} }
\newcommand{\pni}{{\par\noindent}}
\newcommand{\bfq}{{\mbox{\boldmath $q$}} }
\newcommand{\bfz}{{\mbox{\boldmath $z$}} }
\newcommand{\0}{{\mbox{\boldmath $0$}} }
\newcommand{\LE}{\!\!\!&\le&\!\!\!}
\newcommand{\BL}[1]{{\par\smallskip{\bf Lemma #1.}}}
\newcommand{\BT}[1]{{\par\smallskip{\bf Theorem #1.}}}
\newcommand{\Ln}{[\!|}
\newcommand{\Rn}{|\!]}
\newcommand{\n}[1]{{\Ln{#1}\Rn}} 
\newcommand{\nq}[1]{{\Ln{#1}\Rn}_{q}} 
\newcommand{\nqr}[1]{{\Ln{#1}\Rn}_{q,r}} 
\newcommand{\Nq}[1]{{\langle{#1}\rangle}_{q}} 
\newcommand{\Nql}[1]{{\langle{#1}\rangle}_{\lambda,q}} 
\newcommand{\Nqr}[1]{{\langle{#1}\rangle}_{q,r}}
\newcommand{\N}[1]{{|\!\!|\!\!|\,{#1}\,|\!|\!\!|_2}}
\newcommand{\EA}[2]{$$#1$$%
\vspace{-6.mm}
\begin{equation}
\end{equation}
\vspace{-6.mm}
$$
#2
\setlength{\belowdisplayskip}{3mm}
\setlength{\belowdisplayshortskip}{3mm}
$$
}
\newcommand{\A}[2]{$$#1$$%
\vspace{-4.mm}
$$
#2
\setlength{\belowdisplayskip}{3mm}
\setlength{\belowdisplayshortskip}{3mm}
$$
}
\newcommand{\BF}{\begin{footnotesize}}
\newcommand{\EF}{\end{footnotesize}}
\setlength{\jot}{.15in}
\newcommand{\pde}[2]{{\displaystyle \frac{\mbox{$\partial #1$}}{\mbox{$\partial #2$}}}}
\newcommand{\ode}[2]{{\displaystyle \frac{\mbox{$d #1$}}{\mbox{$d #2$}}}}
\newcommand{\f}[2]{\frac{\mbox{$#1$}}{\mbox{$ #2$}}}
\newcommand{\bi}{\begin{itemize}}
\newcommand{\ei}{\end{itemize}}
\newcommand{\ed}{\end{document}}
\newcommand{\be}{\begin{equation}}
\newcommand{\ba}{\begin{array}}
\newcommand{\ea}{\end{array}}
\newcommand{\ee}{\end{equation}}
\newcommand{\eeq}[1]{\label{eq:#1}\end{equation}}
\newcommand{\real}{{\mathbb R}}
\newcommand{\compl}{{\mathbb C}}
\def\Id{\mbox{\boldmath $1$}}
\def\zero{\mbox{\boldmath $0$}}
\newcommand{\PP}{{\rm I\!\!\,P}}
\newcommand{\nat}{{\mathbb N}}
\newcommand{\bfpsi}{\mbox{\boldmath $\psi$}}
\newcommand{\bfomega}{\mbox{\boldmath $\omega$}}
\newcommand{\bfvaromega}{\mbox{\boldmath $\varpi$}}
\newcommand{\bfOmega}{\mbox{\boldmath $\Omega$}}
\newcommand{\bfTheta}{\mbox{\boldmath $\Theta$}}
\newcommand{\bfmu}{\mbox{\boldmath $\mu$}}
\newcommand{\bfx}{\mbox{\boldmath $x$}}
\newcommand{\bfy}{\mbox{\boldmath $y$}}
\newcommand{\bfPsi}{\mbox{\boldmath $\Psi$}}
\newcommand{\bfxi}{\mbox{\boldmath $\xi$}}

\newcommand{\bfphi}{\mbox{\boldmath $\varphi$}}
\newcommand{\bfhi}{\mbox{\boldmath $\phi$}}
\newcommand{\bfPhi}{\mbox{\boldmath $\Phi$}}
\newcommand{\bfv}{{\mbox{\boldmath $v$}} }
\newcommand{\bfu}{{\mbox{\boldmath $u$}} }
\newcommand{\bfuf}{{\mbox{\footnotesize\boldmath $u$}} }
\newcommand{\bfw}{{\mbox{\boldmath $w$}} }
\newcommand{\bff}{{\mbox{\boldmath $f$}} }
\newcommand{\bfa}{{\mbox{\boldmath $a$}} }
\newcommand{\bfi}{{\mbox{\boldmath $i$}} }
\newcommand{\bfj}{{\mbox{\boldmath $j$}} }
\newcommand{\bfc}{{\mbox{\boldmath $c$}} }
\newcommand{\bfo}{{\mbox{\boldmath $o$}} }
\newcommand{\bfp}{{\mbox{\boldmath $p$}} }
\newcommand{\bfkp}{{\mbox{\footnotesize{\boldmath $k$}}} }
\newcommand{\bfka}{{\mbox{\footnotesize{\boldmath $k^*$}}} }
\newcommand{\bft}{{\mbox{\boldmath $t$}} }
\newcommand{\bfd}{{\mbox{\boldmath $d$}} }
\newcommand{\bfl}{{\mbox{\boldmath $l$}} }
\newcommand{\bfr}{{\mbox{\boldmath $r$}} }
\newcommand{\bfk}{{\mbox{\boldmath $k$}} }
\newcommand{\bfA}{{\mbox{\boldmath $A$}} }
\newcommand{\bfS}{{\mbox{\boldmath $S$}} }
\newcommand{\bfO}{{\mbox{\boldmath $O$}} }
\newcommand{\bfM}{{\mbox{\boldmath $M$}} }
\newcommand{\bfP}{{\mbox{\boldmath $P$}} }
\newcommand{\bfB}{{\mbox{\boldmath $B$}} }
\newcommand{\bfR}{{\mbox{\boldmath $R$}} }
\newcommand{\bfC}{{\mbox{\boldmath $C$}} }
\newcommand{\bfD}{{\mbox{\boldmath $D$}} }
\newcommand{\bfQ}{{\mbox{\boldmath $Q$}} }
\newcommand{\bfZ}{{\mbox{\boldmath $Z$}} }
\newcommand{\bfG}{{\mbox{\boldmath $G$}} }
\newcommand{\bfE}{{\mbox{\boldmath $E$}} }
\newcommand{\bfX}{{\mbox{\boldmath $X$}} }
\newcommand{\bfY}{{\mbox{\boldmath $Y$}} }
\newcommand{\bfH}{{\mbox{\boldmath $H$}} }
\newcommand{\bfI}{{\mbox{\boldmath $I$}} }
\newcommand{\bfJ}{{\mbox{\boldmath $J$}} }
\newcommand{\bfN}{{\mbox{\boldmath $N$}} }
\newcommand{\bfh}{{\mbox{\boldmath $h$}} }
\newcommand{\bfm}{{\mbox{\boldmath $m$}} }
\newcommand{\bfone}{{\mbox{\boldmath $1$}} }
\newcommand{\hs}{{\rm I}\!\!\,{\rm R}^3_+}
\newcommand{\cala}{{\cal A}}
\newcommand{\calb}{{\cal B}}
\newcommand{\calc}{{\cal C}}
\newcommand{\cald}{{\cal D}}
\newcommand{\cale}{{\cal E}}
\newcommand{\calf}{{\cal F}}
\newcommand{\calg}{{\cal G}}
\newcommand{\calh}{{\cal H}}
\newcommand{\cali}{{\cal I}}
\newcommand{\calj}{{\cal J}}
\newcommand{\calk}{{\cal K}}
\newcommand{\call}{{\cal L}}
\newcommand{\calm}{{\cal M}}
\newcommand{\caln}{{\cal N}}
\newcommand{\calo}{{\cal O}}
\newcommand{\calp}{{\cal P}}
\newcommand{\calq}{{\cal Q}}
\newcommand{\calr}{{\cal R}}
\newcommand{\cals}{{\cal S}}
\newcommand{\calt}{{\cal T}}
\newcommand{\calu}{{\cal U}}
\newcommand{\calv}{{\cal V}}
\newcommand{\calx}{{\cal X}}
\newcommand{\caly}{{\cal Y}}
\newcommand{\calw}{{\cal W}}
\newcommand{\calz}{{\cal Z}}
\newcommand{\bfsigma}{\mbox{\boldmath $\sigma$}}
\newcommand{\bfSigma}{\mbox{\boldmath $\Sigma$}}
\newcommand{\bftau}{\mbox{\boldmath $\tau$}}
\newcommand{\bfeta}{\mbox{\boldmath $\eta$}}
\newcommand{\bfT}{{\mbox{\boldmath $T$}} }
\newcommand{\bfV}{{\mbox{\boldmath $V$}} }
\newcommand{\bfU}{{\mbox{\boldmath $U$}} }
\newcommand{\bfW}{{\mbox{\boldmath $W$}} }
\newcommand{\bfF}{{\mbox{\boldmath $F$}} }
\newcommand{\bfK}{{\mbox{\boldmath $K$}} }
\newcommand{\bfL}{{\mbox{\boldmath $L$}} }
\newcommand{\bfb}{{\mbox{\boldmath $b$}} }
\newcommand{\bfg}{{\mbox{\boldmath $g$}} }
\newcommand{\bfn}{{\mbox{\boldmath $n$}} }
\newcommand{\bfs}{{\mbox{\boldmath $s$}} }
\newcommand{\cf}{{\it cf.} }
\newcommand{\io}{\int_\Omega}
\newcommand{\1}{\item[({\it i})]}
\newcommand{\2}{\item[({\it ii})]}
\newcommand{\3}{\item[({\it iii})]}
\newcommand{\4}{\item[({\it iv})]}
\newcommand{\5}{\item[({\it v})]}
\newcommand{\6}{\item[({\it vi})]}
\newcommand{\7}{\item[({\it vii})]}
\newcommand{\8}{\item[({\it viii})]}
\newcommand{\9}{\item[({\it xi})]}
\newcommand{\ido}{\int_{\partial\Omega}}
\newcommand{\half}{\mbox{$\frac{1}{2}$}}
\def\parallel{\|}
\def\mid{|}
\def\Bbb R{\real}
\def\hat{\widehat}
\def\tilde{\widetilde}
\def\bar{\overline}
\newcommand{\threehalves}{3\over 2}
\newcommand{\bfPi}{\mbox{\boldmath $\Pi$}}
\newcommand{\bfchi}{\mbox{\boldmath $\chi$}}
\newcommand{\bfalpha}{\mbox{\boldmath $\alpha$}}
\newcommand{\bfbeta}{\mbox{\boldmath $\beta$}}
\newcommand{\bfgamma}{\mbox{\boldmath $\gamma$}}
\newcommand{\bfdelta}{\mbox{\boldmath $\delta$}}
\newcommand{\bfzeta}{\mbox{\boldmath $\zeta$}}
\newcommand{\bfUpsilon}{\mbox{\boldmath $\Upsilon$}}
\newcommand{\bfGamma}{\mbox{\boldmath $\Gamma$}}
\newcommand{\bfcala}{\mbox{\boldmath ${\cal A}$}}
\newcommand{\bfcalm}{\mbox{\boldmath ${\cal M}$}}
\newcommand{\bfcaln}{\mbox{\boldmath ${\cal N}$}}
\newcommand{\bfcalq}{\mbox{\boldmath ${\cal Q}$}}
\newcommand{\bfcalb}{\mbox{\boldmath ${\cal B}$}}
\newcommand{\bfcalc}{\mbox{\boldmath ${\cal C}$}}
\newcommand{\bfcali}{\mbox{\boldmath ${\cal I}$}}
\newcommand{\bfcalg}{\mbox{\boldmath ${\cal G}$}}
\newcommand{\bfcalh}{\mbox{\boldmath ${\cal H}$}}
\newcommand{\bfcalk}{\mbox{\boldmath ${\cal K}$}}
\newcommand{\bfcalt}{\mbox{\boldmath ${\cal T}$}}
\newcommand{\bfcalx}{\mbox{\boldmath ${\cal X}$}}
\newcommand{\bfcall}{\mbox{\boldmath ${\cal L}$}}
\newcommand{\bfcalf}{\mbox{\boldmath ${\cal F}$}}
\newcommand{\bfcalr}{\mbox{\boldmath ${\cal R}$}}
\newcommand{\bfcals}{\mbox{\boldmath ${\cal S}$}}
\newcommand{\bfcalw}{\mbox{\boldmath ${\cal W}$}}
\newcommand{\bfcalu}{\mbox{\boldmath ${\cal U}$}}
\newcommand{\bfcalv}{\mbox{\boldmath ${\cal V}$}}
\newcommand{\bfcalz}{\mbox{\boldmath ${\cal Z}$}}
\pagenumbering{roman}
\newcommand{\art}[6]{{\I[{\sc #1,}] {#2}, {\it #3}, {\bf #4}, {#5} {[#6]}}}
\newcommand{\ED}{\end{description}}
\newcommand{\I}{\item }
\newcommand{\ra}{\rm a}
\newcommand{\rb}{\rm b}
\newcommand{\rc}{\rm c}
\newcommand{\Hsp}{{\rm I}\!\!\,{\rm R}^n_+}
\newcommand{\Hsn}{{\rm I}\!\!\,{\rm R}^n_-}
\newcommand{\po}[1]{\mbox{$\displaystyle \frac{\mbox{$\partial #1$}}
{\mbox{$\partial x_{1}$}}$}}
\newcommand{\PO}[1]{\mbox{$\displaystyle \frac{\mbox{$\partial #1$}}
{\mbox{$\partial y_{1}$}}$}}
\newcommand{\OP}{\left(\Delta+2\lambda\PO{}\right)}
\newcommand{\op}{\left(\Delta+2\lambda\po{}\right)}
\newcommand{\ft}[1]{
\Frac{1}{(2\pi)^{n/2}}\Int{{\Bbb R}^{n}}{}e^{i{\bf x}\cdot \bfchi}
#1(\xi)d\xi}
\newcommand{\Ft}[1]{
\Frac{1}{2\pi}\Int{{\Bbb R}^{2}}{}e^{i{x}\cdot \xi}
#1(\xi)d\xi}
\newcommand{\Z}{\item[({\it a})]}
\newcommand{\B}{\item[({\it b})]}
\newcommand{\D}{\item[({\it d})]}
\newcommand{\E}{\item[({\it e})]}
\newcommand{\G}{\item[({\it g})]}
\newcommand{\Š}{\`e}
\newcommand{\…}{\`a}
\newcommand{\•}{\`o}
\newcommand{\—}{\`u}
\newcommand{\}{\`{\i}}
\def\tag{\renewcommand{\theequation}}
\newcommand{\Footnote}{~\footnote}
\newcommand{\ie}{{\it i.e.}}
\newcommand{\dist}{\mbox{\rm dist\,}}
\newcommand{\const}{\mbox{\rm const}}
\newcommand{\trace}{\mbox{\rm trace}}
\newcommand{\Bo}{\par\hfill{$\Box$}\par\noindent}
\newcommand{\Nor}[1]{\langle{#1}\rangle_q}
\newcommand{\vs}{\vspace*{.5cm}\par\noindent}
\newcommand{\Vs}{\vspace*{.6cm}\par\noindent}
\newcommand{\Vvs}{\vspace*{.7cm}\par\noindent}
\newcommand{\VVs}{\vspace*{.8cm}\par\noindent}
\newtheorem{definition}{Definition}[section]
\newcommand{\Bd}{\begin{definition}\begin{rm}}
\newcommand{\Ed}{\end{rm}\end{definition}}
\newtheorem{remark}{Remark}[section]
\newcommand{\Br}{\begin{remark}\begin{rm}}
\newcommand{\Er}{\end{rm}\end{remark}}
\newtheorem{proposition}{Proposition}[section]
\newcommand{\Bp}{\begin{proposition}\begin{sl}}
\newcommand{\EP}[1]{\end{sl}\label{proposition:#1}\end{proposition}}
\newcommand{\propref}[1]{{\rm Proposition \ref{proposition:#1}}}
\newcommand{\Bt}{\begin{theorem}\begin{sl}}
\newcommand{\Et}{\end{sl}\end{theorem}}
\newcommand{\Bl}{\begin{lemma}\begin{sl}}
\newcommand{\El}{\end{sl}\end{lemma}}
\newtheorem{theorem}{Theorem}[section]
\newtheorem{lemma}{Lemma}[section]
\newtheorem{lemmaA}{Lemma A.}
\newtheorem{corollary}{Corollary}[section]
\newcommand{\eqref}[1]{{\rm (\ref{eq:#1})}}
\newcommand{\Bc}{\begin{corollary}\begin{sl}}
\newcommand{\Ec}{\end{sl}\end{corollary}}
\newcommand{\ET}[1]{\end{sl}\label{theorem:#1}\end{theorem}}
\newcommand{\EDD}[1]{\end{rm}\label{definition:#1}\end{definition}}
\newcommand{\EL}[1]{\end{sl}\label{lemma:#1}\end{lemma}}
\newcommand{\theoref}[1]{{\rm Theorem \ref{theorem:#1}}}
\newcommand{\ER}[1]{\end{rm}\label{remark:#1}\end{remark}}
\newcommand{\EC}[1]{\end{sl}\label{corollary:#1}\end{corollary}}
\newcommand{\remref}[1]{{\rm Remark \ref{remark:#1}}}
\newcommand{\cororef}[1]{{\rm Corollary \ref{corollary:#1}}}
\newcommand{\lemmref}[1]{{\rm Lemma \ref{lemma:#1}}}
\newcommand{\essup}[1]{{\rm ess}\,{{\displaystyle \sup_{\hspace*{-5mm}{#1}}}}}

\renewcommand{\i}{{\rm i}}

\pagenumbering{arabic}
\newcommand{\QED}{{\par\hfill$\square$\par}}
\renewcommand{\thefootnote}{(\arabic{footnote})}
\title{Mathematical Analysis of Flow-Induced Oscillations of a   Spring-Mounted Body in a Navier-Stokes Liquid} 
\author{ Giovanni P. Galdi 
\thanks{Department of Mechanical Engineering and Materials Science, University of Pittsburgh, PA 15261. 
}}
\date{}
\maketitle
\begin{abstract} We study the motion of a rigid body $\mathscr B$ subject to an undamped elastic restoring force, in the stream  of a viscous liquid $\mathscr L$. The motion of the coupled system $\mathscr B$-$\mathscr L\equiv\mathscr S$ is driven by a  uniform flow of $\mathscr L$ at spatial infinity, characterized by a given, constant dimensionless velocity $\lambda\,\bfe_1$, $\lambda>0$. We show that as long as $\lambda\in(0,\lambda_c)$, with $\lambda_c$ a distinct positive number, there is  a uniquely determined time-independent state of $\mathscr S$ where $\mathscr B$ is in a (locally) stable equilibrium and the flow of $\mathscr L$ is steady. Moreover, in that range of $\lambda$, no oscillatory flow may occur. Successively we prove that if certain suitable spectral properties of the relevant linearized operator are met, there exists a $\lambda_0>\lambda_c> 0$ where  an oscillatory regime for $\mathscr S$ sets in. More precisely, a bifurcating  time-periodic branch stems out of the time-independent solution. The significant feature of this result is that {\em no} restriction is imposed on the frequency, $\omega$, of the bifurcating solution, which may thus coincide with the natural structural frequency, $\omega_{\sf n}$, of $\mathscr B$, or any multiple of it. This implies that a dramatic structural failure cannot take place due to resonance effects. However, our analysis also shows that when $\omega$ becomes sufficiently close to $\omega_{\sf n}$ the amplitude of oscillations can become very large in the limit when the density of $\mathscr L$ becomes negligible compared to that of $\mathscr B$.        

 \end{abstract}
\newpage\vspace*{2.7cm}
{\contentsline{section}{\centerline{{\rm \Large Contents}}}{}\vspace*{1.0cm}
\contentsline {section}{{\rm Introduction}}{{\rm 3}}
\contentsline {section}{\numberline {{\rm 1.}}{\rm Formulation of the Problem}}{{\rm 5}}
\contentsline {section}{\numberline {{\rm 2.}}{\rm Function Spaces and their Relevant Properties}}{{\rm 6}}
\contentsline {section}{\numberline {{\rm 3.}}{\rm Steady-States: Existence and Uniqueness} }{{\rm 13}}
\contentsline {section}{\numberline {{\rm 4.}}{\rm Steady-States: Stability}}{{\rm 14}}
\contentsline {section}{\numberline {{\rm 5}}{\rm On the Absence of Oscillatory Regimes}}{{\rm 21}}
\contentsline {section}{\numberline {{\rm 6.}}{\rm A General Approach to Time-Periodic Bifurcation}}{{\rm 24}}
\contentsline {section}{\numberline {{\rm 7.}}{\rm Existence of an Analytic Steady-State Branch}}{{\rm 27}}
\contentsline {section}{\numberline {{\rm 8.}}{\rm Spectral Properties of the Linearized Operator}}{{\rm 30}}
\contentsline {section}{\numberline {{\rm 9.}}{\rm The Linearized Time-Periodic Operator}}{{\rm 34}}
\contentsline {section}{\numberline {{\rm 10.}}{\rm Reformulation of the Problem in Banach Spaces}}{{\rm 42}}
\contentsline {section}{\numberline {{\rm 11.}}{\rm A Time-Periodic Bifurcation Theorem}}{{\rm 43}}
\contentsline {section}{{\rm Appendix}}{{\rm 45}}
\contentsline {section}{{\rm References}}{{\rm 49}}

\newpage
\renewcommand{\theequation}{\arabic{section}.\arabic{equation}}
\section*{Introduction} \addcontentsline{toc}{section}{Introduction} The flow of a viscous fluid around structures is a  fundamental problem that lies at the heart of the broad research area of Fluid-Solid-Interaction. A main feature of this problem regards the study of the oscillations (vibrations) produced by the fluid on the structure. In fact, they may lead either to useful and pleasant motions, like ringing  wind chimes or Aeolian harps, or else destructive consequences, as damage or even collapse of the structure. In regards to 
latter, of particular significance is the phenomenon of forced oscillation of suspension bridges,
induced by the vortex shedding of the fluid (air), which reflects into an oscillatory regime of
the wake. When the frequency of the wake approaches the natural structural frequency of the
the bridge, a resonant phenomenon may occur that could culminate into structural failure. A very well known and infamous example of this phenomenon is the  collapse of the Tacoma Narrows bridge.
\par
Because of its fundamental importance in many  practical situations, it is not surprising that the problem of flow-induced oscillations of structures has received all along a plethora of contributions by the engineering community, from experimental, numerical and theoretical viewpoints. Already the list of   the  most relevant articles only is too long to be included here, and we direct the reader to the books \cite{Bev, Dyr}, the review article \cite{Will} and the references therein. The structure model typically adopted for this study consists of a rigid body (in most cases, a circular cylinder) subject to a linear restoring elastic force, while the fluid flow is taken to be two-dimensional, with the main stream orthogonal to the axis of the cylinder. 
\par
Notwithstanding, the problem has not yet received a similar, systematic attention from the mathematical community. As a matter of fact, we believe that a {\em rigorous} mathematical approach to the problem   will, on the one hand, lead to a deeper and more
clear understanding of the underlying physics of the phenomenon, and, on the other hand, propose new interesting
and challenging questions that might  be of great interest to both mathematician and engineer. With this in mind, in the joint paper \cite{BGG} we have very recently started a rigorous investigation of flow-induced oscillations. There, we have dealt with the simple, though significant, model problem where a two-dimensional
rectangular structure is subject to a unidirectional restoring elastic force, while immersed in the
two-dimensional channel flow of a Navier-Stokes liquid, driven by a time-independent Poiseuille flow.
The main objective in \cite{BGG} is rather basic and regards the  existence of possible equilibrium configurations of
the structure, at least for ``small" data. Successively, several other papers have been dedicated to the investigation of further relevant properties of this model, such as well-posedness of the relevant initial-boundary value problem \cite{Patri}, large-time behavior and existence of a global attractor \cite{GPP}, and sharp threshold for  uniqueness of the equilibrium configuration \cite{GazP}.       
\par
Objective of this article is to provide an additional, rigorous contribution to the flow-induced oscillations problem. The model we shall consider is somewhat more general than that in \cite{BGG}, and is constituted by a rigid (finite) body, $\mathscr B$, of {\em arbitrary} shape,  subject to a linear {\em undamped} restoring force and immersed in the stream of a Navier--Stokes liquid, $\mathscr L$, that fills the entire {\em three-dimensional} space outside $\mathscr B$. The choice of an undamped spring is because we want to investigate the occurrence of possible resonance phenomena in the worst case scenario. The motion of the coupled system $\mathscr B$-$\mathscr L\equiv\mathscr S$ is driven by a  uniform flow of $\mathscr L$ at large distances from $\mathscr B$ (spatial infinity) that is characterized by a {\em time-independent} velocity field $\bfV$. Denote by $\lambda\, (>0)$ a dimensionless value of the magnitude of $\bfV$ (Reynolds number). Our focus is to study the behavior of $\mathscr S$ as the parameter $\lambda$ increases. Precisely, we begin to show that, for any given $\lambda>0$, there is always a corresponding time-independent (steady-state) motion, ${\sf s}(\lambda)$, of $\mathscr S$, namely, $\mathscr B$ is in equilibrium and the flow of $\mathscr L$ is steady; see Section \ref{sec:exst}.   Successively, we are able to single out a positive value, $\lambda_c$, and show the following properties. If $\lambda<\lambda_c$, then ${\sf s}(\lambda)$   is uniquely determined and (locally) asymptotically stable; see Section \ref{sec:stab}. Moreover, in this range of $\lambda$'s,  oscillatory motions of $\mathscr S$ cannot occur; see Section \ref{sec:nobif}. Therefore, flow-induced oscillations must take place at some $\lambda\equiv\lambda_0>\lambda_c$. We thus aim at finding sufficient conditions on $(\lambda_0,{\sf s}(\lambda_0))$ ensuring that a time-periodic branch does exist in a neighborhood of $(\lambda_0,{\sf s}(\lambda_0))$. This question is investigated in the framework of a general time-periodic bifurcation theory we introduced in \cite{GaBif}, which, for completeness,  we recall in Section \ref{sec:bif}, with full details presented in the Appendix. Unlike classical approaches \cite[Theorem 79.F]{Z1}, ours  is able to handle flows in unbounded regions and, in particular, overcome the notorious issue of 0 being in the essential spectrum of the linearized operator \cite{Bab0,Bab,FN,Saz}. Nevertheless,  the application of the main theorem in \cite{GaBif} to the case at hand is not straightforward, and requires a careful study of the functional and spectral properties of  certain linear operators obtained by linearizing the original equations around $(\lambda_0,{\sf s}(\lambda_0))$. This study is performed in Sections \ref{sec:spectrum} and \ref{sec:time_per}. The {\em crucial} result established in Section \ref{sec:time_per} (see \theoref{3.1}) is that the ``purely oscillatory" linear operator is Fredholm of index 0,   {\em whatever} the (finite) value of the material constants. In other words, the Fredholm property holds also in case of ``resonance," namely, when the frequency of the oscillation coincides with the natural frequency of $\mathscr B$ or a multiple of it. It should be noticed that we also show that this important property no longer holds in the limit situation when the ratio of the density of $\mathscr L$ to that of $\mathscr B$ tends to 0, as somehow expected on physical ground. With these preparatory results in hand, we are thus able to reformulate the bifurcation problem in a suitable functional setting (Section \ref{sec:rif}), and, with the help of the theory given  in \cite{GaBif},  to show our main finding, collected in Section \ref{sec:Bif} (see \theoref{8.1}). In a nutshell, the latter states that, under certain conditions ensuring the existence of an analytic steady-state branch in the neighborhood of $(\lambda_0,{\sf s}(\lambda_0))$ --see Section \ref{steady}-- there is a family of analytic time-periodic solutions emanating from $(\lambda_0,{\sf s}(\lambda_0))$, provided the linearized operator at $(\lambda_0,{\sf s}(\lambda_0))$ has a ``non-resonant" simple, purely imaginary eigenvalue crossing the imaginary axis at non-zero speed, when $\lambda$ crosses $\lambda_0$. Also, as expected, the bifurcation is supercritical. It is worth emphasizing that, because of the range of validity of the Fredholm property mentioned above, no restriction is imposed on the frequency of the bifurcating branch, thus ruling out the occurrence of  ``disruptive" resonance.\footnote
{We recall that, in our model, the restoring force is {\em undamped}, which means that, in absence of liquid, the oscillation amplitude will become unboundedly large, whenever the frequency of the forcing mechanism gets closer and closer to the natural frequency of $\mathscr B$ or a multiple of it.}
\par
Of course, the analysis performed here  is by no means exhaustive and leaves out a number of fundamental questions. First and foremost, the stability properties of the bifurcating solution. On physical ground, it is expected that the time-periodic bifurcating branch becomes stable, while the steady-state looses its stability. However, with our current knowledge, a rigorous proof of this assertion remains out of our reach and will be the object of future thoughts.
\par
We conclude this introductory section by briefly outlining  the plan of the article. In Section \ref{sec:Form} we give the relevant equations and furnish  the mathematical formulation of the problem. In Section \ref{spaces}, after recalling some notation, we introduce the appropriate function spaces and present some of their important characteristics. In the next two sections we investigate the properties of steady-state solutions. Precisely, in Section \ref{sec:exst} we prove their existence for arbitrary $\lambda>0$, and uniqueness for $\lambda\in (0,\lambda_1)$, $\lambda_1>0$, whereas in Section \ref{sec:stab} we show that their are stable for $\lambda\in (0,\lambda_2)$, $0<\lambda_2\le \lambda_1$. These results are complemented with the one given in Section \ref{sec:nobif}, where we show that in the range $\lambda\in (0,\lambda_2)$ no time-periodic motion may occur. The remaining part of the paper is devoted to finding sufficient conditions for the existence of time-periodic flow. As mentioned earlier on, this is done with the approach introduced in \cite{GaBif}, that we summarize in Section \ref{sec:bif} and fully detail in the Appendix. In Section \ref{steady}, we show the existence of an analytic branch of steady-state solutions around $(\lambda_0,{\sf s}(\lambda_0))$, provided the linearization at that point is trivial.  The next two sections are dedicated to the study of the properties of the relevant linear operators. Precisely, in Section \ref{sec:spectrum} we show that the intersection of the spectrum of the (time-independent) linearization at $(\lambda_0,{\sf s}(\lambda_0))$ with the imaginary axis is constituted (at most) by a bounded sequence of eigenvalues of finite algebraic multiplicity that can cluster only at 0. In Section \ref{sec:time_per} we show that the linearized ``purely oscillatory" operator, suitably defined, is Fredholm of index 0, whatever the value of the physical parameters. With all the above results in hand, in Section \ref{sec:rif} we may then formulate the bifurcation problem in a functional setting that fits the requirement of the general theory given in \cite{GaBif}. Thus, in the final Section \ref{sec:Bif}, we prove a theorem that ensures the existence of a supercritical family of time-periodic solutions in the neighborhood of $(\lambda_0,{\sf s}(\lambda_0))$, stemming out of the analytic branch of steady-states.       
\section{Formulation of the Problem}\label{sec:Form}
Let $\mathscr B$ be a rigid body moving in a Navier-Stokes liquid that fills the entire space, $\Omega$, outside $\mathscr B$, and whose flow becomes uniform at ``large" distances from $\mathscr B$, characterized by a  constant velocity, $\bfV$. On $\mathscr B$ an elastic restoring force $\bfF$ acts, applied to its center of mass $G$, while a suitable active torque  prevents it from rotating. Therefore, the motion of $\mathscr B$ is translatory. In this situation, the  governing equations of   motion of the coupled system body-liquid when  referred to a body-fixed frame $\calf\equiv\{G,\bfe_i\}$  are given by \cite[Section 1]{Gah}   
\be\ba{cc}\medskip\left.\ba{ll}\medskip
\partial_t\bfv+(\bfv-{\bfgamma})\cdot\nabla\bfv=\nu\Delta\bfv-\nabla p\\
\Div\bfv=0\ea\right\}\ \ \mbox{in $\Omega\times(0,\infty)$}\,,\\ \medskip
\bfv(x,t)={\bfgamma}(t)\,, \ \mbox{ $(x,t)\in\partial\Omega\times(0,\infty)$}\,;\ \ 
\Lim{|x|\to\infty}\bfv(x,t)=\bfV\,,\ t\in(0,\infty)\,,\\
M\dot{\bfgamma}+\rho\Int{\partial\Omega}{} \mathbb T(\bfv,p)\cdot\bfn=\bfF \ \ \mbox{in $(0,\infty)$}\,.
\ea
\eeq{01}
In \eqref{01}, $\bfv$ and $\rho\,p$ represent velocity and pressure fields of the liquid,  $\rho$ and $\nu$ its density and kinematic viscosity,  while $M$ and $\bfgamma=\bfgamma(t)$ denote mass of $\mathscr B$  and velocity of $G$, respectively. Moreover, 
$$
\mathbb T(\bfz,\psi):=2\nu\,\mathbb D(\bfz)-\psi\,\mathbb I\,,\ \ \ \mathbb D(\bfz):=\half\left(\nabla\bfz+(\nabla\bfz)^\top\right)\,,
$$
with $\mathbb I$ identity matrix, is the Cauchy stress tensor, and, finally, $\bfn$ is the unit outer normal at $\partial\Omega$.\par
We shall assume that $\bfF$ is a linear function\footnote{As will become clear from their proof, our findings (appropriately modified) will continue to hold also for a suitable class of nonlinear forces.}
 of the displacement $\bfchi:=\int\bfgamma(s){\rm d}s$ evaluated with respect to a fixed point, namely$$
\bfF=-\ell\,\bfchi\,,\ \ell>0.
$$
Without loss of generality we take $\bfV=-V\bfe_1$, $V>0$. Thus, scaling  velocity with $V$, length with $L:={\rm diam}\,\mathscr B$,  time with $VL^2/\nu$, and setting $\bfu:=\bfv+\bfe_1$, we may rewrite \eqref{01} in the following form
\be\ba{cc}\medskip\left.\ba{ll}\medskip
\partial_t\bfu-\lambda\,[\partial_1\bfu+(\dot{\bfchi}-\bfu)\cdot\nabla\bfu]=\Delta\bfu-\nabla p\\
\Div\bfu=0\ea\right\}\ \ \mbox{in $\Omega\times(0,\infty)$}\,,\\ \medskip
\bfu(x,t)={\dot{\bfchi}}(t)+\bfe_1\,, \ \mbox{ $(x,t)\in\partial\Omega\times(0,\infty)$}\,;\ \ 
\Lim{|x|\to\infty}\bfu(x,t)=\0\,,\ t\in(0,\infty)\,,\\
\ddot{\bfchi}+\omega_{\sf n}^2\bfchi+\varpi\Int{\partial\Omega}{} \mathbb T(\bfu,p)\cdot\bfn=\0 \ \ \mbox{in $(0,\infty)$}\,,
\ea
\eeq{02}
with
$$
\omega_{\sf n}^2:=\frac{ L^2\ell}{M\nu}\,,\ \ \varpi:=\frac{\rho L^3}{M}\,,\ \ \lambda:=\frac{VL}{\nu}\,,
$$
and where now all the involved quantities are non-dimensional.  
\par
We are interested in the existence of  oscillatory solutions $(\bfu,p,\bfchi)$ to \eqref{02}  bifurcating from a steady-state branch. In order to make this statement more precise,  let ${\sf s}_0=(\bfu_0,p_0,\bfchi_0)$ be a steady-state solution to \eqref{02} corresponding to a given $\lambda$, namely, 
\be\ba{cc}\medskip\left.\ba{ll}\medskip
-\lambda\,(\partial_1\bfu_0-\bfu_0\cdot\nabla\bfu_0)=\Delta\bfu_0-\nabla p_0\\
\Div\bfu_0=0\ea\right\}\ \ \mbox{in $\Omega$}\,,\\ \medskip
\bfu_0(x)=\bfe_1\,, \ \mbox{ $x\in\partial\Omega$}\,;\ \ 
\Lim{|x|\to\infty}\bfu_0(x)=\0\,,\\
\omega_{\sf n}^2\bfchi_0+\varpi\Int{\partial\Omega}{} \mathbb T(\bfu_0,p_0)\cdot\bfn=\0\,.
\ea
\eeq{03}
From the physical viewpoint, $\bfchi_0$ represents the elongation of the spring necessary to keep $\calb$ in place, that is, prevent it from oscillating. As  shown in Sections \ref{sec:exst} and \ref{sec:stab}, for any $\lambda>0$, \eqref{03} has at least one solution ${\sf s}_0(\lambda):=(\bfu_0(\lambda),p_0(\lambda),\bfchi_0(\lambda))$ that is unique and stable provided $\lambda$ remains  below a definite value $\lambda_c$. Moreover, as long as $\lambda<\lambda_c$, any oscillatory regime branching out of a steady-state is excluded; see Section \ref{sec:nobif}. Therefore, we look for a value of $\lambda>\lambda_c$ ($\lambda=\lambda_0$, say) at which a time-periodic motion may indeed set in. 
More precisely, writing
$$
\bfu=\bfu_0(\lambda)+\bfw(t;\lambda)\,,\ p=p_0(\lambda)+{\sf p}(t;\lambda)\,,\ \ \bfchi=\bfchi_0(\lambda)+ \bfeta(t;\lambda)\,,
$$
equations \eqref{02} delivers the first-order in time problem in the unknowns $(\bfw,{\sf p},\bfeta,\bfsigma)$:
\be\ba{cc}\medskip\left.\ba{ll}\medskip
\partial_t\bfw-\lambda\,[\partial_1\bfw-\bfu_0\cdot\nabla\bfw+({\bfsigma}-\bfw)\cdot\nabla\bfu_0 -({\bfsigma}-\bfw)\cdot\nabla\bfw]=\Delta\bfw-\nabla {\sf p}\\
\Div\bfw=0\ea\right\}\ \ \mbox{in $\Omega\times(0,\infty)$}\,,\\ \medskip
\bfw(x,t)={{\bfsigma}}(t)\,, \ \mbox{ $(x,t)\in\partial\Omega\times(0,\infty)$}\,;\ \ 
\Lim{|x|\to\infty}\bfw(x,t)=\0\,,\ t\in(0,\infty)\,,\\
\dot{\bfsigma}+\omega_{\sf n}^2\bfeta+\varpi\Int{\partial\Omega}{} \mathbb T(\bfw,{\sf p})\cdot\bfn=\0\,;\ \ \dot{\bfeta}(t)=\bfsigma(t) \ \ \mbox{in $(0,\infty)$}\,.
\ea
\eeq{04}
Then, our objective is to determine sufficient conditions for the existence
of a non-trivial family of time-periodic solutions to \eqref{04}, $(\bfw(\lambda), {\sf p}(\lambda),\bfeta(\lambda),\bfsigma(\lambda))$,
$\lambda\in U(\lambda_0)$, of (unknown) frequency $\zeta =\zeta (\lambda)$,  such that $(\bfw(\lambda), {\sf p}(\lambda),\bfeta(\lambda),\bfsigma(\lambda))\to (\0,\0,\0,\0)$ as $\lambda\to\lambda_0$. We will achieve this objective by employing the general approach to time-periodic bifurcation presented in Section \ref{sec:bif}. The implementation of this approach to the case at hand will require a detailed study of the properties of several linearized problems associated to \eqref{04}, which will be established in Sections \ref{steady}, \ref{sec:spectrum} and \ref{sec:time_per}.
\setcounter{equation}{0} 
\section{Function Spaces and their Relevant Properties}\label{spaces}
We begin to recall some basic notation. By $\Omega$ we indicate a domain of $\real^3$,  complement of the closure of a  bounded domain $\Omega_0$  of class $C^2$. With the origin of coordinates in the interior of $\Omega_0$, we set  $B_r:=\{x\in\real^3:\,|x|<r\,,\ r>0\}$, $\Omega_R:=\Omega\cap B_R$,  $\Omega^R:=\Omega\backslash\bar{\Omega_R}$, $\Omega_{R_1,R_2}:=\Omega^{R_1}\backslash\bar{\Omega^{R_2}}$, with $R,R_1,R_2>R_*:={\rm diam}\,\Omega_0$, $R_2>R_1$. As customary, for $A$ a domain of $\real^3$,  $L^q=L^q(A)$ denotes the Lebesgue space with norm $\|\cdot\|_{q,A}$, and  $W^{m,2}=W^{m,2}(A)$, $m\in\nat$,  the Sobolev space with norm $\|\cdot\|_{m,2,A}$. By $(\,\ ,\,\ )_A$ we indicate the $L^2(A)$-scalar product.  Furthermore, $D^{m,q}=D^{m,q}(A)$ is the homogeneous Sobolev space with semi-norm $\sum_{|l|=m}\|D^lu\|_{q,A}$, whereas $D_0^{1,2}=D_0^{1,2}(A)$ is the completion of $C_0^\infty(A)$ in the norm $\|\nabla(\cdot)\|_{2,A}$. 
In all the above notation we shall typically omit the subscript ``$A$", unless confusion arises.
\par
We shall next introduce  certain function classes along with some of their most important  properties. 
Let
$$\ba{ll}\medskip
\calk=\mathcal K(\real^3):=\big\{\bfphi\in C_0^\infty(\bar{\real^3}):\ \mbox{in $\real^3$}\,;
\bfphi(x)=\hat{\bfphi}, \ \mbox{some $\hat{\bfphi}\in\real^3$}\, \ \mbox{in a neighborhood of $\Omega_0$}\big\}
\\ \medskip
\calc=\mathcal C(\real^3):=\{\bfphi\in\calk(\real^3):\ \Div\bfphi=0\ \mbox{in $\real^3$}\}\,,\\
\calc_0=\mathcal C_0(\Omega):=\{\bfphi\in\calc(\real^3): \hat{\bfphi}=\0\}\,.
\ea
$$
In $\calk$ we consider the following scalar product
\be
\langle \bfphi,\bfpsi\rangle:=
\varpi^{-1}\,\hat{\bfphi}_1\cdot\hat{\bfpsi}_1+(\bfphi,\bfpsi)\,,\ \ \bfphi,\bfpsi\in\calk\,,
\eeq{0.0}
and define the following spaces
\be\ba{ll}\medskip
\mathcal L^2=\call^2(\real^3):= \,\big\{\mbox{completion of $\calk(\real^3)$ in the norm induced by \eqref{0.0}}\big\}\,,\\ \medskip
\mathcal H=\calh(\real^3):=\,\big\{\mbox{completion of $\calc(\real^3)$ in the norm induced by \eqref{0.0}}\big\}\\ 
\calg=\calg(\real^3):=\big\{\bfh\in \call^2(\real^3): \, \mbox{there is $p\in D^{1,2}(\Omega)$ such that $\bfh=\nabla p$ in $\Omega$}\,,
\mbox{and}\ \bfh=-\varpi\int_{\partial\Omega}p\,\bfn   \ \mbox{in}\ \Omega_0\,\big\}\,.
\ea
\eeq{spazi}
\Bl The following characterizations hold
$$\ba{ll}\medskip
\call^2(\real^3)=\{\bfu\in L^2(\real^3): \ \bfu=\hat{\bfu}\ \mbox{in}\ \Omega_0, \ \mbox{for some $\hat{\bfu}\in \real^3$}\}\\
\calh(\real^3)=\{\bfu\in \call^2(\real^3): \ \Div\bfu=0\,\}\,.
\ea
$$
\EL{aria}
{\em Proof.} See  \cite[Theorem 3.1 and Lemma 3.2]{ALS}. \hfill$\square$\smallskip\par
We also have
\Bl The following orthogonal decomposition holds
$$
\call^2(\real^3)=\calh(\real^3)\oplus\calg(\real^3)\,.
$$
\EL{0.1}
{\em Proof.} A proof can be deduced from  \cite[Theorem 3.2]{ALS}. However, for completeness and since this result plays a major role in our analysis, we would like to reproduce it here. Let $\bfu\in \calh$ and $\bfh\in\calg$. Then,
$$
\langle\bfu,\bfh\rangle=\int_{\Omega}\bfu\cdot\nabla p-\int_{\partial\Omega}p\,\hat{\bfu}\cdot\bfn\,.
$$
Therefore, integrating by parts and using $\Div\bfu=0$ we deduce
$$
\langle\bfu,\bfh\rangle=-\int_{\Omega}p\,\Div\bfu+\int_{\partial\Omega}p\,\hat{\bfu}\cdot\bfn-\int_{\partial\Omega}p\,\hat{\bfu}\cdot\bfn=0
$$
which proves $\calh^{\perp}\supset \calg$. Conversely, let $\bfv\in\calh^{\perp}$, and assume  
\be
\varpi^{-1}\hat{\bfv}\cdot\hat{\bfu}+\int_\Omega\bfv\cdot\bfu=0\,,\ \ \mbox{for all $\bfu\in\calh$.}
\eeq{dip}
Since $\calc_0\subset \calh$, by picking $\bfu\in \calc_0$ from the preceding  we find, in particular,
$$
\int_\Omega\bfv\cdot\bfu=0\,,\ \ \mbox{for all $\bfu\in\calc_0$,}
$$ 
so that, by well-known results on the Helmholtz decomposition \cite[Lemma III.1.1]{Gab}, we infer $\bfv=\nabla p$ with $p\in D^{1,2}(\Omega)$. Replacing the latter into \eqref{dip} and integrating by parts, we get
$$
\left(\varpi^{-1}\hat{\bfv}+\int_{\partial\Omega}p\,\bfn\right)\cdot\hat{\bfu}=0\ \ \mbox{for all $\hat{\bfu}\in\real^3$}\,,
$$ 
from which we conclude that $\bfv\in\calg$, that is, $\calh^\perp\subset \calg$. The proof of the lemma is completed.
\par\hfill$\square$\par
We shall denote by $\mathscr P$ the self-adjoint, orthogonal projection of $\mathcal L^2(\real^3)$ onto $\calh(\real^3)$.\smallskip\par
We next define
$$\ba{ll}\medskip
\cald^{1,2}=\cald^{1,2}(\real^3):=\,\big\{\mbox{completion of $\calc(\real^3)$ in the norm $\|\mathbb D(\cdot)\|_2$}\big\}\,,\\
\cald_0^{1,2}=\cald_0^{1,2}(\Omega):=\,\big\{\mbox{completion of $\calc_0(\Omega)$ in the norm $\|\mathbb D(\cdot)\|_2$}\big\}
\ea
$$
and set
\be
Z^{2,2}:=W^{2,2}(\Omega)\cap \cald^{1,2}(\real^3)\,.
\eeq{zetaa}
Notice that, clearly, $\cald_0^{1,2}(\Omega)\subset\cald^{1,2}(\real^3)$. 
The basic properties of the space $\cald^{1,2}$ are collected in the next lemma, whose proof is given in \cite[Lemmas 9--11]{Gah}.
\Bl  Let $\tilde{\cald^{1,2}}$ denote either $\cald^{1,2}$ or $\cald_0^{1,2}$. Then, $\tilde{\cald^{1,2}}$ is a Hilbert space when equipped with the scalar product
$$
(\mathbb D(\bfu_1),\mathbb D(\bfu_2))\,,\ \ \bfu_i\in \tilde{\cald^{1,2}}\,, \, \ i=1,2\,.
$$
Moreover, we have the characterization:
\be
\tilde{\cald^{1,2}}=\big\{\bfu\in L^6(\real^3)\cap D^{1,2}(\real^3)\,;\ \Div\bfu=0\,;\,
\bfu=\hat{\bfu} \ \mbox{in $\Omega_0$} \big\}\,,
\eeq{1.7}
with some $\hat{\bfu}\in\real^3$, if $\tilde{\cald^{1,2}}\equiv\cald^{1,2}$, and $\hat{\bfu}=\0$, if\, $\tilde{\cald^{1,2}}\equiv\cald_0^{1,2}$.
Also, for each $\bfu\in\tilde{\cald^{1,2}}$, it holds
\be
\|\nabla\bfu\|_2=\sqrt{2}\|\mathbb D(\bfu)\|_2\,,
\eeq{1.8}
and
\be
\|\bfu\|_6\le \kappa_0\,\|\mathbb D(\bfu)\|_2\,,\ \ \bfu\in\tilde{\cald^{1,2}}\,,
\eeq{1.9}
for some numerical constant $\kappa_0>0$. Finally, 
there is another positive constant $\kappa_1$ such that
\be
|\hat{\bfu}|\le \kappa_1\,\|\mathbb D(\bfu)\|_2\,.
\eeq{1.10}
\EL{1.1}
Along with the spaces $\call^2,\calh$, and  $\cald^{1,2}$ defined above, we shall need also their ``restrictions" to the ball $B_R$. Precisely, we set
$$\ba{ll}\medskip
\call^2(B_R):= \{\bfphi\in L^2(B_R): \, \bfphi|_{\Omega_0}=\hat{\bfphi}\,\ \mbox{for some $\hat{\bfphi}\in\real^3$}\}\\ \medskip
\calh(B_R):=\{\bfphi\in \call^2(B_R):\, \Div\bfphi=0\,,\ \bfphi\cdot\bfn|_{\partial B_R}=0\}\\
\cald^{1,2}(B_R):= \{\bfphi\in W^{1,2}(B_R): \, \Div\bfphi=0\,,\ \bfphi|_{\Omega_0}=\hat{\bfphi}\,\ \mbox{for some $\hat{\bfphi}\in\real^3$}\,,\ \bfphi|_{\partial B_R}=\0 \}\,.
\ea
$$
Then $\calh(B_R)$ and $\cald^{1,2}(B_R)$ are Hilbert spaces with scalar products
$$
\varpi^{-1}\,\hat{\bfphi}_1\cdot\hat{\bfphi}_2+(\bfphi_1,\bfphi_2)_{\Omega_R}\,,\ \ \bfphi_i\in\calh(\Omega_R)\,;\ \ (\mathbb D(\bfpsi_1),\mathbb D(\bfpsi_2))\,,\ \bfpsi_i\in \cald^{1,2}(\Omega_R)\,,\ \ i=1,2.
$$
\Br The space $\cald^{1,2}(B_R)$ can be viewed as a subspace of $\cald^{1,2}(\real^3)$, by extending to 0 in $\real^3\backslash B_R$ its generic element. Therefore, all the properties mentioned in \lemmref{1.1} continue to hold for $\cald^{1,2}(B_R)$. 
\ER{2.1}
\par
We next introduce  certain subspaces of $\cald_0^{1,2}(\Omega)$ and investigate their relevant embedding properties. 
To this end, let $\cald_0^{-1,2}(\Omega)$ be the dual space of $\cald_0^{1,2}(\Omega)$, endowed with the norm
$$
|\bff|_{-1,2}=\Sup{\mbox{\footnotesize $\ba{c}\bfphi\in \calc_0(\Omega)\\ \|\nabla\bfphi\|_2=1\ea $}}\left|\left(\bff,\bfphi\right)\right|\,,
$$
and set
$$
\caly:=\cald_0^{-1,2}(\Omega)\cap \calh(\real^3)\,,\ \ Y(\Omega):=\cald_0^{-1,2}(\Omega)\cap L^2(\Omega)\,,
$$
with associated norms
$$
\|\bfg\|_{\caly}:=\|\bfg\|_{2,\Omega}+|\bfg|_{-1,2}+|\hat{\bfg}|\,,\ \ \ \ \|\bff\|_{Y}:=\|\bff\|_{2}+|\bff|_{-1,2}\,.
$$
We then define
$$
X(\Omega):=\{\bfu\in\cald^{1,2}_0(\Omega):\, \partial_1\bfu\in \cald_0^{-1,2}(\Omega)\}\,,\ \ X^2=X^2(\Omega):=\left\{\bfu\in X(\Omega): D^2\bfu\in L^2(\Omega)\right\}\,.
$$
It  can be proved \cite[Proposition 65]{GaCe} that   $X(\Omega)$ and $X^2(\Omega)$ are (reflexive, separable) Banach space when equipped with the norms\footnote{Keep in mind that, by \lemmref{1.1} the norms $\|\nabla(\cdot)\|_2$ and $\|\mathbb D(\cdot)\|_2$ are equivalent in $\cald^{1,2}_0$.}
$$
\|\bfu\|_X:=\|\nabla\bfu\|_2+|\partial_1\bfu|_{-1,2}\,,\ \ 
\|\bfu\|_{X^2}:=\|\bfu\|_{X}+\|D^2\bfu\|_{2}\,. 
$$ 
\Bl The following continuous embedding properties hold
\be\ba{ll}\medskip
X^2(\Omega)\subset W^{2,2}(\Omega_R)\,,\ \ \mbox{for all $R>R_*$}\,,\\
X^2(\Omega)\subset  L^\infty(\Omega)\cap D^{1,q}(\Omega)\,,\ \ \mbox{all $q\in[2,6]$}.
\ea
\eeq{X2e}\EL{Xemb}
{\em Proof.} By \lemmref{1.1}, the first property is obvious. From \cite[Theorem II.6.1(i)]{Gab} it follows that $X^2(\Omega)\subset D^{1,6}(\Omega)$ which, in turn, by \cite[Theorem II.9.1]{Gab} and simple interpolation allows us to deduce
also the second stated property.
\par\hfill$\square$\par
The  embedding property  for the space $X(\Omega)$ is given in the following \propref{3.1}.  A proof of this result was first furnished in \cite[Proposition 1.1]{GaFu}. Here, for completeness, we shall give a different  proof that, in fact, turns out to be  also more direct. We need a  preparatory result, showed in the next lemma.
\Bl Let $\bfphi\in \cald_0^{1,2}(\Omega)\cap W_0^{1,\frac65}(\Omega_R)$, all $R>R_*$, with  $\partial_1\bfphi\in L^{\frac65}(\Omega)$. Assume, further,
\be
\int_{\partial B_r}|\bfphi|^2\le c\,r^{-1}\,,\ \ \mbox{all $r>R_*$}\,.
\eeq{G_0}
Then, there is $\{\bfphi_k\}\subset \calc_0(\Omega)$ such that
\be
\lim_{k\to\infty}\|\partial_1(\bfphi-\bfphi_k)\|_{\frac65}=\lim_{k\to\infty}\|\nabla(\bfphi-\bfphi_k)\|_{2}=0\,.
\eeq{G_00}
\EL{pr} 
{\em Proof.} For all sufficiently large $R$, let $\psi_R=\psi_R(x)$ be the  ``cut-off" function such that  \cite[Theorem II.7.1]{Gab}
$$ 
\psi_R(x)=\left\{\ba{ll}\medskip 1\ &\ \mbox{if $|x|\le {\rm e}$}\\
\psi\left(\frac{\ln\ln|x|}{\ln\ln R}\right)\ &\ \mbox{if $|x|>{\rm e}$}\,,\ea\right.
$$
where $\psi=\psi(r)$, $r>0$, is a smooth, non-increasing function that is 1 if $r\le1/2$ and 0 if $r\ge1$. Thus,
\be
\psi_R(x)=\left\{\ba{ll}\medskip
1\ &\ \mbox{if $|x|\le {\rm e}^{\sqrt{\ln R}}$}\\
0\ &\ \mbox{if $|x|\ge  R$}\,,
\ea\right.
\eeq{G1}
and 
\be
|\nabla\psi_R(x)|\le \frac{c}{\ln\ln R}(|x|\ln|x|)^{-1}\,, \ \ \,x\in \hat{\Omega}_R:=\{x\in \Omega:\,\ {\rm e}^{\sqrt{\ln R}}< |x|< R\}\,.
\eeq{G2}
Consider the problem
\be
\Div\bfw_R=-\nabla\psi_R\cdot\bfphi \ \mbox{in $\hat{\Omega}_R$}\,;\ \ \|\nabla\bfw_R\|_{1,q}\le c\,\|\nabla\psi_R\cdot\bfphi\|_{q}\,,\ \ q=\mbox{$\frac65$},2\,.
\eeq{G3}
Using the properties of $\bfphi$ and $\psi_R$, we easily show
\be
\int_{\hat{\Omega}_R}\nabla\Psi_R\cdot\bfphi=0\,.
\eeq{G3_1}
With the change of variables $x_i\to y_i=x_i/R$, it follows $\hat{\Omega}_R\to \tilde{\Omega}_R:=\{y\in\real^3:\, R^{-1}{\rm e}^{\sqrt{\ln R}}<|y|<1\}$, $R>R_*$. Since ${\rm diam}\,(\tilde{\Omega}_R)\le 1$, and problem \eqref{G3}$_1$--\eqref{G3_1} remains invariant under the above change, from \cite[Theorem III.3.1 and Remark III.3.12]{Gab} we deduce that \eqref{G3}--\eqref{G3_1} has a solution with a constant $c$ independent of $R$.
We next write
$$
\bfphi=(\psi_R\bfphi+\bfw_R)+((1-\psi_R)\bfphi-\bfw_R):=\bfphi^{(1)}_R+\bfphi^{(2)}_R\,,
$$
and observe that, by virtue of \eqref{G1}, \eqref{G3}$_2$,
\be\ba{ll}\medskip
\|\nabla\bfphi_R^{(2)}\|_2\le c\,\left(\|\nabla\bfphi\|_{2,\Omega^{{\rm e}^{\sqrt{\ln R}}}}+\|\nabla\psi_R\cdot\bfphi\|_2\right)\\
\|\partial_1\bfphi_R^{(2)}\|_{\frac65}\le c\,\left(\|\partial_1\bfphi\|_{\frac65,\Omega^{{\rm e}^{\sqrt{\ln R}}}}+\|\nabla\psi_R\cdot\bfphi\|_{\frac65}\right)
\,.
\ea
\eeq{G4}
By \eqref{G2} and \eqref{G_0}, we get
\be
\|\nabla\psi_R\cdot\bfphi\|_2^2\le\Frac{c}{(\ln\ln R)^2}\Int{{\rm e}^{\sqrt{\ln R}}}R r^{-2}(\ln r)^{-2}\Int{\partial B_r}{}|\bfphi|^2{\rm d}r\le \Frac{c}{(\ln\ln R)^2}\Int{{\rm e}^{\sqrt{\ln R}}}R r^{-3}(\ln r)^{-2}{\rm d}r\,,
\eeq{G5}
and similarly, also with the help of H\"older inequality,
\be
\|\nabla\psi_R\cdot\bfphi\|_{\frac65}^{\frac65}\le\Frac{c}{(\ln\ln R)^{\frac65}}\Int{{\rm e}^{\sqrt{\ln R}}}R r^{-\frac65}(\ln r)^{\frac65}\Int{\partial B_r}{}|\bfphi|^{\frac65}{\rm d}r\le \Frac{c}{(\ln\ln R)^{\frac65}}\Int{{\rm e}^{\sqrt{\ln R}}}R r^{-1}(\ln r)^{-\frac65}{\rm d}r
\eeq{G6}
As a consequence, from \eqref{G4}--\eqref{G6} and the assumptions, it follows that for any given $\varepsilon>0$ there exists $R_0=R_0(\varepsilon)>0$ such that
\be
\|\partial_1\bfphi^{(2)}_{R_0}\|_{\frac65}+\|\nabla\bfphi^{(2)}_{R_0}\|_{2}<\half\varepsilon\,.
\eeq{G7}
We now observe that, by assumption and  \eqref{G5}, \eqref{G6} and \eqref{G3}, it easily follows that $\bfphi^{(1)}\in \cald_0^{1,2}(\Omega_{R_0})\cap W_0^{1,\frac65}(\Omega_{R_0})$. Therefore, extending $\bfphi^{(1)}_{R_0}$ to 0 outside $\Omega_{R_0}$, by \cite[Exercise III.6.2]{Gab} we find a sequence $\{\bfphi_k\}\subset \{\bfphi\in C_0^\infty(\Omega_{R_0}): \Div\bfphi=0\}\subset \mathcal C_0(\Omega)$ with the property that for every $\varepsilon>0$, there is $k_0=k_0(\varepsilon)$ such that
\be
\|\nabla(\bfphi^{(1)}_{R_0}-\bfphi_k)\|_{\frac65}+\|\nabla(\bfphi^{(1)}_{R_0}-\bfphi_k)\|_{2}<\half\varepsilon\,,\ \ \mbox{for all $k\ge k_0$}\,. 
\eeq{G8}
Thus, since $\bfphi-\bfphi_k=\bfphi^{(1)}_{R_0}-\bfphi_k+\bfphi^{(2)}_{R_0}$ the statement \eqref{G_00} follows from \eqref{G7} and \eqref{G8}.
\par\hfill$\square$\par
\Bp $X(\Omega)$ is embedded in $L^4(\Omega)$. Moreover, there is $c=c(\Omega)$ such that
\be
\|\bfu\|_4\le
c\,\|\bfu\|_{X}\,.
\eeq{3_2}
\EP{3.1}
{\em Proof.} 
For a given $\bff\in C_0^\infty(\Omega)$, consider the following problem 
\be\ba{cc}\medskip\left.\ba{ll}\medskip
\Delta\bfphi-\partial_1{\bfphi}=\nabla p+\bff\ \ \mbox{in $\Omega$}\,,\\
\Div\bfphi=0\ea\right\}\ \ \mbox{in $\Omega$}\,,\\
\bfphi=\0\ \ \mbox{at $\partial\Omega$}\,.
\ea
\eeq{1_21}
This problem has a (unique) solution such that 
\be\bfphi\in L^{s_1}(\real^3)\cap D^{1,s_2}(\real^3)\cap D^{2,s_3}(\real^3)\,,\ \ \partial_1\bfphi\in L^{s_3}(\Omega)\,,\ \ p\in L^{s_4}(\real^3)\cap D^{1,s_3}(\real^3)\,,
\eeq{1_22}
for all $s_1>2$, $s_2>4/3$, $s_3>1$, $s_4>3/2$,
which, in particular, satisfies the following estimate
\be
\|\nabla\bfphi\|_{2}\le c\|\bff\|_{4/3}\,,
\eeq{1_23}
with $c=c(\Omega,s_i))$; see \cite[Theorem VII.7.1]{Gab}. Moreover, from \cite[Theorem VII.6.2 and Exercise VII.3.1]{Gab} it also follows that
\be 
\int_{\partial B_r}|\bfphi|^2\le c\,r^{-1}\,,\ \ \mbox{all $r>R_*$}\,.
\eeq{ER1}
From \eqref{1_21}$_1$ we obtain
\be
(\bfu,\bff)=
(\bfu,\Delta\bfphi-\partial_1\bfphi- \nabla p)\,.
\eeq{1_24}
Notice that the term on the right-hand side is well defined, since, by \lemmref{1.1},  $\bfu\in L^6(\Omega)$ and $\bfphi, p$ satisfy
 \eqref{1_22}.
In particular, $\nabla p\in L^{\frac65}(\Omega)$,  which implies $(\bfu,\nabla p)=0$, whereas, by an easily justified integration by parts, one shows that $(\bfu,\Delta\bfphi)=-(\nabla\bfu,\nabla\bfphi)$. 
From \eqref{1_22} and \eqref{ER1}, we see that $\bfphi$ satisfies the assumptions of \lemmref{pr}, and 
therefore, with $\{\bfphi_k\}\subset \calc_0(\Omega)$ the sequence introduced in that lemma, from \eqref{1_24} we get
$$ 
(\bfu,\bff)=
-(\bfu,\partial_1\bfphi_k)-(\bfu,\partial_1(\bfphi-\bfphi_k))- (\nabla\bfu,\nabla\bfphi)=(\partial_1\bfu,\bfphi_k)-(\bfu,\partial_1(\bfphi-\bfphi_k))- (\nabla\bfu,\nabla\bfphi)\,.
$$
This relation, in turn, delivers
\be
|(\bfu,\bff)|\le |\partial_1\bfu|_{-1,2}\|\nabla\bfphi_k\|_2+\|\bfu\|_6\|\partial_1(\bfphi-\bfphi_k)\|_{\frac65}+\|\nabla\bfu\|_2\|\nabla\bfphi\|_2\,
,
\eeq{InIt}
so that passing to the limit $k\to\infty$ in \eqref{InIt} and employing \eqref{G_00}, \eqref{1_23} we conclude
$$
|(\bfu,\bff)|\le \left(|\partial_1\bfu|_{-1,2}+\|\nabla\bfu\|_2\right)\|\nabla\bfphi\|_2\le c\,\|\bfu\|_X\|\bff\|_{\frac43}\,.
$$
The latter, by the arbitrariness  of $\bff\in C_0^\infty(\Omega)$ and the density of $C_0^\infty(\Omega)$ in $L^{\frac43}(\Omega)$, completes the proof of the lemma. 
\par\hfill$\square$
\par
We  end this  section by introducing certain spaces of time-periodic functions. A function $u:\Omega\times \real\mapsto \real^3$ is 
{\em $2\pi$-periodic},  if $u(\cdot,t+2\pi)=u(\cdot\,t)$, for a.a. $t\in \real$,
and we set
$
{\bar u(\cdot)}:=\frac{1}{2\pi}\int_{0}^{2\pi}u(\cdot,t){\rm d}t\,,
$
whenever the integral is meaningful.
Let $B$ be a function space  with seminorm $\|\cdot\|_B$. By $L^2(0,2\pi;B)$ we denote the class of functions
$u:(0,2\pi)\rightarrow B$ such that 
$$
\|u\|_{L^2(B)}:= \big(\Int{0}{2\pi}\|u(t)\|_B^2 \big)^{\frac 12}<\infty
$$
Likewise, we put
$$
W^{1,2}(0,2\pi;B)=\Big\{u\in L^{2}(0,2\pi;B): \partial_tu\in L^{2}(0,T;2\pi)\Big\}\,.
$$
For simplicity, we write $L^2(B)$ for $L^2(0,2\pi;B)$, etc. Moreover, we define the Banach spaces
$$\ba{rl}\medskip
L^{q}_\sharp&\!\!\!\!:=\{\bfxi\in L^{q}(0,2\pi), \ \mbox{$\bfxi$ is $2\pi$-periodic with $\bar{\bfxi}=\0$}\}\,,\ \ q\in[1,\infty]\,,\\\medskip
W^{k}_\sharp&\!\!\!\!:=\{\bfxi\in L^{2}_\sharp(0,2\pi), \ d^l{\bfxi}/dt^l\in L^2(0,2\pi)\,,\ \ l=1,\ldots,k\}\,,\\\medskip
\mathcal L_\sharp^{2}&\!\!\!\!:=\{\bfu\in L^{2}(L^2(\Omega)); \ \mbox{$\bfu$ is $2\pi$-periodic, with $\bar{\bfu}=\0$}\}\,,
\\
\mathcal W_\sharp^{2}&\!\!\!\!:=\{\bfu\in W^{1,2}(L^2(\Omega))\cap L^2(W^{2,2}(\Omega)); \ \mbox{$\bfu$ is $2\pi$-periodic, with $\bar{\bfu}=\0$}\}\,,
\ea
$$
with associated norms
$$\ba{ll}\medskip
\|\bfxi\|_{L^q_\sharp}:=\|\bfxi\|_{L^q(0,2\pi)}\,,\ \ \ \|\bfxi\|_{W^k_\sharp}:=\|\bfxi\|_{W^{k,2}(0,2\pi)}\,,\\
\|\bfu\|_{\mathcal L^{2}_\sharp}:=\|\bfu\|_{L^2(L^2(\Omega))}\,,\ \ \
\|\bfu\|_{\mathcal W_\sharp^{2}}:=\|\bfu\|_{W^{1,2}(L^2(\Omega))}+\|\bfu\|_{L^2(W^{2,2}(\Omega))}\,.
\ea
$$
Finally, we define the Banach spaces
$$\ba{rl}\medskip
\textsf{{W}$_\sharp^{2}$}&\!\!\!\!:=\left\{\bfu\in L^2(Z^{2,2})\cap W^{1,2}(\calh):\, \mbox{$\bfu$ is $2\pi$-periodic,\,with\, $\bar{\bfu}|_{\Omega}=\bar{\hat{\bfu}}=\0$}\right\}\,,
\\
\textsf{{L}$_\sharp^{2}$}&\!\!\!\!:=\left\{\bfu\in L^2(\calh):\, \mbox{$\bfu$ is $2\pi$-periodic,\,with\, $\bar{\bfu}|_{\Omega}=\bar{\hat{\bfu}}=\0$}\right\}
\ea
$$
with corresponding norms
$$
\|\bfu\|_{\mbox{\scriptsize $\W$}}:=\|\partial_t\bfu\|_{L^2(\Omega)}+\|\bfu\|_{L^2(W^{2,2}(\Omega))}+\|\hat{\bfu}\|_{W_\sharp^1}\,,\ \ \
\|\bfu\|_{\mbox{\scriptsize $\M$}}:=\|\bfu\|_{L^2(L^2(\Omega))}+\|\hat{\bfu}\|_{L^2_\sharp}
\,.
$$
\setcounter{equation}{0} 
\section{Steady-States: Existence and Uniqueness }
\label{sec:exst}
Our first goal is to investigate  existence and uniqueness  properties of solutions to the steady-state problem \eqref{03}. \par 
We begin to state a general existence result  in a suitable function class, followed by a  corresponding uniqueness result. Both findings are, in fact, obtained as a corollary to classical results regarding steady-state Navier-Stokes problem in exterior domains.   
Precisely, we have the following theorem.
\Bt For any $\lambda>0$, problem \eqref{03} has at least one corresponding solution \be{\sf s}_0(\lambda):=(\bfu_0(\lambda),p_0(\lambda),\bfchi_0(\lambda))\in [L^q(\Omega)\cap D^{1,r}(\Omega)\cap D^{2,s}(\Omega)]\times[L^{\sigma}(\Omega)\cap D^{1,s}(\Omega)]\times\real^3\,,
\eeq{sfaco}
for all $q\in (2,\infty]$, $r\in(\frac43,\infty]$, $s\in (\frac32,\infty]$, and $\sigma\in (1,\infty)$. Moreover, set
\be\frac1{\lambda_1}:=-\max_{{\mbox{\footnotesize $\bfu\in\cald^{1,2}_0(\Omega)$}}}\frac{(\bfu\cdot\nabla\bfu_0,\bfu)}{\|\nabla\bfu\|_2^2}\,.
\eeq{sfaco1}
Then $1/\lambda_1$ exists and ${\sf s}_0(\lambda)$ is unique, provided $\lambda<\lambda_1$.\footnote{We suppose $\lambda_1>0$,  otherwise, as it is clear from the proof, no restriction is needed on $\lambda$.}
\ET{exi}
{\em Proof.} From \cite[Theorem X.6.4]{Gab} we know that for any $\lambda>0$  problem \eqref{03}$_{1-4}$ has  one corresponding  solution $(\bfu_0,p_0)$ in the class \eqref{sfaco}. We then set
\be
\bfchi_0:= -\frac{\varpi}{\omega_{\sf n}^2}\int_{\partial\Omega}\mathbb T(\bfu_0,p_0)\cdot\bfn\,,
\eeq{2.3_00}
which is well defined by standard trace theorems, and complete the proof of the existence. We now turn to the uniqueness proof. The existence of $1/\lambda_1$ follows by the summability properties of $\bfu_0$ given in \eqref{sfaco} and standard arguments about maxima of quadratic forms in exterior domains \cite{Gam}. Let $(\bfu+\bfu_0,p+p_0,\bfchi+\bfchi_0)$ be another solution to \eqref{03} in the class \eqref{sfaco} corresponding to the same $\lambda$. Then $(\bfu,p,\bfchi)$
satisfies the following equations
\be\ba{cc}\medskip\left.\ba{ll}\medskip
-\lambda\,(\partial_1\bfu-(\bfu_0+\bfu)\cdot\nabla\bfu-\bfu\cdot\nabla\bfu_0)=\Delta\bfu-\nabla p\\
\Div\bfu=0\ea\right\}\ \ \mbox{in $\Omega$}\,,\\ \medskip
\bfu(x)=\0\,, \ \mbox{ $x\in\partial\Omega$}\,,\\
\omega_{\sf n}^2\bfchi+\varpi\Int{\partial\Omega}{} \mathbb T(\bfu,p)\cdot\bfn=\0\,.
\ea
\eeq{03_111}
By dot-multiplying \eqref{03_111}$_1$ by $\bfu$, integrating by parts over $\Omega$ and using \eqref{sfaco}, \eqref{03_111}$_{2,3}$ and \eqref{sfaco1} we easily deduce
$$
\|\nabla\bfu\|_2^2=-\lambda(\bfu\cdot\nabla\bfu_0,\bfu)\le \frac{\lambda}{\lambda_1}\|\nabla\bfu\|_2^2\,,
$$
from which it follows $\bfu\equiv\0$ for $\lambda<\lambda_1$. From the latter and \eqref{03_111}$_1$, \eqref{sfaco}, \eqref{03_111}$_4$ we infer $\bfchi=\0$, which ends the proof of the theorem. 
\par\hfill$\square$
\setcounter{equation}{0} 
\section{Steady-States: Stability}
\label{sec:stab}
Our next task is to find sufficient conditions for the stability of solutions determined in \theoref{exi}, in a suitable class of ``perturbations." In this regard, let ${\sf s}_0(\lambda)$ be the steady-state solution given in \eqref{sfaco} and let $(\bfu,p,\bfchi)$ be a corresponding time-dependent perturbation. By \eqref{02} we then have that $(\bfu,p,\bfchi)$ must satisfy the following set of equations
\be\ba{cc}\medskip\left.\ba{ll}\medskip
\partial_t\bfu-\lambda\,[\partial_1\bfu-\bfu_0\cdot\nabla\bfu+(\dot{\bfchi}-\bfu)\cdot\nabla\bfu_0 -(\dot{\bfchi}-\bfu)\cdot\nabla\bfu]=\Delta\bfu-\nabla {p}\\
\Div\bfu=0\ea\right\}\ \ \mbox{in $\Omega\times(0,\infty)$}\,,\\ \medskip
\bfu(x,t)={\dot{\bfchi}}(t)\,, \ \mbox{ $(x,t)\in\partial\Omega\times(0,\infty)$}\,,\\ \medskip
\ddot{\bfchi}+\omega_{\sf n}^2\bfchi+\varpi\Int{\partial\Omega}{} \mathbb T(\bfu,{p})\cdot\bfn=\0\,, \ \ \mbox{in $(0,\infty)$}\\
\bfu(x,0)=\bfu^0\,,\ x\in\Omega\,,\ \ \bfchi(0)=\bfchi^{0}\,,\ \ \dot{\bfchi}(0)=\bfchi^1 
\,.
\ea
\eeq{04_111}
In order to study  problem \eqref{04_111}, we need some preliminary results.
\Bl Let either $\{A,D\}\equiv\{\real^3,\Omega\}$ or $\{A,D\}\equiv\{B_R,\Omega_R\}$,  and let $(\bfu,p)\in [\cald^{1,2}(A)\cap W^{2,2}(D)]\times D^{1,2}(D)$.
Then, there is a constant $C$ depending only on the regularity of $\Omega$ such that
$$
\|D^2\bfu\|_{2,D}+\|\nabla p\|_{2,D}\le C\,(\|\Div \mathbb T(\bfu,p)\|_{2,D}+\|\nabla\bfu\|_{2,D}+|\hat{\bfu}|)\,,
$$
where, we recall, $\hat{\bfu}=\bfu|_{\Omega_0}$.
\EL{stoces}
{\em Proof.} See \cite[Lemma 1]{Hey}\footnote{In that lemma, the domain is requested to be of class $C^3$. However, $C^2$ would suffice.\label{foo:hey}}.\par\hfill$\square$\par
\Bl Let $\bfu\in \cald^{1,2}(A)\cap W^{2,2}(D)$, with $A$ and $D$ as in \lemmref{stoces}, and  let $\bfv\in L^2(D)$. Then, for any $\varepsilon>0$ there exists a positive constant $C$ depending only on $\varepsilon$, $\bfu_0$ and the regularity of $\Omega$ such that
$$\ba{ll}\medskip
\left|\left(\partial_1\bfu-\bfu_0\cdot\nabla\bfu+(\hat{\bfu}-\bfu)\cdot\nabla\bfu_0 -(\hat{\bfu}-\bfu)\cdot\nabla\bfu,\bfv\right)_D\right|\\
\hspace*{3cm}\le C\,(\|\nabla\bfu\|_{2,D}^2 +\|\nabla\bfu\|_{2,D}^4+\|\nabla\bfu\|_{2,D}^6)+\varepsilon(\|D^2\bfu\|_{2,D}^2+\|\bfv\|_{2,D}^2)\,\ea
$$
\EL{prep}
{\em Proof.} Let us denote by $I_i$, $i=1,\ldots6$, in the order, the six terms in the scalar product. Taking  into account that, by \theoref{exi}, $\bfu_0\in L^\infty(\Omega)\cap D^{1,q}(\Omega)$, $q=2,3$, and using \eqref{1.9}, \eqref{1.10},  H\"older and Cauchy-Schwarz inequalities we readily get
$$
\sum_{i=1}^5\left|I_i\right|\le C\,( \|\nabla\bfu\|_{2,D}^2 +\|\nabla\bfu\|_{2,D}^4)+\half\varepsilon\|\bfv\|_{2,D}^2\,.
$$Moreover, again by H\"older inequality,  \cite[Lemma 1]{Hey}\footnote{See Footnote \ref{foo:hey}.},  \eqref{1.9} and \remref{2.1},
$$\ba{rl}\medskip
|I_6|\le \|\bfu\|_6\|\nabla\bfu\|_3\|\bfv\|_2&\!\!\!\le C\,\|\nabla\bfu\|_2\left(\|D^2\bfu\|_{2,D}^{\frac12}\|\nabla\bfu\|_{2,D}^{\frac12}+\|\nabla\bfu\|_{2,D}\right)\|\bfv\|_{2,D}\\
&\!\!\!\le C\,(\|\nabla\bfu\|_{2,D}^4+\|\nabla\bfu\|_{2,D}^6)+\varepsilon\,\|D^2\bfu\|_{2,D}^2+\half\varepsilon\,\|\bfv\|_{2,D}^2\,.
\ea
$$
The lemma is proved.
\par\hfill$\square$\par
\Bl 
Let $y\in C^1(0,\infty)$, $y(t)\ge 0$, satisfy
\be
\dot{y}(t)\le a+b[y(t)+y^\alpha(t)]\,,\ \ \alpha\ge 1,\ \ t>0\,,
\eeq{a}
where $a,b\in L^\infty(0,\infty)$.
Assume, also, $y\in L^1(0,\infty)$, and set
$$
{\sf a}:=\essup{t\in (0,\infty)}\,|a(t)|\,,\ {\sf b}:=\essup{t\in (0,\infty)}\,|b(t)|\,.
$$ 
Then, there exists a constant $\delta>0$, such that if
\be
y(0)\le \delta\,,\ \ \int_0^\infty y(s)\,ds\le\delta^2
\eeq{b}
it follows that:
$$\mbox{$y(t)< M\,\delta$ for all $t\in (0,\infty)\,,$ $M:=3\max\{1,2{\sf a},2{\sf b}\}$}\,.
$$
\EL{gronwa}
{\em Proof.} Let
$$
Y:=y^2\,,\ \beta:=(1+\alpha)/2\,.
$$
Multiplying both sides of \eqref{a} by $y$  we get
\be
\dot{Y}\le 2{\sf a}y+2{\sf b}[Y+Y^\beta]
\eeq{a1}
Contradicting (i) means that there there exists $t_0>0$ such that \be y(0)\le\delta\,,\,\ y(t)<M\,\delta\,, \ \ \mbox{for all $t\in (0,t_0)$\,,\, and  \,$y(t_0)=M\,\delta.$}\eeq{c} Now, integrating both sides of \eqref{a1} from 0 to $t_0$ and using the latter and \eqref{b}$_2$, we deduce, in particular
$$\ba{rl}\medskip
Y(t_0)&\!\!\!\le Y(0) +2{\sf a}\Int0{\infty}y(s)\,{\rm d}s+2{\sf b}\Int0{t_0} Y(s)\,{\rm d}s +2{\sf b}\Int0{t_0}Y^\beta(s)\,{\rm d}s\,\\
&\!\!\!\le \delta^2+2{\sf a}\delta^2+2{\sf b}M\delta^3
+2{\sf b}\delta^2(M\delta)^\alpha\le\mu \delta^2 \,(2+M\delta+(M\delta)^\alpha)\,, 
\ea 
$$
where $\mu=\max\{1,2{\sf a},2{\sf b}\}$. Therefore, choosing $\delta>0$ in such a way that
$$
2+M\delta+(M\delta)^\alpha<4\mu
$$
we deduce $y(t_0)\le\frac{2}3\,M\,\delta$, which contradicts \eqref{c}. \par\hfill$\square$\par
We are now in a position to prove the following stability result.
\Bt Suppose $\bfu_0$ is such that
\be\frac1{\lambda_2}:=-\sup_{\mbox{\footnotesize $\ \bfu\in\cald^{1,2}(\real^3)$}}\frac{((-\hat{\bfu}+\bfu)\cdot\nabla\bfu_0,\bfu)}{\|\nabla\bfu\|_2^2}\,,
\eeq{sfaco2}
exists finite,
where, we recall, $\hat{\bfu}=\bfu|_{\Omega_0}$, and assume $\lambda<\lambda_2$.\footnote{We suppose $\lambda_2>0$,  otherwise, as it is clear from the proof, no restriction is needed on $\lambda$.} Then, there is $\delta=\delta(\Omega,\lambda,\omega_{\sf n},\varpi)>0$ such that if
$$
\|\bfu^0\|_{1,2}+|\bfchi^0|+|\bfchi^1|\le\delta
$$
problem \eqref{04_111} has one and only one solution such that
\be\left.\ba{ll}\medskip
\bfu\in C([0,T]; \cald^{1,2}(\real^3))\cap L^2(0,T; W^{2,2}(\Omega))\cap W^{1,2}(0,T; \call^2(\real^3))\,,\\ p\in L^2(0,T;D^{1,2}(\Omega))\,,\ \ \bfchi\in W^{2,2}(0,T)\,, 
\ea\right.
\ \ \mbox{for all $T>0$.}
\eeq{cLaSS}
Moreover,
\be
\lim_{t\to\infty}\left(\|\bfu(t)\|_6+\|\nabla\bfu(t)\|_2+|\dot{\bfchi}(t)|+|\bfchi(t)|\right)=0
\eeq{asym}
\ET{5.1_01} 
{\em Proof.} 
To show the existence of solutions to \eqref{04_111}, 
we shall follow the arguments introduced and developed in \cite{GaSi1,GaSi2}. Let $\{\Omega_R\}$, $R\in\nat$, be an increasing sequence such that $\Omega=\cup_{R\in\nat}\Omega_R$ and, for each fixed $R$, consider the  problem  
\be\ba{cc}\medskip\left.\ba{ll}\medskip
\partial_t\bfu_R-\lambda\,[\partial_1\bfu_R-\bfu_0\cdot\nabla\bfu_R+({\bfsigma}_R-\bfu_R)\cdot\nabla\bfu_0 -({\bfsigma}_R-\bfu_R)\cdot\nabla\bfu_R]\\ \medskip
\hspace*{7.2cm}=\Div\mathbb T(\bfu_R,{p}_R)\\
\Div\bfu_R=0\ea\right\}\, \ \mbox{in $\Omega_R\times(0,\infty)$}\,,\\ \medskip
\bfu_R|_{\partial\Omega}={{\bfsigma}_R}(t)\,, \ \ \bfu_R|_{\partial B_R}=\0\,, \ \ \mbox{in $(0,\infty)$}\\ \medskip
\dot{\bfsigma}_R+\omega_{\sf n}^2\bfchi_R+\varpi\Int{\partial\Omega}{} \mathbb T(\bfu_R,{p}_R)\cdot\bfn=\0\,, \ \ \dot{\bfchi}_R=\bfsigma_R\ \ \mbox{in $(0,\infty)$}\\
\bfu_R(x,0)=\bfu^0\,,\ x\in\Omega_R\,,\ \ \bfchi_R(0)=\bfchi^{0}\,,\ \ {\bfsigma}_R(0)=\bfchi^1 
\,.
\ea
\eeq{04_11R}
Our approach to existence develops in two steps. In the first one, we show by the classical Galerkin method that  \eqref{04_11R} has a solution in the class \eqref{cLaSS}. This will be accomplished with the help of a suitable base, constituted by eigenvectors of a modified Stokes problem. This procedure will also lead to the proof of estimates for $(\bfu_R,p_R,\bfchi_R)$ with bounds that are independent of $R$. In this way, in the second step, we will pass to the limit $R\to\infty$ and show that the limit functions $(\bfu,p,\bfchi)$ solve the original problem \eqref{04_111} along with the asymptotic property \eqref{asym}. We begin to put \eqref{04_11R} in a ``weak" form. If we multiply \eqref{04_11R}$_1$ by $\bfpsi\in\cald^{1,2}(B_R)$, integrate by parts and use \eqref{04_11R}$_{2,3,4}$ we deduce
\be\ba{ll}\medskip
(\partial_t\bfu,\bfpsi)-\lambda\,\left(\partial_1\bfu-\bfu_0\cdot\nabla\bfu+({\bfsigma}-\bfu)\cdot\nabla\bfu_0 -({\bfsigma}-\bfu)\cdot\nabla\bfu,\bfpsi\right)\\ \medskip\hspace*{1.5cm}=(\mathbb D(\bfu),\mathbb D ({\bfpsi}))-\varpi^{-1}(\dot{\bfsigma}+\omega^2_{\sf n}\bfchi)\cdot{\hat{\bfpsi}}\,,\ \ \mbox{for all $\bfpsi\in\cald^{1,2}(B_R)$}\,,\\
\hspace*{1.15cm}\dot{\bfchi}=\bfsigma\,,
\ea
\eeq{wefo}
where, for simplicity, the subscript $``R$" has been omitted and $(\cdot,\cdot)\equiv(\cdot,\cdot)_{\Omega_R}$. By using  standard procedures \cite{Gah}, it is easily seen that if $(\bfu,p,\bfchi,\bfsigma)$ is  a smooth solution to \eqref{wefo}, then  it also satisfies \eqref{04_11R}$_{1-6}$.  
In \cite{GaSi1} it is shown that the problem  
\be
\begin{array}{l}
\left. 
\begin{array}{l}
\hspace{-0.2cm} \displaystyle  - \nabla \cdot \mathbb T(\bfpsi,\phi)=\lambda  \,\bfpsi \medskip \\ 
\hspace{-0.2cm} \displaystyle \Div   \bfpsi = 0 \medskip
\end{array}
\right\} \mbox{ in } \Omega_R \,,\\  \medskip
\displaystyle \bfpsi  = \hat{\bfpsi} \ \   \mbox{ at } \partial\Omega\,, \ \  
\displaystyle \bfpsi= 0\ \   \mbox{ at } \partial B_R \,, \\ 
\displaystyle \lambda\,  \hat{\bfpsi}=\varpi\int_{\partial\Omega} \mathbb T(\bfpsi,\phi)\cdot \bfn \,,
\end{array}
\eeq{eipr}
with the natural extension $\bfpsi(x)=\hat{\bfpsi}$ in $\Omega_0$, 
admits a denumerable number of positive eigenvalues $\{\lambda_{Ri}\}$ clustering at infinity, 
and  corresponding eigenfunctions $\{\bfpsi_{Ri}\} \subset {\mathcal D}^{1,2}(B_R)\cap W^{2,2}(\Omega_R)$    forming an orthonormal basis of ${\mathcal{H}}(B_R)$ that is also orthogonal in $\cald^{1,2}(\Omega_R)$. Also, the correspondent ``pressure" fields satisfy $\phi_{Ri}\in W^{1,2}(\Omega_R)$, $i\in\nat$. Thus, for each fixed $R\in\nat$, we look for ``approximated" solutions to \eqref{wefo} of the form
\be
\bfu_N(x,t)=\sum_{k=1}^Nc_{kN}(t)\bfpsi_{Rk}(x)\,,\ \ \bfsigma_N(t)=\sum_{k=1}^Nc_{kN}(t)\hat{\bfpsi}_{Rk}\,,\ \ \bfchi_N(t)\,,
\eeq{gal0}
where the vector functions $\bfc_N(t):=\{c_{1N}(t),\ldots c_{NN}(t)\}$ and $\bfchi_N(t)$ satisfy the following system of  equations  
\be
\ba{ll}\medskip
(\partial_t\bfu_N,\bfpsi_{Ri})-\lambda\,\left(\partial_1\bfu_N-\bfu_0\cdot\nabla\bfu_N+({\bfsigma}_N-\bfu_N)\cdot\nabla\bfu_0 -({\bfsigma}_N-\bfu_N)\cdot\nabla\bfu_N,\bfpsi_{Ri}\right)\\ \medskip\hspace*{1.5cm}=-(\mathbb D(\bfu_N),\mathbb D ({\bfpsi}_{Ri}))-\varpi^{-1}(\dot{\bfsigma}_N+\omega^2_{\sf n}\bfchi_N)\cdot{\hat{\bfpsi}}_{Ri}\,,\ \ i=1,\ldots,N\,,
\\
\hspace*{.83cm}\dot{\bfchi}_N=\bfsigma_N
\ea
\eeq{gal}
where $P_N$ is the orthogonal projector (in both $\calh(B_R)$ and $\cald^{1,2}(B_R)$) onto ${\rm span}\, \{\bfpsi_{R1},\ldots,\bfpsi_{RN}\}$ \cite[Lemma 3.5]{GaSi2}. Indeed, \eqref{gal} is a system of first order differential equations in normal form in the unknowns $\bfc_N,\bfchi_N$. To this end it suffices to observe that
\be
\langle\bfpsi_{Ri},\bfpsi_{Rj}\rangle:=(\bfpsi_{Ri},\bfpsi_{Rj})+\varpi^{-1}\hat{\bfpsi}_{Ri}\cdot\hat{\bfpsi}_{Rj}=\delta_{ij}\,,
\eeq{ortg} 
so that \eqref{gal} is equivalent to
\be\ba{ll}\medskip
\dot{c}_{iN}= F_i(\bfc_N,\bfchi_N)\,,\ \ i=1,\ldots,N\,,\ \ \dot{\bfchi}_N=\bfsigma_N\,,\\ \medskip
\hspace*{2mm}F_i:=\sum_{k=1}^Nc_{kN}\left[\lambda\big(\partial_1\bfpsi_{Rk}-\bfu_0\cdot\nabla\bfpsi_{Rk}+(\hat{\bfpsi}_{Rk}-\bfpsi_{Rk})\cdot\nabla\bfu_0,\bfpsi_{Ri}\big)+\big(\mathbb D(\bfpsi_{Rk}),\mathbb D(\bfpsi_{Ri})\big)\right]\\
\hspace{1.4cm}-\lambda\sum_{k,m=1}^Nc_{kN}c_{mN}\left((\hat{\bfpsi}_{Rk}-{\bfpsi}_{Rk})\cdot\nabla{\bfpsi}_{Rm},{\bfpsi}_{Ri}\right)-\frac{\omega_{\sf n}^2}{\varpi}\bfchi_N\cdot\hat{\bfpsi}_{Ri}\,,
\ea
\eeq{ODE}
which we equip with the following initial conditions:
\be
c_{iN}(0)=(\bfu^0,\bfpsi_{Ri})+\varpi^{-1}\bfchi^1\cdot\hat{\bfpsi}_{Ri}\,,\ \ \bfchi(0)=\bfchi^0\,.
\eeq{inco}
From \eqref{gal0} and \eqref{ortg}, it follows that
\be
\|\bfu_N(0)\|_{2,\Omega_R}+\varpi^{-1}|\bfsigma(0)|^2\le\|\bfu^0\|_{2,\Omega}^2+\varpi^{-1}|\bfchi^1|^2\,.
\eeq{incou}
Likewise, since
$$
2 (\mathbb D(\bfpsi_{Ri}),\mathbb D(\bfpsi_{Rj}))= \lambda_{Ri} \big[ \varpi^{-1} \hat{\bfpsi}_{Ri} \cdot \hat{\bfpsi}_{Rj} +(\bfpsi_{Ri},\bfpsi_{Rj})\big] = \lambda_{Ri} \delta_{ij} 
$$
we have
\begin{eqnarray*}
\mathbb D(\bfu_N(0))&= &\sum_{j=1}^{N}c_{jN}(0)\mathbb D(\bfpsi_{Rj})= 2 \sum_{j=1}^{k}  \frac{1}{\lambda_{Rj}} (\mathbb D(\bfu^0), \mathbb D(\bfpsi_{Rj}))_{\Omega_R} \mathbb  D(\bfpsi_{Rj})\\&&= \sum_{j=1}^{N}   \frac{(\mathbb D(\bfu^0), \mathbb D(\bfpsi_{Rj})_{\Omega_R}}{\| \mathbb D(\bfpsi_{Rj}) \|_{2,\Omega_R}^{2}} \,   \mathbb D(\bfpsi_{Rj})\,,
\end{eqnarray*}
and therefore
\be
\|\mathbb D( \bfu_N(0))\|_{2,\Omega_R}\leq \|\mathbb D(\bfu^{0})\|_{2,\Omega}. 
\eeq{incod}
We shall now derive three basic ``energy estimates."
Multiplying both sides of \eqref{gal}$_1$ by $c_{iN}$, summing over $i$,  integrating by parts over $\Omega_R$ and using \eqref{gal0} we get
$$
\half\ode{}t\left[\|\bfu_N\|_2^2+\varpi^{-1}(|\bfsigma_N|^2+\omega_{\sf n}^2|\bfchi_N|^2)\right]=-2\|\mathbb D(\bfu_N)\|_2^2+\lambda\,((\bfsigma_N-\bfu_N)\cdot\nabla\bfu_0,\bfu_N)\,,
$$
which, by \eqref{04_111}, \eqref{1.8}, \remref{2.1} and the assumption, in turn furnishes
\be
\half\ode{}t\left[\|\bfu_N\|_2^2+\varpi^{-1}(|\bfsigma_N|^2+\omega_{\sf n}^2|\bfchi_N|^2)\right]\le-\gamma \,\|\nabla\bfu_N\|_2^2\,,\ \ \gamma:=1-\lambda/\lambda_2>0\,.
\eeq{gal1}
We next multiply both sides of \eqref{gal}$_1$ by $\dot{c}_{iN}$, sum over $i$,  and integrate by parts over $\Omega_R$ as necessary. Taking again into account \eqref{gal0}, \eqref{1.8} and \remref{2.1}, we show
\be\ba{ll}\medskip
\ode{}t\|\nabla\bfu_N\|_2^2+\|\partial_t\bfu\|_2^2+\varpi^{-1}|\dot{\bfsigma}_N|^2\\=\lambda\,\left(\partial_1\bfu_N-\bfu_0\cdot\nabla\bfu_N+({\bfsigma}_N-\bfu_N)\cdot\nabla\bfu_0 -({\bfsigma}_N-\bfu_N)\cdot\nabla\bfu_N,\partial_t\bfu_N\right)-\frac{\omega^2_{\sf n}}{\varpi}\bfchi_N\cdot\dot{\bfsigma}_N\,.
\ea
\eeq{gal2}
Finally, we multiply both sides of \eqref{gal1} by $\lambda_{Ri}c_{iN}$ and sum over $i$. Integrating by parts over $\Omega_R$ and employing \eqref{eipr} and, one more time, \eqref{04_111}, \eqref{1.8}, and  \remref{2.1} we show
\be\ba{ll}\medskip
\half\ode{}t\|\nabla \bfu_N\|_2^2+\|\Div\mathbb T(\bfu_N,p_N)\|_2^2+\varpi\,|\bfS_N|^2\\ 
\hspace*{.8cm}=\lambda\,\left(\partial_1\bfu_N-\bfu_0\cdot\nabla\bfu_N+({\bfsigma}_N-\bfu_N)\cdot\nabla\bfu_0 -({\bfsigma}_N-\bfu_N)\cdot\nabla\bfu_N,\Div\mathbb T(\bfu_N,p_N)\right)
+{\omega^2_{\sf n}}\bfchi_N\cdot\bfS_N
\ea
\eeq{gal3}
where
$$
\bfS_N:=\int_{\partial\Omega}\mathbb T(\bfu_N,p_N)\cdot\bfn\,,\ \ p_N:=\sum_{k=1}^Nc_{kN}\phi_{Rk}\,.
$$
We shall now derive a number of estimates for the approximated solutions,
paying attention that the constant involved are independent of $N$ and $R$. Such a generic constant will be denote by  $C$, which can thus depend, at most, on $\Omega$, $\bfu_0$ and the physical constant involved in the problem. Moreover, without specification, its value may change from a line to the next one (e.g. $2C\le C$). From \eqref{gal1}, \eqref{incou} and \eqref{1.10} we get
\be\ba{ll}\medskip
\Sup{t\in(0,\infty)}\left[\|\bfu_N(t)\|_2^2+\varpi^{-1}(|\bfsigma_N(t)|^2+\omega_{\sf n}^2|\bfchi_N(t)|^2)\right]+\gamma\Int0\infty(2\kappa_1^{-2}|\bfsigma_N(s)|^2+\|\nabla\bfu_N(s)\|_2^2)\,{\rm d}s\\
\hspace*{5cm}\le\|\bfu^0\|_2^2+\varpi^{-1}(|\bfchi^1|^2+\omega_{\sf n}^2|\bfchi^0|^2)\,.
\ea
\eeq{est1}
Such an estimate implies, in particular, that the initial-value problem \eqref{ODE}--\eqref{inco} has a (unique) solution in the whole interval $(0,\infty)$.
Moreover, from \eqref{gal2}, Cauchy-Schwarz inequality and \lemmref{prep} with $\bfv\equiv\partial_t\bfu_N$ and $\varepsilon\equiv\varepsilon_1<\half$, we infer
\be\ba{rl}\medskip
\ode{}t\|\nabla\bfu_N\|_2^2+\half\|\partial_t\bfu_N\|_2^2+|\dot{\bfsigma}_N|^2 &\!\!\!\le C\left(\|\nabla\bfu_N\|_2^2+\|\nabla\bfu_N\|_2^4+\|\nabla\bfu_N\|_2^6+|\bfchi_N|^2\right)+\varepsilon_1\,\|D^2\bfu_N\|_2^2\\
&\!\!\!\le C\left(\|\nabla\bfu_N\|_2^2+\|\nabla\bfu_N\|_2^6+|\bfchi_N|^2\right)+\varepsilon_1\,\|D^2\bfu_N\|_2^2
\ea
\eeq{est2}
Likewise, employing \lemmref{prep}, this time with $\bfv\equiv\Div \mathbb T(\bfu_N,p_N)$, from \eqref{gal3} we obtain
$$
\half \ode{}t\|\nabla\bfu_N\|_2^2+\|\Div\mathbb T(\bfu_N,p_N)\|_2^2+\varpi\,|\bfS_N|^2 
\le C\left(\|\nabla\bfu_N\|_2^2+\|\nabla\bfu_N\|_2^6+|\bfchi_N|^2\right)+\varepsilon_2\,\|D^2\bfu_N\|_2^2
$$
which, in turn, combined with \lemmref{stoces} implies, by taking $\varepsilon_2$ small enough
\be
\half\ode{}t\|\nabla\bfu_N\|_2^2+C\,\|D^2\bfu_N\|_2^2+\varpi\,|\bfS_N|^2 
\le C\left(\|\nabla\bfu_N\|_2^2+\|\nabla\bfu_N\|_2^6+|\bfchi_N|^2\right)\,.
\eeq{est3}
Summing side-by-side \eqref{est2} and \eqref{est3} by choosing $\varepsilon_1$ sufficiently small we deduce 
\be
\ode{}t\|\nabla\bfu_N\|_2^2+C\,(\|D^2\bfu_N\|_2^2+\|\partial_t\bfu_N\|_2^2+|\dot{\bfsigma}_N|^2+\varpi\,|\bfS_N|^2) 
\le C\left(\|\nabla\bfu_N\|_2^2+\|\nabla\bfu_N\|_2^6+|\bfchi_N|^2\right)\,,
\eeq{est4}
which furnishes,
in particular,
\be
\ode{}t\|\nabla\bfu_N\|_2^2
\le C\left(\|\nabla\bfu_N\|_2^2+\|\nabla\bfu_N\|_2^6+|\bfchi_N|^2\right)\,,
\eeq{est5}
We wish to apply \lemmref{gronwa} to \eqref{est5}, with $y\equiv \|\nabla\bfu_N\|_2^2$, $a\equiv C\,|\bfchi_N|^2$ and $b\equiv C$. By \eqref{est1}, it follows that both $a$ and $b$ satisfy the assumption of the lemma. Furthermore, if we choose initial data such that
$$
\|\nabla\bfu^0\|_2^2\le \delta\,,\ \
\|\bfu^0\|_2^2+\varpi^{-1}(|\bfchi^1|^2+\omega_{\sf n}^2|\bfchi^0|^2)\le\gamma^{-1}\delta^2
$$ 
from  \eqref{incod} and \eqref{est1} it follows that also the assumption \eqref{b} is met. As a result, \lemmref{gronwa} entails the existence of a positive constant $C_0$ independent of $N$ and $R$ such that
\be
\sup_{t\in(0,\infty)}\|\nabla\bfu_N(t)\|_2\le C_0\,.
\eeq{fiest}
Employing \eqref{fiest} in \eqref{est4} and keeping in mind \eqref{est1} we conclude
\be
\int_0^T\left(\|D^2\bfu_N(s)\|_2^2+\|\partial_t\bfu_N(s)\|_2^2\right)\,{\rm d}s\le C_1\,T\,,\ \ \mbox{for all $T>0$}\,,
\eeq{scest}
with $C_1$ another positive constant independent of $N$ and $R$. Thanks to \eqref{est1} and \eqref{scest}, we can now use a standard argument (e.g., \cite{GaSi1}) to prove the existence of a subsequence $\{(\bfu_{N_k},\bfchi_{N_k},\bfsigma_{N_k})\}$ converging in suitable topology to some $(\bfu_{R},\bfchi_{R},\bfsigma_{R})$ in the class \eqref{cLaSS} (with $\Omega\equiv\Omega_R$ and $\real^3\equiv B_R$)  and satisfying \eqref{wefo}. Since, clearly, $(\bfu_{R},\bfchi_{R},\bfsigma_{R})$ continue to obey the bounds \eqref{est1} and \eqref{scest}, we can similarly select a subsequence $(\bfu_{R_m},\bfchi_{R_m},\bfsigma_{R_m})$ converging (again, in suitable topology) to a certain $(\bfu,\bfchi,\bfsigma)$ that is in the class \eqref{cLaSS} and obeys \eqref{est1}, \eqref{scest}, and \eqref{04_111} for a.a. $x\in \Omega$ and $t\in(0,\infty)$. The demonstration of this convergence is rather typical and we will omit it, referring to \cite[Step 3 at p. 141]{GaSi1} for details. Thus, the proof of existence is completed. 
We now pass to show the validity of \eqref{asym}. In this regard, we begin to observe that the solution just constructed satisfies, in particular,
\be
\Sup{t\in(0,\infty)}\left(\|\bfu(t)\|_2+\|\nabla\bfu(t)\|_2+|\bfchi(t)|+|\dot{\bfchi}(t)|\right)+\Int0\infty(|\dot{\bfchi}(s)|^2+\|\nabla\bfu(s)\|_2^2)\,{\rm d}s\le K\,,
\eeq{asy1}
where $K>0$ is a constant depending only on the data. By dot-multiplying both sides of \eqref{04_111}$_1$ by $\partial_t\bfu$ and proceeding as in the proof of \eqref{est2} we obtain
$$\ba{ll}\medskip
\ode{}t(\|\nabla\bfu\|_2^2+\frac{\omega^2_{\sf n}}{\varpi}\bfchi\cdot\dot{\bfchi})+\|\partial_t\bfu\|_2^2+\varpi^{-1}|\ddot{\bfchi}|^2\\\hspace*{2cm}=\lambda\,\left(\partial_1\bfu-\bfu_0\cdot\nabla\bfu+(\dot{\bfchi}-\bfu)\cdot\nabla\bfu_0 -({\bfsigma}-\bfu)\cdot\nabla\bfu,\partial_t\bfu\right)+\mbox{$\frac{\omega^2_{\sf n}}{\varpi}$}|\dot{\bfchi}|^2\,.
\ea
$$
We now use, on the right-hand side of this relation, \lemmref{prep} with $\bfv\equiv\partial_t\bfu$, $\varepsilon\equiv\varepsilon_1<\half$,  along with the uniform bound on $\|\nabla\bfu\|_2$ in \eqref{asy1} to get
\be 
\ode{}t(\|\nabla\bfu\|_2^2+\mbox{$\frac{\omega^2_{\sf n}}{\varpi}$}\bfchi\cdot\dot{\bfchi})+\half\|\partial_t\bfu\|_2^2+\varpi^{-1}|\ddot{\bfchi}|^2\le C\,\|\nabla\bfu\|_2^2+\varepsilon_1\,\|D^2\bfu\|_2^2+\mbox{$\frac{\omega^2_{\sf n}}{\varpi}$}|\dot{\bfchi}|^2
\eeq{asy2}
Finally, we test both sides of \eqref{04_111}$_1$ by $\Div\mathbb T(\bfu,p)$  and apply Cauchy-Schwarz inequality to deduce
$$
\|\Div\mathbb T(\bfu,p)\|_2^2\le 
-2\lambda\,\left(\partial_1\bfu-\bfu_0\cdot\nabla\bfu+(\dot{\bfchi}-\bfu)\cdot\nabla\bfu_0 -(\dot{\bfchi}-\bfu)\cdot\nabla\bfu,\Div\mathbb T(\bfu,p)\right)+\|\partial_t\bfu\|_2^2\,.
$$
Employing in this inequality \lemmref{prep} with $\bfv\equiv\Div\mathbb T(\bfu,p)$ along with the bound \eqref{asy1} on $\|\nabla\bfu\|_2$,  we infer
$$
\|\Div\mathbb T(\bfu,p)\|_2^2\le C\,\|\nabla\bfu\|_2^2 +\varepsilon_2\,\|D^2\bfu\|_2^2+\|\partial_t\bfu\|_2^2\,,
$$
which, in turn, with the help of \lemmref{stoces} and by selecting $\varepsilon_2$ small enough, entails
\be
\|D^2\bfu\|_2^2+\|\nabla p\|_2^2\le C\,(\|\nabla\bfu\|_2^2 +\|\partial_t\bfu\|_2^2++|\dot{\bfchi}|^2)\,.
\eeq{asy3}
Next, we utilize \eqref{asy3} on the right-hand side of \eqref{asy2} and pick $\varepsilon_2$ suitably, which enables us to find
\be  
\ode{}t(\|\nabla\bfu\|_2^2+\mbox{$\frac{\omega^2_{\sf n}}{\varpi}$}\bfchi\cdot\dot{\bfchi})+\mbox{$\frac14$}\|\partial_t\bfu\|_2^2+\varpi^{-1}|\ddot{\bfchi}|^2\le C\,(\|\nabla\bfu\|_2^2+|\dot{\bfchi}|^2)\,.
\eeq{asy4}
Integrating over time both sides of \eqref{asy4} and taking into account \eqref{asy1} it follows  that
\be
\partial_t\bfu\in L^2(0,\infty; L^2(\Omega))\,,\ \ \ddot{\bfchi}\in L^2(0,\infty)\,,
\eeq{asy5}
which once replaced in \eqref{asy4}, again with the help of \eqref{asy1},  furnishes
\be
D^2\bfu,\nabla p\in L^2(0,\infty; L^2(\Omega))\,.
\eeq{asy6}
By possibly adding a suitable function of time to $p$, we may get \cite[Theorem II.6.1]{Gab}
\be
p\in L^6(\Omega)\ \ \mbox{and}\,\ \|p(t)\|_6\le C\,\|\nabla p(t)\|_2\,,\ \ \mbox{a.a. $t>0$\,.} 
\eeq{p_p}
On the other hand,
 from \eqref{04_111}$_4$ and standard trace theorems, we have
$$
|\bfchi(t)|^2\le \half\,\left(|\ddot{\bfchi}(t)|^2+\left|\int_{\partial\Omega}\mathbb T(\bfu,p)\cdot\bfn\right|^2\right)\le C\left(|\ddot{\bfchi}(t)|^2+\|\nabla\bfu\|_{2,2,\Omega_\rho}^2+\|p\|^2_{1,2,\Omega_\rho}\right)\,,
$$
for some fixed $\rho$, which, in vew of \eqref{asy5}--\eqref{p_p} allows us to conclude.
\be
\bfchi\in L^2(0,\infty)\,.
\eeq{asy7}
Combining \eqref{asy1}, \eqref{asy5} and \eqref{asy7} we get at once
\be
\lim_{t\to\infty}\left(|\bfchi(t)|+|\dot{\bfchi}(t)|\right)=0\,.
\eeq{lim1}
From \eqref{asy1} it follows that there exists at least one unbounded sequence $\{t_n\}\in (0,\infty)$ such that
\be
\lim_{n\to\infty}\|\nabla\bfu(t_n)\|_2=0\,.
\eeq{lim2}
Thus, integrating both sides of \eqref{asy4} between $t_n$ and $t>t_n$ we infer, in particular
$$
\|\nabla\bfu(t)\|_2^2\le \|\nabla\bfu(t_n)\|_2^2 +C\,\left(|\bfchi(t)|\,|\dot{\bfchi}(t)|+|\bfchi(t_n)|\,|\dot{\bfchi}(t_n)|+\int_{t_n}^\infty(\|\nabla\bfu(s)\|_2^2+|\dot{\bfchi}(s)|^2)\,{\rm d}s\right)\,, 
$$
which by \eqref{lim1} and \eqref{lim2} entails
$$
\lim_{t\to\infty}\|\nabla\bfu(t)\|_2=0\,.
$$
The latter and \eqref{1.9} complete the proof of \eqref{asym}.
The proof of uniqueness is quite standard and we will only sketch it here. Let $(\bfu_i,p_i,\bfchi_i)$, $i=1,2$, be two solutions to \eqref{04_111} corresponding to the same initial data, and set $\bfu:=\bfu_1-\bfu_2$, $p=p_1-p_2$, $\bfchi=\bfchi_1-\bfchi_2$. We thus have
\be\ba{cc}\medskip\left.\ba{ll}\medskip
\partial_t\bfu-\lambda\,[\partial_1\bfu-\bfu_0\cdot\nabla\bfu+(\dot{\bfchi}-\bfu)\cdot\nabla\bfu_0 -(\dot{\bfchi}_1-\bfu_1)\cdot\nabla\bfu-(\dot{\bfchi}-\bfu)\cdot\nabla\bfu_2]\\ \medskip\hspace*{7.5cm}=\Delta\bfu-\nabla {p}\\
\Div\bfu=0\ea\right\}\ \ \mbox{in $\Omega\times(0,\infty)$}\,,\\ \medskip
\bfu(x,t)={\dot{\bfchi}}(t)\,, \ \mbox{ $(x,t)\in\partial\Omega\times(0,\infty)$}\,,\\ \medskip
\ddot{\bfchi}+\omega_{\sf n}^2\bfchi+\varpi\Int{\partial\Omega}{} \mathbb T(\bfu,{p})\cdot\bfn=\0\,, \ \ \mbox{in $(0,\infty)$}\\
\bfu(x,0)=\0\,,\ x\in\Omega\,,\ \ \bfchi(0)=\0\,,\ \ \dot{\bfchi}(0)=\0 
\,.
\ea
\eeq{04_11U}
We dot-multiply \eqref{04_11U}$_1$ by $\bfu$, integrate by parts over $\Omega$ and use \eqref{04_11U}$_{2-4}$ to get
\be
\half\ode{}t\left[\|\bfu\|_2^2+\varpi^{-1}(|\dot{\bfchi}|^2+\omega_{\sf n}^2|\bfchi|^2)\right]+\gamma \|\nabla\bfu\|_2^2\le \lambda\,((\dot{\bfchi}-\bfu)\cdot\nabla\bfu_2,\bfu)\,,
\eeq{un1}
where $\gamma$ is defined in \eqref{est1}. From \eqref{1.10} and Cauchy--Schwarz inequality we get
$$
|(\bfchi\cdot\nabla\bfu_2,\bfu)|\le \half\gamma\|\nabla\bfu\|_2^2+c\,\|\nabla\bfu_2\|_2^2\|\bfu\|_2^2\,,
$$
whereas from \eqref{1.9}, H\"older,  Sobolev and Cauchy--Schwarz inequalities,
$$
|(\bfu\cdot\nabla\bfu_2,\bfu)|\le \|\bfu\|_6\|\nabla\bfu_2\|_3\|\bfu\|_2\le \half\gamma \|\nabla\bfu\|_2^2+c\,\|\bfu_2\|^2_{2,2}\|\bfu\|_2^2\,.
$$
Replacing the last two displayed relations in \eqref{un1} we thus conclude
\be
\ode Et\le c\, g(t)\, E(t)
\eeq{EN}
where $g:=\|\bfu(t)\|_{2,2}^2$, $E:=
\|\bfu\|_2^2+\varpi^{-1}(|\dot{\bfchi}|^2+\omega_{\sf n}^2|\bfchi|^2)$. Since  $\bfu_2$ is in the class \eqref{cLaSS}, we have $g\in L^1(0,T)$, for all $T>0$ and also, by assumption, $E(0)=0$. Uniqueness then follows by using Gronwall's lemma to \eqref{EN}. The proof of the theorem is completed.
\par\hfill$\square$\par
\setcounter{equation}{0}
\section{On the Absence of Oscillatory Regimes}\label{sec:nobif}
We begin with a simple observation concerning the two quantities defined in \eqref{sfaco1} and \eqref{sfaco2}. Precisely, since $\cald_0^{1,2}(\Omega)\subset \cald^{1,2}(\real^3)$, it follows that $\lambda_1^{-1}\le \lambda_2^{-1}$. This, by \theoref{exi} and \theoref{5.1_01}, implies that if $\lambda<\lambda_2$, the corresponding steady-state solution determined in \theoref{exi} is unique {\em and} stable. The objective of this section is to show, in addition, that,  if $\lambda<\lambda_2$ no oscillatory motion can stem out of the steady-state branch in a suitable function class of solutions $\calc$, which means to say that time-periodic bifurcation may occur only at some $\lambda_0>\lambda_2$. More precisely, let ${\sf s}_0(\lambda)$ be the steady-state solution given in \eqref{sfaco}. A generic $T$-periodic solution to \eqref{02} can then always be written as 
$$
\bfu(x,t)+\bfu_0(x)\,,\ \ p(x,t)+p_0(x)\,,\ \ \bfchi(t)+\bfchi_0\,,
$$
where $(\bfu,p,\bfchi)$, after the scaling  $\tau=t/T$, is a $2\pi$-periodic solution to the following equations
\be\ba{cc}\medskip\left.\ba{ll}\medskip
\frac1T\partial_\tau\bfu-\lambda\,[\partial_1\bfu-\bfu_0\cdot\nabla\bfu+(\dot{\bfchi}-\bfu)\cdot\nabla\bfu_0 -(\dot{\bfchi}-\bfu)\cdot\nabla\bfu]=\Delta\bfu-\nabla {p}\\
\Div\bfu=0\ea\right\}\ \ \mbox{in $\Omega\times(0,2\pi)$}\,,\\ \medskip
\bfu(x,t)={\dot{\bfchi}}(t)\,, \ \mbox{ $(x,t)\in\partial\Omega\times(0,2\pi)$}\,,\\ \medskip
\ddot{\bfchi}+\omega_{\sf n}^2\bfchi+\varpi\Int{\partial\Omega}{} \mathbb T(\bfu,{p})\cdot\bfn=\0\,, \ \ \mbox{in $(0,2\pi)$} 
\,.
\ea
\eeq{per0}
We now introduce the  class
$$
\calc:=\left\{(\bfu=\bar{\bfu}+\bfw,p=\bar{p}+{\sf p},\bfchi=\bar{\bfchi}+\bfxi):\, \bar{\bfu}\in X(\Omega)\,,  \bfw\in\W;\,\ \bar{p}\in W^{1,2}(\Omega)\,, {\sf p}\in L^2(D^{1,2});\,\ \bfxi\in W^2_\sharp\right\}\,,
$$
which constitutes  the functional framework  where, later on, we shall prove the occurrence of time-periodic bifurcation.
\par
The following result holds.
\Bt
Let $(\bfu,p,\bfchi)\in\calc$. Then, if $\lambda<\lambda_2$, necessarily $(\bfu,p,\bfchi)\equiv(\0,0,\0)$.
\ET{noper}
{\em Proof.} From \eqref{per0} we  get (with $\zeta:=\frac1T$)
\be\ba{cc}\medskip\left.\ba{ll}\medskip
-\lambda(\partial_1\bar{\bfu}-\bfu_0\cdot\nabla\bar{\bfu}-\bar{\bfu}\cdot\nabla\bfu_0)-\Delta\bar{\bfu}+\nabla\bar{p}=-\lambda(\bar{\bfu}\cdot\nabla\bar{\bfu}+\bar{(\zeta\dot{\bfxi}-\bfw)\cdot\nabla\bfw})\\
\Div\bar{\bfu}=0\ea\right\}\ \ \mbox{in $\Omega$\,,}\\ \medskip
\bar{\bfu}=\0\ \ \mbox{at $\partial\Omega$}\,,\\
\omega_{\sf n}^2\bar{\bfchi}+\varpi\Int{\partial\Omega}{}\mathbb T(\bar{\bfu},\bar{p})\cdot\bfn=\0\,,
\ea
\eeq{per1}
and
\be\ba{cc}\medskip\left.\ba{ll}\medskip
\zeta\partial_\tau\bfw-\lambda(\partial_1\bfw-\bfu_0\cdot\nabla\bfw+(\zeta\dot{\bfxi}-\bfw)\cdot\nabla\bfu_0)-\Delta\bfw+\nabla{\sf p}=\\ \medskip
\hspace*{3cm}\lambda[\bar{(\zeta\dot{\bfxi}-\bfw)\cdot\nabla\bfw}+{(\zeta\dot{\bfxi}-\bfw)\cdot\nabla\bfw}-\bar{\bfu}\cdot\nabla\bfw+(\zeta{\bfxi}-\bfw)\cdot\nabla\bar{\bfu}]\\
\Div\bfw=0\ea\right\}\ \ \mbox{in $\Omega\times [0,2\pi]$}\\ \medskip
\bfw=\zeta\dot{\bfxi}\ \ \mbox{at $\partial\Omega\times[0,2\pi]$}\\
\mbox{${\zeta^2}$}\ddot{\bfxi}+\omega_{\sf n}^2\bfxi+\varpi\Int{\partial\Omega}{}\mathbb T(\bfw,{\sf p})\cdot\bfn=\0\,,\ \ \mbox{in $[0,2\pi]$}\,.
\ea
\eeq{per2}
We test both sides of \eqref{per2}$_1$ by $\bfw$, and integrate by parts over $\Omega\times(0,2\pi)$. Taking into account the summability properties of members in the class $\calc$ and the definition of $\lambda_2$, we readily show
\be
\Int0{2\pi}\|\nabla\bfw(s)\|_2^2\,{\rm d}s=\lambda\left(\bar{(\zeta\dot{\bfxi}-\bfw)\cdot\nabla(\bfu_0+\bar{\bfu}),\bfw}\right)\le \mbox{$\frac{\lambda}{\lambda_2}$}\Int0{2\pi}\|\nabla\bfw(s)\|_2^2\,{\rm d}s-\lambda\left(
\bar{(\zeta\dot{\bfxi}-\bfw)\cdot\nabla\bfw,\bar{\bfu}}\right)\,.
\eeq{per3}
From  \cite[Lemma II.6.4]{Gab} we know that that there exists a ``cut-off" function, $\psi_R\in C_0^\infty(\real^3)$, $R\in(0,\infty)$,  with the following properties
\begin{itemize}
  \item [(i)] $\psi_R(x)\in [0,1]$, for all $x\in\real^3$ and $R>0$;
  \item [(ii)] $\Lim{R\to\infty}\psi_R(x)=1$ for all $x\in\real$;
  \item [(iii)] $\psi_R(x)=1$ for all $x\in B_R$, and $\supp(\psi_R)\subset B_{2R^2}$;
  \item [(iv)]  $\supp(\nabla\psi_R)\subset B^R\backslash B^{2R^2}:=S_R$, $R\ge 1$;
 \item [(v)] $\|\bfu|\nabla\psi_R|\|_2\le c\,\|\nabla\bfu\|_{2,\Omega^{\frac{R}{\sqrt2}}}$,\, with $c$ independent of $R$;
  \item [(vi)] $\partial_1\psi_R\in L^{2}(\Omega)$\,.
\end{itemize}
Testing both sides of \eqref{per1}$_1$ by $\psi_R\bar{\bfu}$ and integrating by parts over $\Omega$ as needed, we get
\be\ba{ll}\medskip
\|\psi_R^{\frac12}\nabla\bar{\bfu}\|_2^2=\half\lambda\,[-(\partial_1\psi_R\bar{\bfu},\bar{\bfu})+(|\bar{\bfu}|^2(\bfu_0+\bar{\bfu}),\nabla\psi_R)+2(\bar{\psi_R(\zeta\dot{\bfxi}-\bfw)\cdot\nabla\bfw,\bar{\bfu}})]
\\
\hspace*{3cm}
-\lambda(\psi_R\bar{\bfu}\cdot\nabla\bfu_0,\bar{\bfu})+(\bar{p}\,\bar{\bfu},\nabla\psi_R):=\Sum{k=1}5\,I_k\,.
\ea
\eeq{per4}
Using (iv)--(vi) above along with H\"older inequality, we show
$$
|I_1|+|I_2|\le \half\lambda\,\left(\|\partial_1\psi_R\|_2+\|(\bfu_0+\bar{\bfu})|\nabla\psi_R|\|_2\right)\|\bar{\bfu}\|_{4,S_R}^2\le c\,(1+\|\nabla\bfu_0\|_2+\|\nabla\bar{\bfu}\|_2)\|\bar{\bfu}\|_{4,S_R}^2
$$
and also
$$
|I_5|\le \|\bar{\bfu}|\nabla\psi_R|\|_2\|\bar{p}\|_{2,S_R}\le c\,\|\nabla\bar{\bfu}\|_2\|\bar{p}\|_{2,S_R}\,.
$$
We now pass to the limit $R\to\infty$ in \eqref{per4}. With the help of the last two displayed inequalities along with (ii) above,  the fact that $\bfu,\bfw\in \calc$ and the definition of $\lambda_1$, we easily show
\be
\|\nabla\bar{\bfu}\|_2^2=-\lambda\,(\bar{\bfu}\cdot\nabla\bfu_0,\bar{\bfu})+\lambda\bar{(\zeta\dot{\bfxi}-\bfw)\cdot\nabla\bfw,\bar{\bfu}})\le\mbox{$\frac\lambda{\lambda_1}$}\|\nabla\bar{\bfu}\|_2^2+\lambda\bar{(\zeta\dot{\bfxi}-\bfw)\cdot\nabla\bfw,\bar{\bfu}})\,.
\eeq{per5}
Summing side-by-side \eqref{per3} and \eqref{per5} and imposing $\lambda<\lambda_2$ $(\le \lambda_1)$ we thus show 
$$
\nabla{\bfu}(x,t)= \0\,,\ \ \mbox{for all $(x,t)\in \Omega\times [0,2\pi]$\,,} 
$$
which, by \eqref{1.9} and \eqref{1.10} concludes the proof 
of the theorem. 
\par\hfill$\square$\par

\setcounter{equation}{0}
\section{A General Approach to Time-Periodic Bifurcation}\label{sec:bif}
Our next objective is to formulate sufficient conditions for the occurrence of oscillatory motion of $\mathscr S$. As we mentioned in the end of Section \ref{sec:Form}, the oscillation will be induced by the time-periodic motion bifurcating from the steady-state branch. 
In order to address this problem,  we shall employ the theory developed in \cite{GaBif} for a general class of  operator equations.  The theory and its associated main result will be recalled in the current section whereas, for the reader's sake, full proofs will be reported in the Appendix. To this end, it is worth  making some introductory remarks that will also provide the appropriate rationale of our approach.
\par
Consider an evolution problem described by the equation
\be
v_t=M(v,\lambda)\,,
\eeq{Ar.0}
where $\lambda$ is a positive, real parameter and $M$ a nonlinear  operator. Suppose that in a neighborhood of $\lambda_0$, $I(\lambda_0)$, there exists a branch of steady-state solutions $v_0=v_0(\lambda)$ to \eqref{Ar.0}, namely, satisfying
$$
M(v_0,\lambda)=0\ \ \lambda\in I(\lambda_0)\,.
$$
The bifurcation problem consists then in finding sufficient conditions on $(v(\lambda_0),\lambda_0)$ guaranteeing the existence of a family of time-periodic solution around $(v_0,\lambda_0)$ that coincides with $(v_0,\lambda_0)$ as $\lambda\to\lambda_0$. This problem can be equivalently reformulated in the following, more familiar way.   
Set $u=v-v_0(\lambda)$, $\mu=\lambda-\lambda_0\in I(0)$. From \eqref{Ar.0} we then derive the following equation for $(u,\mu)$
\be
u_t+L(u)=N(u,\mu)\,,
\eeq{Ar.1}
where $L$ is the linearization of $M$ at $(v_0,\lambda_0)$  and $N$ is a nonlinear operator depending on the parameter $\mu\in I(0)$, such that $N(0,\mu)=D_uN(0,\mu)=0$ for all such $\mu$. Therefore, the bifurcation problem  at $(v_0(\lambda_0),\lambda_0)$  $(\equiv (u=0,\mu=0))$ is equivalent to show that, in addition to the trivial solution $(u=0,\mu=0)$, there exists a family of non-trivial time-periodic solutions to \eqref{Ar.1}, $u=u(\mu;t)$, of (unknown) period $T=T(\mu)$ ($T$-{\em periodic} solutions) in a neighborhood of $\mu=0$, and such that $(u(\mu;t),\mu)\to (0,0)$ as $\mu\to 0$. Setting $\tau:=2\pi\,t/T\equiv \zeta\, t$, \eqref{Ar.1} becomes
\be 
\zeta\,u_\tau+L(u)=N(u,\mu)
\eeq{Ar.2}
and the problem reduces to find a family of $2\pi$-periodic solutions to \eqref{Ar.2} with the above properties. If we  write $u=\bar{u}+(u-\bar{u}):=v+w$, where the bar denotes average over a $2\pi$-period,  it follows that \eqref{Ar.2} is formally equivalent to the following two equations
\be\ba{ll}\medskip 
L(v)=\bar{N(v+w,\mu)}:=N_1(v,w,\mu)\,,\\ \zeta\,w_\tau+L(w)=N(v+w,\mu)-\bar{N(v+w,\mu)}:=N_2(v,w,\mu)\,.\ea
\eeq{Ar.3}
We are thus led to solve for the nonlinear ``elliptic-parabolic" system \eqref{Ar.3}. It turns out that, in many circumstances,  the elliptic ``steady-state" problem is better investigated in function spaces with quite less regularity\footnote{Here `regularity' is meant in the sense of behavior at large spatial distances.} than the space where the parabolic ``oscillatory" problem should be studied. This especially occurs  when the physical system evolves in an {\em unbounded spatial region}, in which case the natural framework for the study of \eqref{Ar.3}$_1$ is a {\em homogeneous} Sobolev space, whereas that of \eqref{Ar.3}$_2$ is a classical Sobolev space \cite{GaArma,GaMaH}.
This suggests that, in general, it is more appropriate to study the two equations in \eqref{Ar.3} in two {\em distinct} function classes. As a consequence, even though being {\em formally}  the same operator, the operator $L$ in \eqref{Ar.3}$_1$ acts on and ranges into spaces  different than those where the operator $L$ in \eqref{Ar.3}$_2$ acts and ranges. With this in mind, \eqref{Ar.3} becomes
$$
L_1(v)=N_1(v,w,\mu)\,;\ \ \zeta\,w_\tau+L_2(w)=N_2(v,w,\mu)\,.
$$   
\par
The general theory proposed in \cite{GaBif},  and that we are about to recall in its main aspects, takes its cue  exactly from the these  considerations. 
\par  
Let ${\calu }, {\calv}$,  be real Banach spaces with norms $\|\cdot\|_{\mathcal U}$, $\|\cdot\|_{\mathcal V}$, respectively, and let  $\mathcal Z$ be a real Hilbert space with norm $\|\cdot\|_{\mathcal Z}$.
Moreover,  let 
$$
L_1:\mathcal U\mapsto \mathcal V\,
$$
be a bounded linear operator, and let \footnote{Given an operator $M$, we denote by ${\sf D}[M]$ and ${\sf R}[M]$ its domain and range, respectively,  and by ${\sf N}[M]:=\{u: M(u)=0\}$ its null space.}
$$
L_2:{\sf D}\,[L_2]\subset \mathcal Z\mapsto \mathcal Z\,,
$$
be a densely defined, closed linear operator, with a non-empty resolvent set ${\sf P}(L_2)$. For a fixed (once and for all) $\theta\in {\sf P}(L_2)$ we denote by $\mathcal W$ the linear subspace of $\mathcal Z$ closed under the (graph) norm $\|w\|_{\mathcal W}:=\|(L_2+\theta\,I)w\|_{\mathcal Z}$, where $I$ stands for the identity operator. We also define the following spaces
$$\ba{rl}\medskip
\mathcal Z_{2\pi,0}&\!\!\!:=\left\{w:\real\to\calz,\, \mbox{$2\pi$-periodic with} \ \bar{w}=0, \,\mbox{and}\,\int_{-\pi}^\pi\|w(s)\|_{\calz}^2{\rm d}s<\infty \right\} 
\\\medskip\calw_{2\pi,0}&\!\!\!:=\left\{w\in \calz_{2\pi,0}:\,\|w\|_{\calw_{2\pi,0}}:=\left(\int_{-\pi}^\pi\left(\|w_t(s)\|_{\calz}^2+\|w(s)\|_{\calw}^2\right){\rm d}s\right)^{\frac12}<\infty\right\}\,.
\ea$$ 
Next, let
$$
N: \calu\times \calw_{2\pi,0}\times \real\mapsto \calv\oplus \calz_{2\pi,0}   
$$
be a (nonlinear) map satisfying the following properties:
\be\ba{rl}\medskip
N_1&\!\!\!: (v,w,\mu)\in\calu\times \calw_{2\pi,0}\times \real\mapsto \bar{N(v,w,\mu)}\in \calv
\\
N_2&\!\!\!:=N-N_1:\calu\times \calw_{2\pi,0}\times \real\mapsto \calz_{2\pi,0}\,.
\ea
\eeq{Ar.4} 
The Bifurcation Problem is then formulated as follows.\smallskip\par 
{\em Find a neighborhood of the origin $U(0,0,0)\subset \calu\times \calw_{2\pi,0}\times \real$ such that the equations
\be
L_1(v)=N_1(v,w,\mu)\,,\ \mbox{in $\calv$}\,;\ \ \zeta\, w_\tau +L_2(w)=N_2(v,w,\mu)\,,\ \mbox{in $\calz_{2\pi,0}$}\,,
\eeq{Ar.5}
possess there a family of non-trivial $2\pi$-periodic solutions $(v(\mu),w(\mu;\tau))$ for some $\zeta=\zeta(\mu)>0$,  such that $(v(\mu),w(\mu;\cdot))\to 0$ in $\calu\times\calw_{2\pi,0}$ as $\mu\to0$.}\smallskip\par
Whenever the Bifurcation Problem is solvable, we say that  $(u=0,\mu=0)$ is a {\em bifurcation point}. Moreover, the bifurcation is called {\em supercritical} [resp. {\em subcritical}] if the family of solutions $(v(\mu),w(\mu;\tau))$ exists only for $\mu>0$ [resp. $\mu<0$]. 
\medskip\par
In order to provide sufficient conditions for the resolution of the above problem, we make the following  assumptions (H1)--(H4) on the involved operators.\footnote{As customary, For $X$ a real Banach space, we denote by $X_\comps: = X + {\rm i} X$ its complexification. If $L$ is
a linear operator on $X$, we continue to denote by $L$ its extension to a linear operator on
$X_\comps$, while ${\sf N}_\comps[L]$ and
${\sf R}_\comps[L]$ stand for the null space and range in $X_\comps$. The spectrum, $\sigma(L)$, is computed
with respect to $X_\comps$, so $\nu\in\sigma(L)$ if and only if its complex conjugate, $\nu^*$, is in $\sigma(L)$.\label{foot:complex}}
\begin{itemize}
  \item[(H1)] $L_1$ is a homeomorphism\,;
  \item[(H2)] There exists $\nu_0:={\rm i}\,\zeta_0$, $\zeta_0>0$ such that   $L_2-\nu_0I$ is Fredholm of index 0,  and   ${\rm dim}\,{\sf N}_\comps[L_2-\nu_0I]=1$ with ${\sf N}_\comps[L_2-\nu_0I]\cap {\sf R}_{\comps}[L_2-\nu_0I]=\{0\}$. Namely, $\nu_0$ is a simple  eigenvalue of $L_2$. Moreover,  $k\,\nu_0\in {\sf P}(L_2)$, for all $k\in\nat\backslash\{0,1\}$. 
  \item[(H3)] The operator
$$
\mathscr Q:w\in \calw_{2\pi,0}\mapsto \zeta_0\,w_\tau+L_2(w)\in \calz_{2\pi,0}\,,
$$  
is Fredholm of index 0\,;
\item[(H4)] The nonlinear operators $N_1,N_2$
are analytic in the neighborhood $U_1(0,0,0)\subset \calu\times \calw_{2\pi,0}\times \real$, namely, there exists $\delta>0$ such that for all $(v,w,\mu)$ with $\|v\|_{\calu}+\|w\|_{\calw_{2\pi,0}}+|\mu|<\delta$, the Taylor series 
$$\ba{ll}\medskip
N_1(v,w,\mu)=\Sum{k,l,m=0}{\infty}R_{klm}v^kw^l\mu^m\,,\\
N_2(v,w,\mu)=\Sum{k,l,m=0}{\infty}S_{klm}v^kw^l\mu^m\,,
\ea
$$
are absolutely convergent in $\calv$ and $\calz_{2\pi,0}$, respectively, for all $(v,w,\mu)\in U_1$. Moreover, we assume that the multi-linear operators $R_{klm}$ and $S_{klm}$ satisfy $R_{klm}=S_{klm}=0$ whenever $k+l+m\le1$, and $R_{011}=R_{00m}=S_{00m}=0$, all $m\ge2$.
\end{itemize}
\par\noindent
Let
\be
L_2(\mu):=L_2+\mu\,S_{011}\,,
\eeq{L2mu}
and observe that, by (H2), $\nu_0$ is a simple eigenvalue of $L_2(0)\equiv L_2$. Therefore, denoting by $\nu(\mu)$ the eigenvalues of $L_2(\mu)$, it follows (e.g.  \cite[Proposition 79.15 and Corollary 79.16]{Z1}) that in a neighborhood of $\mu=0$ the map $\mu\mapsto\nu(\mu)$ is well defined and of class $C^\infty$.
\renewcommand{\theequation}{\arabic{section}.\arabic{equation}}
\smallskip\par\setcounter{equation}{5}
We may now state the following bifurcation result.
\Bt\label{3.1_ar} Suppose  {\rm (H1)--(H4)} hold and, in addition,
\be
\Re[\nu'(0)]\neq 0\,,
\eeq{nupr}
namely, the  eigenvalue $\nu(\mu)$ crosses the imaginary axis with ``non-zero speed." Moreover, let $v_0$ be a normalized eigenvector of $L_2$ corresponding to the eigenvalue $\nu_0$, and set
$
v_1:=\Re[v_0\,{\rm e}^{-{\rm i}\,\tau}].$ 
Then, the following properties are valid. \smallskip\\
{\rm (a)} {\rm Existence.} There are analytic families
\be
\big(v(\varepsilon),w(\varepsilon),\zeta(\varepsilon),\mu(\varepsilon)\big)\in \calu\times \calw_{2\pi,0}\times \real_+\times\real
\eeq{fam}
satisfying \eqref{Ar.5}, for all $\varepsilon$ in a neighborhood $\mathcal I(0)$ of $0\in\real$, and such that
\be
\big(v(\varepsilon),w(\varepsilon)-\varepsilon\,v_1,\zeta(\varepsilon),\mu(\varepsilon)\big)\to (0,0,\zeta_0,0)\ \ \mbox{as $\varepsilon\to 0$}\,.
\eeq{Ar.10}
\par\noindent
{\rm (a)} {\rm Uniqueness.}
There is a neighborhood  $$U(0,0,\zeta_0,0)\subset \calu\times \calw_{2\pi,0}\times \real_+\times \real$$ such that every (nontrivial) $2\pi$-periodic solution to \eqref{Ar.5},  lying in $U$ must coincide, up to a phase shift, with a member of the family \eqref{fam}.
\smallskip\par\noindent
{\rm (a)} {\rm Parity.}  The functions $\zeta(\varepsilon)$ and $\mu(\varepsilon)$ are even:
$$
\zeta(\varepsilon)=\zeta(-\varepsilon)\,,\ \ \mu(\varepsilon)=\mu(-\varepsilon)\,,\ \ \mbox{for all $\varepsilon\in\cali(0)$\,.} 
$$
Consequently, the bifurcation due to these solutions is either subcritical or supercritical, a two-sided bifurcation being excluded.\footnote{Unless $\mu\equiv 0$.}
\ET{3.1_ar}
\par
A proof of this result is given in \cite[Theorem 3.1]{GaBif}. However, for completeness and reader's sake, we deem it appropriate to reproduce it in the Appendix. 

\setcounter{equation}{0}
\section{Existence of an Analytic  Steady-State Branch}\label{steady}
Our next main objective is to reformulate  the bifurcation problem presented in Section~\ref{sec:Form} in a suitable functional setting, which will eventually allow us to apply the theory presented in Section~\ref{sec:bif}. The major  challenge in reaching this goal is the choice of the `right' setting, where the operators involved possess the properties required by \theoref{3.1_ar}. This will require a careful study of appropriate linearized problems obtained from \eqref{03} and \eqref{04} that will be performed in Sections \ref{steady} through \ref{sec:time_per}. \par
In this section we will provide sufficient conditions for the existence of an analytic branch of solutions to \eqref{03} in a neighborhood of some $\lambda_0$.  
To this end, we begin to define the map:
$$
\hat{\Delta}:(\bfu,\bfchi)\in X^{2}(\Omega)\times\real^3\mapsto \tilde{\Delta}(\bfw,\bfxi)\in \call^2(\real^3) 
$$
where 
\be
\hat{\Delta}(\bfu,\bfchi)=\left\{\ba{ll}\medskip
 -\Delta\bfw\ \mbox{in $\Omega$}\,,\\ 
\omega_{\sf n}^2\bfchi+2\varpi\Int{\partial\Omega}{}\mathbb D(\bfw)\cdot\bfn\ \ \mbox{in $\Omega_0$\,,}\ea\right.
\eeq{Deltah}
and set
$\hat{\mathscr A}:=\mathscr P\,\hat{\Delta}$, where, we recall, $\mathscr P$ is the orthogonal projection of $\call^2(\real^3)$ onto $\calh(\real^3)$.  
We also define
$$
\hat{\partial_1}:\bfu\in X^{2}(\Omega)\mapsto \hat{\partial_1}(\bfu)\in \calh(\real^3) 
$$
where
\be
\hat{\partial_1}(\bfu)=\left\{\ba{ll}\medskip
 -\partial_1\bfu\ \mbox{in $\Omega$}\,,\\ 
\0\ \ \mbox{in $\Omega_0$\,.}\ea\right.
\eeq{D1h}
It is readily checked that, as stated, $\hat{\partial}_1(\bfu)\in \calh(\real^3)$. By \lemmref{0.1} and \eqref{D1h}, this amounts to show that
\be
\int_{\Omega}\partial_1\bfu\cdot\nabla p=0\,,\ \ \mbox{for all $p\in D^{1,2}(\Omega)$}\,.
\eeq{wwh}  
Since $\bfu\in\cald^{1,2}_0(\Omega)$, there is a sequence $\{\bfu_k\}\subset \calc_0(\Omega)$ such that $\|\nabla(\bfu_k-\bfu)\|_2\to 0$ as $k\to\infty$. Then, by an integration by parts combined with the condition $\Div\bfu_k=0$, we show that \eqref{wwh} holds for each $\bfu_k$ and hence for $\bfu$ after passing to the limit $k\to\infty$. Let $\bsfu_0$ be the velocity field associated to the solution to \eqref{03} determined in \theoref{exi}, corresponding to $\lambda=\lambda_0$, and set
$$
\mathscr C:\bfu\in X^{2}(\Omega)\mapsto \mathscr C(\bfu)\in \call^2(\real^3) 
$$
where
\be
\mathscr C(\bfu)=\left\{\ba{ll}\medskip
 \bsfu_0\cdot\nabla\bfu+\bfu\cdot\nabla\bsfu_0\ \mbox{in $\Omega$}\,,\\ 
\0\ \ \mbox{in $\Omega_0$\,.}\ea\right.\,.
\eeq{Kh}
\par
Finally,
let ${\bsfU}:=(\bfu,\bfchi)$ and consider  the following operator
\be
\mathscr L_1:\bsfU\in X^2(\Omega)\times \real^3\mapsto\hat{\mathscr A}(\bsfU) +\lambda_0[\hat{\partial_1}\bfu-\mathscr P\mathscr C(\bfu)]\in \calh(\real^3)\cap\cald_0^{-1,2}(\Omega)\,.
\eeq{2.2}
\Bl
The operator $\mathscr L_1$ is well-defined.
\EL{L1}
{\em Proof.} Clearly $\mathscr P\bfphi=\bfphi$,  for  $\bfphi\in \mathcal C_0(\Omega)$. Therefore, for any such a $\bfphi$, integrating by parts we get
$$
\left(\mathscr L_1(\bsfU),\bfphi\right)=\lambda_0(\partial_1\bfu,\bfphi)-\left(\lambda_0(\bsfu_0\otimes\bfu+\bfu\otimes\bsfu_0)+\nabla\bfu,\nabla\bfphi\right).
$$
Employing  H\"older inequality and \propref{3.1} we deduce that since $\bfu\in X(\Omega)$, then $\mathscr L_1(\bsfU)\in\cald_0^{-1,2}(\Omega)$.
Moreover,
$$
\|\bsfu_0\cdot\nabla\bfu+\bfu\cdot\nabla\bsfu_0\|_2\le \|\bsfu_0\|_\infty\|\nabla\bfu\|_2+\|\bfu\|_\infty\|\nabla\bsfu_0\|_2\,,
$$
which implies, by \lemmref{Xemb} and \theoref{exi},  that since $\bfu\in X^2(\Omega)$, then $\mathscr L_1(\bsfU)\in \calh(\real^3)$.\par\hfill$\square$
\par
\Bl  The operator $\mathscr L_1$ is Fredholm of index 0.
\EL{fredo}
{\em Proof.} We have $\mathscr L_1=\mathscr L_1^0+\mathscr P\mathscr C$ where
$\mathscr L_1^0:=\hat{\mathscr A}+\lambda_0\hat{\partial_1}.$  
The operator $\mathscr L_1^{0}$ is a homeomorphism, that is, for any $(\bff,\bfF)\in Y(\Omega)\times\real^3$ there exists unique $(\bfu,\bfchi)\in X^2(\Omega)\times \real^3$ solving
\be\ba{cc}\medskip\left.\ba{ll}\medskip
\lambda_0\,\partial_1\bfu+\Delta\bfu=\nabla p+\bff\\ \medskip
\Div\bfu=0\ea\right\}\ \ \mbox{in $\Omega$}\,,\\ \medskip
\bfu(x)=\0\ \ \mbox{at $\partial\Omega$}\,,\\
\omega_{\sf n}^2\bfchi+\varpi\Int{\partial\Omega}{}\mathbb T(\bfu,p)\cdot\bfn=\bfF
\ea
\eeq{chia_p}
and continuously depending on the data. To show the validity of this property, 
we notice that from \cite[Theorem 2.1]{GaEi} it follows that, corresponding to the given $\bff\in Y(\Omega)$, there exists a unique solution $(\bfu,p)\in X^2(\Omega)\times [D^{1,2}(\Omega)\cap L^6(\Omega)]$ to \eqref{chia_p}$_{1,2,3}$  such that
\be
\|\bfu\|_{X^2}+\|p\|_6+\|\nabla p\|_{2}\le c\,\|\bff\|_Y\,.
\eeq{EiG}
By the trace theorem we get 
\be
\left|\int_{\partial\Omega}\mathbb T(\bfu,p)\cdot\bfn\right|\le c\,(\|\bfu\|_{2,2,\Omega_R}+\|p\|_{1,2,\Omega_R})<\infty\,,
\eeq{2.3_s}
The associated and uniquely determined displacement $\bfchi$ is then  obtained from \eqref{chia_p}$_4$. In view of \eqref{EiG} and \eqref{2.3_s} we have also
$$
|\bfchi|\le c\,(\|\bff\|_Y+|\bfF|)
$$ 
which concludes the proof of the homeomorphism property of $\mathscr L_1^0$.
We shall next show that the operator $\mathscr C$ is compact. Let $\{\bfu_k\}\subset X^2(\Omega)$ with
\be
\|\bfu_k\|_{X^2}\le M\,,
\eeq{unibou}
and $M$ independent of $k\in\nat$.
This implies, in particular, that there exist  an element $\bfu_*\in L^4(\Omega)\cap D^{1,2}(\Omega)\cap D^{2,2}(\Omega)$ such that along a subsequence (that we continue to denote by $\{\bfu_k\}$)
\be\ba{ll}\medskip
\bfu_k\to \bfu_*\,,\ \ \mbox{weakly in $L^4(\Omega)$}\,,\\
\nabla\bfu_k\to \nabla\bfu_*\,,\ \ \mbox{weakly in $W^{1,2}(\Omega)$}\,.
\ea
\eeq{weakc}
Moreover, by \lemmref{Xemb} and classical compact embedding results, we also have
\be
\bfu_k\to \bfu_*\,,\ \ \mbox{strongly in $W^{1,2}(\Omega_R)$, for all $R>R_*$.}
\eeq{strongc}
Setting $\bfw_k:=\bfu_k-\bfu_*$,  we get
$$
|\bsfu_0\cdot\nabla\bfw_k+\bfw_k\cdot\nabla\bsfu_0|_{-1,2}\le\|\bsfu_0\otimes\bfw_k\|_2\le  \|\bfu_0\|_\infty\|\bfw_k\|_{2,\Omega_R}+\|\bsfu_0\|_{4,\Omega^R}\|\bfw_k\|_4\,,
$$
and
$$
\|\bsfu_0\cdot\nabla\bfw_k+\bfw_k\cdot\nabla\bsfu_0\|_{2}\le \|\bsfu_0\|_{1,\infty}\|\bfw_k\|_{1,2,\Omega_R}+\|\bsfu_0\|_{1,4,\Omega^R}\|\bfw_k\|_{1,4}\,.
$$
From \theoref{exi}, we know $\bsfu_0\in W^{1,\infty}(\Omega)\cap W^{1,4}(\Omega)$, so that letting first $k\to\infty$ in the above two inequalities and using \eqref{unibou}--\eqref{strongc}, and then $R\to\infty$, we show
$$
\lim_{k\to\infty}\left(|\bsfu_0\cdot\nabla\bfw_k+\bfw_k\cdot\nabla\bsfu_0|_{-1,2}+\|\bsfu_0\cdot\nabla\bfw_k+\bfw_k\cdot\nabla\bsfu_0\|_{2}\right)=0\,.
$$
This, in view of \eqref{Kh},  proves that $\mathscr C$ is compact.
As a result,  we then conclude the claimed Fredholm property for $\mathscr L_1$.\par\hfill$\square$\par
We are now in a position to prove the following result, representing the main contribution of this section.
\Bt Let $(\bsfu_0,{\sf p}_0,{\sf \bfchi}_0)$ be the solution corresponding to $\lambda=\lambda_0$ determined in \theoref{exi}, and  assume that the equation
$$
\mathscr L_1(\bsfU)=\0
$$
has only the solution $\bsfU\equiv\0$. Then there exists a neighborhood of $\lambda_0$, $U(\lambda_0)$, such that \eqref{03} has an analytic branch of solutions $(\bfu_0(\lambda),p_0(\lambda),\bfchi_0(\lambda))\in X^{2}(\Omega)\times D^{1,2}(\Omega)\times\real^3$, $\lambda\in U(\lambda_0)$, with the property
$$
(\bfu_0(\lambda),p_0(\lambda),\bfchi_0(\lambda))\to (\bsfu_0,{\sf p}_0,{\bfchi}_0)\ \ \mbox{as $\lambda\to\lambda_0$\,.}
$$
\ET{2.1}
{\em Proof.}   We write the solution $(\bfu_0(\lambda),p_0(\lambda),\bfchi_0(\lambda))$ to \eqref{03} as
$(\bfu(\lambda)+\bsfu_0,p(\lambda)+{\sf p}_0,\bfchi(\lambda)+\bftau_0)$,  set $\bsfU:=(\bfu,\bftau_0)$,  
$\mathcal X:= X^{2}(\Omega)\times\real^3$, $\mathcal Y:= \cald_0^{-1,2}(\Omega)\cap\calh(\real^3)$, 
$\mu:=\lambda-\lambda_0,$ and 
$$
\mathscr M:\bfu\in X^{2}(\Omega)\mapsto \mathscr M(\bfu)\in \call^2(\real^3) 
$$
where
\be
{\mathscr M}(\bfu)=\left\{\ba{ll}\medskip
 \left(\partial_1(\bfu+
\bsfu_0)-(\bsfu_0+\bfu)\cdot\nabla(\bsfu_0+\bfu)\right)\ \mbox{in $\Omega$}\,,\\ 
\0\ \ \mbox{in $\Omega_0$\,.}\ea\right.\,.
\eeq{M}
We then define the map
$$ 
\mathcal F:(\bsfU,\mu)\in \mathcal X\times \cali(0)\mapsto \mathscr L_1(\bsfU)+\mu\mathscr P\mathscr M(\bfu)\in\mathcal Y\,,
$$
where $\cali(0)$ is a neighborhood of $0\in\real$.
The map $\calf$ is well defined. In fact, taking into account \theoref{exi} and that $\bfu\in X^2(\Omega)$, by arguments similar to those employed previously we easily show
\be
\partial_1\bfu-(\bsfu_0+\bfu)\cdot\nabla(\bsfu_0+\bfu)\in \mathcal Y\,,\ \ \partial_1\bsfu_0\in L^2(\Omega)\,.
\eeq{mcsnf}
Moreover, testing \eqref{03}$_1$ (with $\bfu_0\equiv \bsfu_0$, $p_0\equiv {\sf p}_0$, $\lambda\equiv\lambda_0$) by $\bfphi\in \calc_0(\Omega)$ and integrating by parts, we get
$$
\lambda_0(\partial_1\bsfu_0,\bfphi)=-(\bsfu_0\otimes\bsfu_0,\nabla\bfphi)+(\nabla\bsfu_0,\nabla\bfphi)\,.
$$
Since $\bsfu_0\in L^4(\Omega)\cap D^{1,2}(\Omega)$, from the latter we deduce $\partial_1\bsfu_0\in \cald_0^{-1,2}(\Omega)$, which, along with \eqref{mcsnf} proves that $\calf$ is well defined.
Now, from \eqref{03} it is checked at once that   $\bfU$ satisfies the equation
\be
\mathcal F(\bfU,\mu)=\0\in \mathcal Y\,.
\eeq{2.5}
Clearly, \eqref{2.5} has the solution $(\bfU=\0, \mu=0)$. Moreover, $\mathcal F$ is Frech\'et-differentiable with derivative $D_{{\mbox{\tiny $\bfU$}}}\mathcal F(\0,0)=\mathscr L_1$. Since $\mathscr L_1$ is Fredholm of index 0, the assumption made on $\mathscr L_1$ implies that $\mathscr L_1$ is a homeomorphism. As a consequence, by the analytic version of the Implicit Function Theorem we show the  property stated in the second part of the theorem, which is thus completely proved. 
\par\hfill$\square$\par
\setcounter{equation}{0}
\section{Spectral Properties of the Linearized Operator}\label{sec:spectrum}
The main goal of this section is to secure some relevant spectral properties of the  operator obtained by linearizing \eqref{04} around the trivial solution  $\bfw\equiv{\sf p}\equiv\bfxi\equiv\bfsigma\equiv\0$.  In this regard, we begin to 
 define the maps 
$$
\tilde{\Delta}:(\bfw,\bfxi)\in Z^{2,2}\times\real^3\mapsto \tilde{\Delta}(\bfw,\bfxi)\in \call^2(\real^3) 
$$
where $Z^{2,2}$ is introduced in \eqref{zetaa},
\be
\tilde{\Delta}(\bfw,\bfxi)=\left\{\ba{ll}\medskip
 -\Delta\bfw\ \mbox{in $\Omega$}\,,\\ 
\omega_{\sf n}^2\bfxi+2\varpi\Int{\partial\Omega}{}\mathbb D(\bfw)\cdot\bfn\ \ \mbox{in $\Omega_0$\,,}\ea\right.
\eeq{Delta}
and
$$
\tilde{\partial_1}:\bfw\in Z^{2,2}\mapsto \tilde{\partial_1}(\bfw)\in \calh(\real^3) 
$$
where
\be
\tilde{\partial_1}(\bfw)=\left\{\ba{ll}\medskip
 -\partial_1\bfw\ \mbox{in $\Omega$}\,,\\ 
\0\ \ \mbox{in $\Omega_0$\,.}\ea\right.
\eeq{D1}
It is readily checked that, as stated, $\tilde{\partial}_1(\bfw)\in \calh(\real^3)$. This can be proved exactly in the same way we did for the operator \eqref{D1h}. 
Set $\mathscr A:=\mathscr P\,\tilde{\Delta}$, with $\mathscr P$  the orthogonal projection of $\call^2(\real^3)$ onto $\calh(\real^3)$,  ${\bsfW}:=(\bfw,\bfxi)$, $\bfsigma:=\bfw|_{\Omega_0}$, and consider the operator 
\be 
\mathscr L_0:{\bsfW}\in Z^{2,2}\times \real^3\subset \calh(\real^3)\times\real^3\mapsto\mathscr L_0({\bsfW})= (\lambda_0\tilde{\partial}_1\bfw+\mathscr A(\bfw,\bfxi),-\bfsigma)\in \calh(\real^3)\times\real^3\,.
\eeq{L0}
\par
The following result holds.\footnote{See Footnote \ref{foot:complex}.} \Bl Let $\zeta\in\real\backslash\{0\}$. Then, the operator $\mathscr L_0-\i\,\zeta\cali$, with $\cali$ identity in $\calh_{\mbox{\tiny $\comp$}}(\real^3)\times\mathbb C^3$, is a homeomorphism of $Z_{\mbox{\tiny $\mathbb C$}}^{2,2}\times\comp^3$ onto $\calh_{\mbox{\tiny $\comp$}}(\mathbb R^3)\times\comp^3$.  Moreover, there is $c=c(\Omega_0,\omega_{\sf n},\varpi,\lambda_0)>0$, such that
\be
\|D^2\bfw\|_2+|\zeta|^{\frac12}\|\nabla\bfw\|_2+|\zeta|(\|\bfw\|_2+|\bfsigma|+|\bfxi|)\le c\,\|(\mathscr L_0-\i\,\zeta)(\bsfW)\|_2\,,\ \ |\zeta|\ge 1\,.
\eeq{4.2}
\EL{4.1}
{\em Proof.} In view of  \lemmref{0.1},  \eqref{Delta} and \eqref{D1} it readily follows  that, for a given $\bfcalf:=((\bff,\bfF),\bfG)\in\calh_{\mbox{\tiny $\mathbb C$}}(\real^3)\times\comp^3$, the equation $$(\mathscr L_0-\i\,\zeta)(\bsfW)=\bfcalf$$ 
is equivalent to the following
set of equations
\be\ba{cc}\medskip\left.\ba{ll}\medskip\Delta\bfw+\lambda_0\partial_1\bfw-\nabla{\mathfrak p}=-\i\,\zeta\,\bfw -\bff\\
\Div\bfw=0\ea\right\}\,  \mbox{in $\Omega$}\,,\\ \medskip
\bfw=\bfsigma\ \ \mbox{at $\partial\Omega$}\,,\\ 
-\i\,\zeta\,{\bfsigma}+\omega^2_{\sf n}\bfxi+\varpi\Int{\partial\Omega}{}\mathbb T(\bfw,{\mathfrak p})\cdot\bfn=\bfF\,,\ \   
\bfsigma=-\i\,\zeta\bfxi-\bfG\,.
\ea
\eeq{4.1}
Therefore, the homeomorphism property follows if we show that, for any given $((\bff,\bfF),\bfG)$ specified above, problem \eqref{4.1}  
has one and only one solution $(\bfw,\mathfrak p,\bfsigma,\bfxi)\in W_{\mbox{\tiny $\mathbb C$}}^{2,2}(\Omega)\times D^{1,2}_{\mbox{\tiny $\mathbb C$}}(\Omega)\times\mathbb C^3\times\mathbb C^3$. We begin to establish some formal estimates, holding for all $\zeta\neq 0$. In this regards, we recall that, by \lemmref{1.1},
\be 
|\bfsigma|\le c_0\|\mathbb D(\bfw)\|_2
\eeq{tra}
with $c_0=c_0(\Omega_0)>0$. 
Let us dot-multiply both sides of \eqref{4.1}$_1$ by $\bfw^*$ (${}^*=$\,complex conjugate) and integrate by parts over $\Omega$. On account of \eqref{4.1}$_{2-5}$ we infer
\be
2\|\mathbb D(\bfw)\|_2^2-\i\,\zeta\,(\|\bfw\|_2^2+\frac1\varpi|\bfsigma|^2-\frac{\omega_{\sf n}^2}{\varpi}|\bfxi|^2)-\lambda_0(\partial_1\bfw,\bfw^*)=(\bff,\bfw^*)+\frac1\varpi\bfF\cdot\bfsigma^*+\frac{\omega_{\sf n}^2}{\varpi}\bfG^*\cdot\bfxi\,.
\eeq{4.3}
If we take the real part of \eqref{4.2},  use \eqref{1.8} and Schwarz inequality, and observe that $\Re \,(\partial_1\bfw,\bfw^*)=0$, we deduce
\be
\|\nabla\bfw\|_2^2\le c_1 \,\left(\|\bff\|_2\|\bfw\|_2+|\bfF||\bfsigma|+|\bfG||\bfxi|\right)
\eeq{4.4}
where, here and in the rest of the proof, $c_1$ denotes a positive constant depending, at most, on $\Omega_0$, $\omega_{\sf n}$, $\lambda_0$, and $\varpi$.
From \eqref{4.1}$_5$ and \eqref{tra} we deduce
\be
|\bfxi|\le {c_1}{|\zeta|^{-1}}\left(|\bfG|+\|\nabla\bfw\|_2\right)\,,
\eeq{xi}
so that employing in \eqref{4.4} the latter, \eqref{tra}, \eqref{1.8} and Cauchy-Schwarz inequality we show
$$
\|\nabla\bfw\|_2^2\le c_1\left[(|\zeta|^{-1}+|\zeta|^{-2})|\bfG|^2+\|\bff\|_2\|\bfw\|_2+|\bfF||\bfsigma|\right]\,.
$$
If we apply Cauchy inequality on the last two terms on the right-hand side of this inequality we may deduce, on the one hand, 
\be
\|\nabla\bfw\|_2^2\le c_1\left[|\zeta|^{-1}(|\bfF|^2+|\bfG|^2+\|\bff\|_2^2)+|\zeta|^{-2}|\bfG|^2\right]+|\zeta|(|\bfsigma|^2+\|\bfw\|_2^2)
\eeq{Du0}
and, on the other hand,
\be
\|\nabla\bfw\|_2^2\le c_1\left[|\zeta|^{-1}|\bfG|^2+|\zeta|^{-2}\left(|\bfF|^2+|\bfG|^2+\|\bff\|_2^2\right)\right]+\varepsilon |\zeta|^2(|\bfsigma|^2+\|\bfw\|_2^2)\,.
\eeq{Du}
where $\varepsilon>0$ is arbitrarily fixed and $c_1$ depends also on $\varepsilon$. 
We now take the imaginary part of \eqref{4.2} and use Schwarz inequality to get
$$
|\zeta|(|\bfsigma|^2+\|\bfw\|_2^2)\le c_1\left[|\zeta||\bfxi|^2+\|\nabla\bfw\|_2\|\bfw\|_2+\|\bff\|_2\|\bfw\|_2+|\bfG||\bfxi|+|\bfF||\bfsigma|\right]\,.
$$
By utilizing, on the right-hand side of the latter, \eqref{xi} and Cauchy inequality we obtain
\be
|\zeta|(|\bfsigma|^2+\|\bfw\|_2^2)\le {c_1}\,|\zeta|^{-1}\left[\|\nabla\bfw\|_2^2+|\bfG|^2+|\bfF|^2+\|\bff\|_2^2\right]\,.
\eeq{ima}
Thus, combining \eqref{Du} and \eqref{ima}, and taking $\varepsilon$ suitably small we deduce  
\be
|\zeta|^2(|\bfsigma|^2+\|\bfw\|_2^2)\le {c_1}(1+|\zeta|^{-1}+|\zeta|^{-2})\left(|\bfG|^2+|\bfF|^2+\|\bff\|_2^2\right)\,,
\eeq{ima1}
which, in turn, once replaced in \eqref{Du0}, delivers
\be
|\zeta|\|\nabla\bfw\|_2^2\le {c_1}(1+|\zeta|^{-1}+|\zeta|^{-2})\left(|\bfG|^2+|\bfF|^2+\|\bff\|_2^2\right)\,.
\eeq{Du1}
Finally, from \eqref{xi} and \eqref{Du1} we conclude
\be
|\zeta|^2|\bfxi|^2\le c_1(1+|\zeta|^{-1}+|\zeta|^{-2}+|\zeta|^{-3})\left(|\bfG|^2+|\bfF|^2+\|\bff\|_2^2\right)\,.
\eeq{xi1}
Combining estimates \eqref{ima1}--\eqref{xi1} with classical Galerkin method, we can show that for any given $((\bff,\bfF),\bfG)$ in the specified class and $\zeta\neq 0$, there exists a (unique, weak) solution to \eqref{4.1} such that $(\bfw,\bfxi)\in [\cald^{1,2}_{\mbox{\tiny $\mathbb C$}}(\real^3)\cap\call^2_{\mbox{\tiny $\mathbb C$}}(\real^3)]\times\compl^3$, satisfying  \eqref{ima1}--\eqref{xi1}. We next write \eqref{4.1}$_{1-3}$ in the following Stokes-problem form
$$\ba{cc}\medskip\left.\ba{ll}\medskip\Delta\bfw=\nabla{\mathfrak p}+\bfcalg\\
\Div\bfw=0\ea\right\}\,  \mbox{in $\Omega$}\,,\\ \medskip
\bfw=\bfsigma\ \ \mbox{at $\partial\Omega$}
\ea
$$
where
$$
\bfcalg:=-\lambda_0\, \partial_1\bfw +\i\,\zeta\,\bfw-\bff\,.
$$
Since $\bfcalg\in L^2_{\mbox{\tiny $\mathbb C$}}(\Omega)$ and $\bfu\in W^{1,2}_{\mbox{\tiny $\mathbb C$}}(\Omega)$, from classical results \cite[Theorems IV.5.1 and V.5.3]{Gab} it follows that $D^2\bfw\in L^2(\Omega)$, thus completing the existence (and uniqueness) proof. Furthermore, by \cite[Lemma IV.1.1 and V.4.3]{Gab} we get  
\be
\|D^2\bfw\|_2\le c\,\big[\|\bff\|_2+(\lambda_0+1)\|\nabla\bfw\|_2+(|\zeta|+1)\,\|\bfw\|_2+|\bfsigma|\big]\,.
\eeq{4.7}
As a result, if $|\zeta|\ge 1$, the inequality in \eqref{4.2} is a consequence of \eqref{ima1}--\eqref{4.7}, and the proof of the lemma is completed.
\par\hfill$\square$\par
Let $\bsfu_0\in X^{2}(\Omega)$ and consider the operator 
\be
\mathscr K: \bsfW:=(\bfw,\bfxi)\in  Z^{2,2}\times\real^3 \mapsto \mathscr K(\bsfW)\in \call^2(\real^3)
\eeq{3.32}
where ($\bfsigma=\bfw|_{\Omega_0}$), 
\be
\mathscr K(\bsfW)=\left\{\ba{ll}\medskip
\lambda_0\big(\bsfu_0\cdot\nabla\bfw+(\bfw-\bfsigma)\cdot\nabla\bsfu_0\big)\ \ \mbox{in $\Omega$}\,,\\
\0\ \ \mbox{in $\Omega_0$}\,.
\ea\right.
\eeq{3.32_1}
By \propref{3.1} and classical results, we have
\be\bsfu_0\in L^4(\Omega)\cap D^{1,2}(\Omega)\,,\ \ W^{2,2}(\Omega)\subset  L^\infty(\Omega)\cap W^{1,4}(\Omega)\,,
\eeq{u0w} 
so that we easily check that the operator $\mathscr K$ is well defined.
Next, let
\be
\mathscr L_2:\bsfW\in Z^{2,2}\times \real^3\mapsto \mathscr L_2(\bsfW):=\mathscr L_0(\bsfW) +\mathscr P\mathscr K(\bsfW)\in \calh(\real^3)\,.
\eeq{L2}
The main result of this section reads as follows.
\Bt The operator
$ 
\mathscr L_2-\i\,\zeta\cali,
$ with $\cali$ as in \lemmref{4.1},
is Fredholm of index 0, for all $\zeta\neq0$. Moreover, denoting by $\sigma(\mathscr L_2)$  the spectrum of $\mathscr L_2$, we have that  $\sigma(\mathscr L_2)\cap \{\i\,\real\backslash\{0\}\}$
consists, at most, of a finite or countable number of eigenvalues, each of which is isolated and of finite (algebraic) multiplicity, having 0 as the only cluster point.  
\ET{5.1}
{\em Proof.} We begin to prove  that the operator $\mathscr K$ defined in \eqref{3.32}--\eqref{3.32_1} is compact. Let $\{\bsfW_k\}$ be a sequence bounded in $Z^{2,2}\times\real^3$. This implies, in particular, that there is $M>0$  independent of $k$ such that
\be
\|\bfw_k\|_{2,2}+|\bfsigma_k|\le M\,.
\eeq{5.9}
Employing \eqref{5.9} along with the compact embedding $W^{2,2}(\Omega)\subset W^{1,4}(\Omega_R)\cap L^\infty(\Omega_R)$, for all $R>1$,  we secure the existence of $(\bfw_*,\bfsigma_*)\in W^{2,2}(\Omega)\times\real^3$ and subsequences, again denoted by $\{\bfw_k,\bfsigma_k\}$ such that
\be
\bfw_k\to\bfw_*\ \ \mbox{strongly in $W^{1,4}(\Omega_R)\cap L^\infty(\Omega_R)$, for all $R>1$}\,;\ \ \bfsigma_k\to\bfsigma_*\ \ \mbox{in $\real^3$}\,. 
\eeq{JN}
In view of the linearity of $\mathscr K$, we may take, without loss,  $\bfw_*\equiv\bfsigma_*\equiv\0$. By H\"older inequality, we deduce 
$$\ba{ll}\medskip
\|\mathscr K(\bsfW_k)\|_2\le \lambda_0\left[\|\bsfu_0\|_4\|\bfw_k\|_{1,4,\Omega_R}+\|\bsfu_0\|_{4,\Omega^R}\|\bfw_k\|_{1,4}\right.\\
\hspace*{3.5cm}+\left.\|\nabla\bsfu_0\|_2(\|\bfw_k\|_{\infty,\Omega_R}+|\bfsigma_k|)+\|\nabla\bsfu_0\|_{2,\Omega^R}\|\bfw_k\|_{\infty}\right]\,.
\ea
$$
Therefore, letting $k\to\infty$ in this inequality, and using \eqref{u0w},  \eqref{5.9}, and  \eqref{JN} we infer, with a positive constant $C$ independent of $R$, that
$$
\lim_{k\to\infty}\|\mathscr K(\bsfW_k)\|_2\le C\,\left(\|\bsfu_0\|_{4,\Omega^R}+\|\nabla\bsfu_0\|_{2,\Omega^R}\right)\,.
$$
This, in turn, again by \eqref{u0w} and the arbitrariness of $R>1$ shows the desired property. From the property just showed and  \lemmref{4.1} it then follows that the operator
$$ 
\widehat{\mathscr L}_\zeta:=\mathscr L_2-\i\,\zeta\cali
$$ 
is Fredholm of index 0, for all $\zeta\neq0$. The theorem is then a consequence of well-known results (e.g. \cite[Theorem XVII.4.3]{GG}) provided we show that the null space of $\widehat{\mathscr L}_\zeta$ is trivial, for all sufficiently large $|\zeta|$. To this end, we observe that $\widehat{\mathscr L}_\zeta(\bsfW)=0$ means $\mathscr L_0(\bsfW)-\i\,\zeta\,\bsfW=-\mathscr P\mathscr K(\bsfW)$. Applying \eqref{4.2}, we thus get, in particular, the following inequality valid for all  $|\zeta|\ge 1$
$$
|\zeta|^{\frac12}\|\nabla\bfw\|_2+|\zeta|(\,\|\bfw\|_2+|\bfsigma|)\le c\,\|\mathscr K(\bsfW)\|_2\,,
$$ 
where $c$ is independent of $\zeta$. Also, from \eqref{3.32_1} and and H\"older inequality, we get
$$
\|\mathscr K(\bsfW)\|_2\le \lambda_0\, \|\bsfu_0\|_{1,\infty}(\|\bfw\|_{1,2}+|\bfsigma|)\,,
$$
and so, from the last two displayed equations and \theoref{exi} we deduce $\bfw\equiv\0$, provided we choose $|\zeta|$ larger than a suitable positive constant. This completes the proof of the theorem.
\par\hfill$\square$\par

\setcounter{equation}{0}
\section{The Linearized Time-Periodic Operator}\label{sec:time_per}
Objective of this section is to demonstrate some fundamental functional properties of the linear operator obtained by formally setting $(\bfsigma-\bfw)\cdot\nabla\bfw\equiv\0$ in equations \eqref{04}, in the class of periodic solutions with vanishing average.     
In order to do this, we need first to investigate the same question for the more special linear problem to which \eqref{04} reduces by taking $(\bfsigma-\bfw)\cdot\nabla\bfw\equiv\bfu_0\equiv\0$. 
\par
We begin with the following lemma.
\Bl Let $\zeta_0\in\real-\{0\}$, $\lambda_0\in\real$. Then, the boundary-value problems, $i\in\{1,2,3\}$, $k\in \mathbb Z\backslash\{0\}$,
\be
\ba{cc}\medskip\left.\ba{ll}\medskip
{\rm i}\,k\,\zeta_0\,\bfh_k^{(i)}-\lambda_0\partial_1\bfh_k^{(i)}=\Delta \bfh_k^{(i)}-\nabla p_k^{(i)}\\
\Div \bfh_k^{(i)}=0\ea\right\}\ \ \mbox{in $\Omega$}\\
\bfh_k^{(i)}=\bfe_i\ \ \mbox{at $\partial\Omega$}\,,
\ea
\eeq{3.2}
and $(\bfh_0^{(i)},p^{(i)}_0)\equiv (\0,0)$,
have a unique solutions $(\bfh_k^{(i)},p_k^{(i)})\in W^{2,2}(\Omega)\times D^{1,2}(\Omega)$. These solutions satisfy the estimates
\be
\|\bfh_k^{(i)}\|_2\le C\,\ \ \|\nabla\bfh_k^{(i)}\|_2\le C \,(|k|+1)^{\frac12}\,,\ \
\|D^2\bfh_k^{(i)}\|_2\le C \,(|k|+1)\,,
\eeq{3.4}
where $C$ is a constant independent of $k$. Moreover, for fixed $k$, let  $\mathbb K$ be the  matrix defined by ($j,i=1,2,3$):
\be
({\mathbb K})_{ji}=\IdS(\mathbb T(\bfh^{(i)}_k,p^{(i)}_k)\cdot\bfn)_j\,.
\eeq{Matrices}
Then, for any $\mu\in\real$, $\mathbb K+\i\,\mu\mathbb{I}$ is invertible. Finally,  for every $\bfalpha\in \mathbb C^3$, we have
\be
{\rm i}\,k\,\zeta_0\,\|\bfsf h_k\|_2^2+2\|\mathbb D(\bfsf h_k)\|_2^2-\lambda_0(\partial_1\bfsf h_k,\bfsf h_k^*)=\bfalpha^*\cdot\mathbb K\cdot\bfalpha
\eeq{3.A}
where $\bfsf h_k:= \alpha_i\bfh^{(i)}_k$. 
\EL{3.1}
{\em Proof.}  Since the proof is the same
for $i = 1, 2, 3$, we pick $i = 1$ and will omit the superscript. 
Let 
\be
\bsfe(x):=\curl(x_2\phi(|x|)\bfe_3) 
\eeq{ee}
where $\phi=\phi(|x|)$ is a smooth ``cut-off" function that is 1 in $\Omega_{\rho_1}$ and 0 in $\bar{\Omega^{\rho_{2}}}$, $R_*<\rho_1<\rho_2$. Clearly, $\bsfe$ is smooth in $\Omega$, with bounded support and, in addition, $\bsfe(x)=\bfe_1$ in a neighborhood of $\partial\Omega$ and $\Div\bsfe=0$ in $\Omega$. 
From \eqref{3.2} we then deduce that $\bfv_k:=\bfh_k-\bsfe$ solves the following boundary-value problem, for all $|k|\ge1$:
\be
\ba{cc}\medskip\left.\ba{ll}\medskip
{\rm i}\,k\,\zeta_0\,\bfv_k-\lambda_0\,\partial_1\bfv_k=\Delta \bfv_k-\nabla p_k^{(i)}+\lambda_0\partial_1\bsfe-{\rm i}\,k\,\zeta_0\,\bsfe+\Delta\bsfe\\
\Div \bfv_k=0\ea\right\}\ \ \mbox{in $\Omega$}\\
\bfv_k=\0\ \ \mbox{at $\partial\Omega$}\,.
\ea
\eeq{3.5}
We shall next show a number of a priori estimates for solutions $(\bfv_k,p_k)$  to \eqref{3.5} in their stated function class that once combined, for instance, with the classical Galerkin method, will produce the desired existence (and uniqueness) result.
For this, we dot-multiply both sides of \eqref{3.5}$_1$ by $\bfv^*$, and integrate by parts to get
\be
{\rm i}\,k\,\zeta_0\,\|\bfv_k\|_2^2-\lambda_0(\partial_1\bfv,\bfv_k^*)+\|\nabla\bfv_k\|_2^2=(\bfcalf_k,\bfv_k^*)\,,
\eeq{3.6}
where $\bfcalf_k:=\lambda_0\partial_1\bsfe-{\rm i}\,k\,\zeta_0\bsfe+\Delta\bsfe$. By the properties of $\bsfe$,
\be
\|\bfcalf_k\|_2\le c\,(|k|+1)
\eeq{3.7}
where, here and in the rest of the proof, $c$ denotes a generic (positive) constant depending, at most, on $\zeta_0,\lambda_0$ and $\Omega$. Since
\be
\Re\, (\partial_1\bfv_k,\bfv_k^*)=0\,,
\eeq{3.8}
by taking the real part of \eqref{3.6} and using \eqref{3.7} we infer
\be
\|\nabla\bfv_k\|_2^2\le c\,(|k|+1)\|\bfv_k\|_2\,.
\eeq{3.9}
Considering the imaginary part of \eqref{3.6} and employing \eqref{3.7}--\eqref{3.9} along with Schwarz inequality, we obtain 
$$
|k|\|\bfv_k\|_2\le c\,(\|\nabla\bfv_k\|_2+|k|+1)\le c(|k|+1)\|\bfv_k\|_2^{\frac12}\,,
$$
which gives
\be
\|\bfv_k\|_2\le c\,.
\eeq{3.10}
Taking into account that $\bfh_k=\bsfe+\bfv_k$, \eqref{3.10} proves \eqref{3.4}$_1$. Inequality \eqref{3.4}$_2$ then follows by using  \eqref{3.10} into \eqref{3.9}. Also, from classical estimates on the Stokes problem \cite[Lemma 1]{Hey} we get
$$
\|D^2\bfv_k\|_2\le c\,\big(\|\partial_1\bfv_k\|_2+\|\bfcalf_k\|_2+\|\nabla\bfv_k\|_2\big)
$$
and so, recalling that $|k|\ge 1$,  from the latter  inequality, \eqref{3.7}, \eqref{3.9} and \eqref{3.10} we show \eqref{3.4}$_3$.  Let $\bfalpha\in\mathbb C^3$, and, for fixed $k\neq0$,
 set\footnote{Summation over repeated indices.}
$$
\bfsf h_k:= \alpha_i\bfh^{(i)}_k\,,\ \ {\sf p}_k:=\alpha_i\,p^{(i)}_k\,.
$$
From \eqref{3.2}  we then find
\be
\ba{cc}\medskip\left.\ba{ll}\medskip
{\rm i}\,k\,\zeta_0\bfsf h_k-\lambda_0\partial_1\bfsf h_k=\Div\mathbb T(\bfsf h_k,{\sf p}_k)\\
\Div \bfsf h_k=0\ea\right\}\ \ \mbox{in $\Omega$}\\
\bfsf h_k=\bfalpha\,\ \mbox{at $\partial\Omega$}.
\ea
\eeq{3.2_ar}
Dot-multiplying both sides of \eqref{3.2_ar}$_1$ by $\bfsf h^*_k$ and integrating by parts over $\Omega$ we deduce
$$
\i\,k\,\zeta_0\|\bfsf h_k\|_2^2+\|\mathbb D(\bfsf h_k)\|_2^2-\lambda_0(\partial_1\bfsf h_k,\bfsf h^*_k)=\bfalpha^*\cdot\mathbb K\cdot\bfalpha\,.
$$
Now, suppose that there is $\hat{\bfalpha}\in\mathbb C^3$ such that $\mathbb K\cdot\hat{\bfalpha}=-\i\,\mu\,\hat{\bfalpha}$, for some $\mu\in\real$. Then from the previous relation we obtain
$$
\i\left(k\,\zeta_0\|\bfsf h_k\|_2^2+\mu\,|\hat{\bfalpha}|^2\right)-\lambda_0(\partial_1\bfsf h_k,\bfsf h^*_k)=2\|\mathbb D(\bfsf h_k)\|_2^2\,,
$$
which, in turn, recalling that \be
\Re\,(\partial_1\bfsf h_k,\bfsf h_k^*)=0\,, 
\eeq{re0}
allows us to we deduce $\bfsf h_k=\0$ in $W^{2,2}(\Omega)$. The latter implies $\hat{\bfalpha}=\0$ and thus shows the desired property for $\mathbb K$.  The proof of the lemma is completed.
\par\hfill$\square$\par
With the help of \lemmref{3.1}, we are now able to show the next one. 
\Bl Let $\lambda_0,\zeta_0$ be as in \lemmref{3.1}, and $(\omega_{\sf n}^2,\varpi)\in (0,\infty)\times (\real\backslash\{0\})$. Then, for any $\bsfF\in  L^2_\sharp$, the  problem 
\be\ba{cc}\medskip\left.\ba{ll}\medskip\zeta_0\partial_\tau\bsfw-\lambda_0\,\partial_1\bsfw=\Delta\bsfw-\nabla{\sf q}\\
\Div\bsfw=0\ea\right\}\,  \mbox{in $\Omega\times [0,2\pi]$}\,,\\ \medskip
\bfw=\dot{\bfxi}\ \ \mbox{at $\partial\Omega\times [0,2\pi]$}\,,\\ 
\zeta_0^2\ddot{\bfxi}+\omega_{\sf n}^2\bfxi+\varpi\Int{\partial\Omega}{}\mathbb T(\bsfw,{\sf q})\cdot\bfn=\bsfF\,,\, \ \mbox{in $[0,2\pi]$}\,,
\ea
\eeq{3.1}
has one and only one solution $\big(\bsfw,{\sf q},\bfxi\big)\in \mathcal W_\sharp^{2}
\times \mathcal P^{1,2}\times W^{2}_\sharp$. Such a  solution obeys the bound
\be
\|\bsfw\|_{\mathcal W_\sharp^{2}}+\|{\sf q}\|_{\mathcal P^{1,2}}+\|\bfxi\|_{W^{2}_\sharp}\le C\,\|\bsfF\|_{L^2_\sharp}\,,
\eeq{n}
with $C=C(\Omega,\lambda_0,\omega_{\sf n},\zeta_0)>0$.
\EL{3.2_0}
{\em Proof.} 
We formally expand $\bsfw$, ${\sf q}$, and  $\bfxi$, in Fourier series as follows:
\be
\bsfw(x,t)=\Sum{k\in\mathbb Z}{}\bsfw_k(x)\,{\rm e}^{\i k\,t}\,,\ \ {\sf q}(x,t)=\Sum{k\in\mathbb Z}{}{\sf q}_k(x)\,{\rm e}^{\i k\,t}\,,\ \
\bfxi(t)=\Sum{k\in\mathbb Z}{}\bfxi_k \,{\rm e}^{\i k\,t}\,,\ \ \bsfw_0\equiv\nabla{\sf q}_0\equiv\bfxi_0\equiv\0\,,
\eeq{Fou}
where $(\bsfw_k,{\sf q}_k,\bfxi_k)$ solve the problem
 ($k\neq0$)
\be\ba{cc}\medskip\left.\ba{ll}\medskip
\i\,k\,\zeta_0\,\bsfw_k-\lambda_0\partial_1\bsfw_k=\Delta \bsfw_k-\nabla{\sf q}_k\\
\Div\bsfw_k=0\ea\right\}\ \ \mbox{in $\Omega$}\\
\bsfw_k|_{\partial\Omega}={\rm i} k \bfxi_k\,,
\ea
\eeq{2.11}
with the further condition
\be
\left(-k^2\,\zeta_0^2+\omega_{\sf n}^2\right)\bfxi_k+\varpi\Int{\partial\Omega}{}\mathbb T(\bsfw_k,{\sf q}_k)\cdot\bfn=\bsfF_k\,,
\eeq{2.8_1}
where $\{\bsfF_k\}$ are Fourier coefficients of $\bsfF$ with $\bsfF_0\equiv\0$. 
For each fixed $k\in\mathbb Z-\{0\}$, a  solution to \eqref{2.11}--\eqref{2.8_1} is given by\Footnote{No summation over $k$.}
\be
\bsfw_k=\sum_{i=1}^3\i\,k\,\bfxi_{ki}\bfh_k^{(i)}\,,\ \ {\sf q}_k=\sum_{i=1}^3\i\,k\,\bfxi_{ki}p_k^{(i)}\,,
\eeq{1.20}
with $(\bfh_k^{(i)},p_k^{(i)})$ given in \lemmref{3.1}, and where $\bfxi_k$ solve the equations
\be
\left(-k^2\,\zeta_0^2+\omega_{\sf n}^2\right)\bfxi_k+\sum_{i=1}^3{\rm i}\,k\,\varpi\,\xi_{ki}\Int{\partial\Omega}{}\mathbb T(\bfh_k^{(i)},{p}^{(i)}_k)\cdot\bfn=\bsfF_k\,,
\eeq{chi}
which, with the notation of \lemmref{3.1}(ii) can be equivalently rewritten as
\be
\mathbb M\cdot\bfxi_k=\varpi\,\bsfF_k\,,\ \ \mathbb M:=(-k^2\,\zeta_0^2+\omega_{\sf n}^2)\mathbb I+{\rm i}\,k\,\varpi\,\mathbb K\,.
\eeq{1.21}
The matrix $\mathbb M$ is invertible for all $k\neq 0$. In fact, using \eqref{3.A} we show, for all $\bfalpha\in\mathbb C^3$,
$$
\bfalpha^*\cdot\mathbb M\cdot\bfalpha=(-k^2\,\zeta_0^2+\omega_{\sf n}^2)\,|\bfalpha|^2-k^2\zeta_0\,\varpi\,
\|\bfsf h_k\|_2^2-{\rm i}\lambda_0\varpi\,(\partial_1\bfsf h_k,\bfsf h_k^*)+2{\rm i}\,k\,\varpi\,\|\mathbb D(\bfsf h_k\|_2^2\,.
$$
Thus, assuming $\mathbb M\cdot \bfalpha=\0$, by \eqref{re0} it follows 
\be
\mathbb D(\bfsf h_k)\equiv 0\,.
\eeq{D} 
However, by the properties of $\bfh^{(i)}_k$ we obtain that $\bfsf h_k|_{\partial\Omega}=\bfalpha$, which by \lemmref{1.1}, the embedding $W^{1,2}(\Omega)\subset L^6(\Omega)$, and \eqref{D} implies $\bfalpha=\0$, namely,  $0$ is not an eigenvalue of $\mathbb M$. As a result, for the given $\bsfF_k$, equation  \eqref{1.21} has one and only one solution $\bfxi_k$.
If we now dot-multiply both sides of \eqref{1.21} by $\bfxi_k^*$ and use again \eqref{3.A} we deduce
$$\ba{ll}\medskip
(-k^2\,\zeta_0^2+\omega_{\sf n}^2)\,|\bfxi_k|^2-k^2\zeta_0\,\varpi\,
\|\xi_{ki}\bfh_k^{(i)}\|_2^2-{\rm i}\lambda_0\varpi\,(\partial_1(\xi_{ki}\bfh_k^{(i)}),(\xi_{ki}\bfh_k^{(i)})^*)\\
\hspace*{5.8cm}+2{\rm i}\,k\,\varpi\,\|\mathbb D(\xi_{ki}\bfh_k^{(i)})\|_2^2=(\bsfF_k,\bfxi_k^*)\,,
\ea
$$
which, in view of \eqref{re0}, furnishes
\be\ba{rl}\medskip
2k\,\varpi\|\mathbb D(\xi_{ki}\bfh_k^{(i)})\|_2^2&\!\!\!\!=\Im\,(\bsfF_k,\bfxi_k^*)\,,
\\
(-k^2\,\zeta_0^2+\omega_{\sf n}^2)\,|\bfxi_k|^2-k^2\zeta_0\varpi\,
\|\xi_{ki}\bfh_k^{(i)}\|_2^2-{\rm i}\lambda_0\varpi\,(\partial_1(\xi_{ki}\bfh_k^{(i)}),(\xi_{ki}\bfh_k^{(i)})^*)&\!\!\!\!=\Re\,(\bsfF_k,\bfxi_k^*)\,.
\ea
\eeq{2.17}
Recalling that $\xi_{ki}\bfh_k^{(i)}|_{\partial\Omega}=\bfxi_k$, by \eqref{2.17}$_1$, Schwarz inequality and \lemmref{1.1} we show the crucial estimate
\be 
|\bfxi_k|+\|\nabla(\xi_{ki}\bfh_k^{(i)})\|_2\le c\,|k|^{-1}\,|\bsfF_k|\,,\ \ |k|\ge1\,, 
\eeq{2.18}
where, here and in the following, $c$ denotes a generic positive constant independent of $k$.
Again by Schwarz inequality, from \eqref{2.17}$_2$ we get
$$
k^2\,(|\bfxi_k|^2+\|\xi_{ki}\bfh_k^{(i)}\|_2^2)\le c\,\left(|\bfxi_k|^2+\|\xi_{ki}\bfh_k^{(i)}\|_2\|\nabla(\xi_{ki}\bfh_k^{(i)})\|_2+|\bsfF_k||\bfxi_k|\right)
$$
from which, using \eqref{2.18} and Cauchy inequality
we deduce
$$
k^2\,(|\bfxi_k|^2+\|\xi_{ki}\bfh_k^{(i)}\|_2^2)\le c\,|k|^{-2}|\bsfF_k|^2\,,\ \ |k|\ge1\,,
$$
that allows us to conclude
\be
k^4|\bfxi_k|\le c\,|\bsfF_k|^2\,, \ \ |k|\ge1.
\eeq{2.190}
From \eqref{2.190} it immediately follows that
\be
\|\bfxi\|^2_{W^{2}_\sharp}=\sum_{|k|\ge1}(|k|^4+|k|^2+1)|\bfxi_k|^2\le c\sum_{|k|\ge1}|\bsfF_k|^2= c\|\bsfF\|_{L^2_\sharp}^2\,.
\eeq{1.46}
Moreover, from \eqref{1.20},
\eqref{1.46} and \eqref{3.4} we  infer
\be
\|\bsfw\|_{\calw_\sharp^2}^2=\!\sum_{|k|\ge1}\left[(|k|^2+1)\|\bsfw_k\|_2^2+\|\nabla\bsfw_k\|_2^2+\|D^2\bsfw_k\|_2^2\right]\le c\!\sum_{|k|\ge1}(|k|^4+|k|^2+1)|\bfsigma_k|^2\le c\,\|\bsfF\|_{L^2_\sharp}^2\,,
\eeq{1.49}
so that, combining \eqref{1.46},  \eqref{1.49}, and  \eqref{3.1}$_{1}$, we obtain
\be
\|\bsfw\|_{\calw_\sharp^2}+\|\bfxi\|_{W^{2}_\sharp}+\|\nabla{\sf q}\|_{L^2(D^{1,2})}\le c\,\|{\bsfF}\|_{L^2_\sharp}\,.
\eeq{2.22}
and the proof of existence is completed. The uniqueness property amounts to show that the problem
\be\ba{cc}\medskip\left.\ba{ll}\medskip
\zeta_0\partial_\tau\bsfw-\lambda_0\,\partial_1\bsfw=\Delta\bsfw-\nabla {\sf q}\\ 
\Div\bsfw=0\ea\right\}\ \ \mbox{in $\Omega\times[0,2\pi]$}\\ \medskip
\bsfw|_{\partial\Omega}=\dot{\bfxi}\,;\\\zeta_0^2 \ddot{\bfxi}+\omega_{\sf n}^2\bfxi+\varpi\Int{\partial\Omega}{}\mathbb T(\bsfw,{\sf q})\cdot\bfn=\0\
\ea
\eeq{3.31}
has only the zero solution in the specified function class. If we dot-multiply \eqref{3.31}$_1$ by $\bsfw$, integrate by parts over $\Omega$ and use \eqref{3.31}$_3$, we get
$$
\half\ode{}t(\zeta_0\|\bsfw(t)\|_2^2+\zeta_0^2|\dot{\bfxi}(t)|^2+\omega_{\sf n}^2|\bfxi(t)|^2)+2\|\mathbb D(\bsfw(t))\|_2^2=0\,.
$$
Integrating both sides of this equation from $0$ to $2\pi$ and employing the $2\pi$-periodicity of the solution, we easily obtain  $\|\mathbb D(\bsfw(t))\|_2\equiv 0$ which in turn immediately furnishes\footnote{Recall that $\bsfw(t)\in L^2(\Omega)$.}
$\bsfw\equiv\nabla{\sf q}\equiv\0$. The proof of the lemma is completed.\par\hfill$
\square$\par
\Br The invertibility of the matrix $\mathbb M$, defined in \eqref{1.21} is crucial to the resolution of \eqref{chi} and, therefore, to the result stated in the lemma. The remarkable fact is that this property holds  for {\em all} values of $\zeta_0,\omega_{\sf n}>0$, and $|k|>1$. In physical terms, this means that the possibility of a ``disruptive" resonance is always ruled out. Notice, however, that the same property is lost if $k^2\zeta_0=\omega_{\sf n}^2$, for some $k\neq 0$, and we allow $\varpi\to 0$, that is,  the density of the liquid vanishingly small compared to that of the body.
\ER{9.1}
\par
The next lemma proves well-posedness of the linear problem obtained from \eqref{04} by setting $(\bfsigma-\bfw)\cdot\nabla\bfw\equiv\bfu_0\equiv\0$, in the class of $2\pi$-periodic solutions with zero average. The crucial aspect of this result is that it does {\em not} impose {\em any} restriction on the natural frequency $\omega_{\sf n}$ of the spring. 
\Bl Let $\lambda_0,\zeta_0, \omega_{\sf n}^2, \varpi \in (0,\infty)$.\footnote{As is clear from the proof, we could  replace the assumption $\lambda_0>0$  with $\lambda\in\real$, but this is irrelevant to our purposes and also at odds with the physical meaning of the parameter $\lambda_0$. Likewise, we could, more generally,  assume $\zeta_0\in\real\backslash\{0\}$ which, again, would be immaterial; see \cite[Lemma 1.1]{You}.} Then, for any $(\bff,\bfF,\bfG)\in \call_\sharp^{2}\times L^2_\sharp\times W_\sharp^1$, the  problem 
\be\ba{cc}\medskip\left.\ba{ll}\medskip\zeta_0\partial_\tau\bfw-\lambda_0\partial_1\bfw=\Delta\bfw-\nabla{p}+\bff\\
\Div\bfw=0\ea\right\}\,  \mbox{in $\Omega\times [0,2\pi]$}\,,\\ \medskip
\bfw=\dot{\bfxi}-\bfG\ \ \mbox{at $\partial\Omega\times [0,2\pi]$}\,,\\ 
\,\zeta_0^2\ddot{\bfxi}+\omega_{\sf n}^2\bfxi+\varpi\Int{\partial\Omega}{}\mathbb T(\bfw,{p})\cdot\bfn=\bfF\,, \ \ \mbox{in $[0,2\pi]$}\,,
\ea
\eeq{3.1}
has one and only one solution $\big(\bfw,p,\bfxi\big)\in \mathcal W_\sharp^{2}
\times \mathcal P^{1,2}\times W^{2,2}_\sharp$. This solution satisfies the estimate
\be
\|\bfw\|_{\mathcal W_\sharp^{2}}+\|{p}\|_{\mathcal P^{1,2}}+\|\bfxi\|_{W^{2,2}}\le C\,\Big(\|\bff\|_{\call_\sharp^{2}}+\|\bfF\|_{L^2_\sharp}+\|\bfG\|_{W^1_\sharp}\Big)\,,
\eeq{n}
where $C=C(\Omega,\lambda_0,\zeta_0,\omega_{\sf n},\varpi)$.
\EL{3.2}
{\em Proof.} Let $\bfw=\bfz+\bfu$ where $\bfz$ and $\bfu$ satisfy the following problems
\be\ba{cc}\medskip\left.\ba{ll}\medskip
\zeta_0\partial_\tau\bfz-\lambda_0\partial_1\bfz=\Delta\bfz-\nabla {\sf r}+\bff\\ 
\Div\bfz=0\ea\right\}\ \ \mbox{in $\Omega\times [0,2\pi]$}\\
\bfz|_{\partial\Omega}=-\bfG
\ea
\eeq{2.3} 
and
\be\ba{cc}\medskip\left.\ba{ll}\medskip\zeta_0^2
\partial_\tau\bfu-\lambda_0\partial_1\bfu=\Delta\bfu-\nabla {\sf q}\\ 
\Div\bfu=0\ea\right\}\ \ \mbox{in $\Omega\times[0,2\pi]$}\\ \medskip
\bfu|_{\partial\Omega}=\dot{\bfxi}\,;\\ \medskip\zeta_0^2\ddot{\bfxi}+\omega_{\sf n}^2\bfxi+\varpi\Int{\partial\Omega}{}\bfT(\bfu,{\sf q})\cdot\bfn=\bfF-\varpi\Int{\partial\Omega}{}\bfT(\bfz,{\sf r})\cdot\bfn:=\bsfF\,,\ \ \mbox{in $[0,2\pi]$}\,.
\ea
\eeq{3.17}
Set 
\be   
\bfW(t):=x_3G_2(t)\bfe_1+x_1G_3(t)\bfe_2+x_2G_1(t)\bfe_3\,.
\eeq{W}
Clearly,
\be
\curl\bfW=\bfG(t)\,.
\eeq{W1}
Let $\phi(x)$ be the function defined in the beginning of the proof of \lemmref{3.1} and define
$$
\bfw(x,t):=\curl\big(\phi(x)\bfW(t)\big)\,.
$$
In view of \eqref{W1} we deduce
\be
\bfw(x,t)=\phi(x)\bfG(t)-\bfW\times\nabla\phi(x)
\eeq{03_1}
so that $\bfw$ is a $2\pi$-periodic solenoidal vector function 
that is equal to $\bfG(t)$ for $|x|\ge \rho_2$ and equal to 0 for $|x|\le \rho_1$. Therefore, from \eqref{2.3} we deduce that the field
\be
\bsfz(x,t):=\bfz(x,t)-\bfw(x,t)\,,
\eeq{zeta}
obeys the following problem
\be\ba{cc}\medskip\left.\ba{ll}\medskip
\zeta_0\partial_\tau\bsfz-\lambda_0\partial_1\bsfz=\Delta\bsfz-\nabla {\sf r}+\bsff\\ 
\Div\bsfz=0\ea\right\}\ \ \mbox{in $\Omega\times [0,2\pi]$}\\
\bsfz|_{\partial\Omega}=\0\,,
\ea
\eeq{2.3_1}
where
\be
\bsff:=\bff-\zeta_0\partial_\tau\bfw+\lambda_0\partial_1\bfw+\Delta\bfw\,.
\eeq{EDB}
From \eqref{W}--\eqref{03_1},  \eqref{EDB} and the assumption on $\bfG$ it follows that $\bsff\in \call^2_\sharp$ and that
\be
\|\bsff\|_{\call_\sharp^2}\le c\,(\|\bff\|_{\call_\sharp^2}+\|\bfG\|_{W^2_\sharp})\,.
\eeq{mds}
Employing \cite[Theorem 12]{GaMaH}, we then deduce  that there exists a unique solution $(\bsfz,{\sf r})\in  \mathcal W_\sharp^{2}\times\calp^{1}$   that, in addition, obeys the inequality
$$
\|\bsfz\|_{\mathcal W_\sharp^{2}}+\|{\sf r}\|_{\mathcal P^{1,2}}\le c\,\|\bsff\|_{\mathcal L_\sharp^{2}}\,.
$$
The latter, in combination with \eqref{mds} and \eqref{zeta}, allows us to conclude $\bfz\in\calw_\sharp^2$ and
\be
\|\bfz\|_{\mathcal W_\sharp^{2}}+\|{\sf r}\|_{\mathcal P^{1,2}}\le c\,(\|\bff\|_{\mathcal L_\sharp^{2}}+\|\bfG\|_{W^2_\sharp})\,.
\eeq{3.18}
Now, by the trace theorem\Footnote{Possibly, by modifying $\sf r$ by adding to it a suitable function of time.} and \eqref{3.18} we get
$$
\|\Int{\partial\Omega}{}\mathbb T(\bfz,{\sf \sf r})\cdot\bfn\|_{L^2}\le c\,\left(\|\bfz\|_{\mathcal W^{2}_\sharp}+\|{\sf \sf r}\|_{\mathcal P^{1,2}}\right)\le c\,(\|\bff\|_{\mathcal L_\sharp^{2}}+\|\bfG\|_{W^2_\sharp})\,,
$$
so that  $\bsfF$ is in $L^2_\sharp(0,2\pi)$ and satisfies
\be
\|\bsfF\|_{L^2_\sharp}\le c(\|\bff\|_{\call^2_\sharp}+\|\bfF\|_{L^2_\sharp}+\|\bfG\|_{W^2_\sharp}))
\,.
\eeq{3.19} 
Thus, from \lemmref{3.2_0} it follows that there is one and only one solution $(\bfu,{\sf q},\bfxi)\in \calw_\sharp^2\times \calp^{1,2}\times W^{2,2}_\sharp$ to \eqref{3.17} that, in addition, satisfies the estimate
$$
\|\bfu\|_{\mathcal W_\sharp^{2}}+\|{\sf q}\|_{\mathcal P^{1,2}}+\|\bfxi\|_{W^{2}_\sharp}\le c\,\|\bsfF\|_{L^2_\sharp}\,.
$$
As a result, combining the latter with \eqref{3.19} and \eqref{3.18}, we
 complete the existence proof. The stated uniqueness property amounts to show that every solution $(\bfw,{\sf p},\bfxi)$ to \eqref{3.31} in the stated function class vanishes identically. However, this has already been showed at the end of the proof of \lemmref{3.2_0}.  
\par\hfill$\square$\par
We are now in a position to define a suitable linearized operator and to state its main functional properties. To this end, let
$$
\mathscr Q_0:\bsfW\in {\W}\times W^2_\sharp\mapsto \mathscr Q_0(\bsfW)=\zeta_0\partial_t\bsfW+\mathscr L_0(\bsfW)\in \M\times W^1_\sharp\,,
$$
where $\mathscr L_0$ is defined in \eqref{L0}. Then, given $\bsfF:=(\bff,\bfG)\in \M\times W^1_\sharp$,  
the operator equation
\be
\mathscr Q_0(\bsfW)=\bsfF\,.
\eeq{chiappa}
is {\em equivalent} to the following problem 
\be \ba{cc}\medskip\left.\ba{ll}\medskip
\zeta_0\,\partial_\tau\bfw-\lambda_0\partial_1\nabla\bfw-{\Delta}\bfw=\nabla \phi+\bff\\ 
\Div\bfw=0\ea\right\}\ \ \mbox{in $\Omega\times[0,2\pi]$}\,,\\ \medskip
\bfw=\bfsigma\ \ \mbox{at $\partial\Omega\times[0,2\pi]$}\,,\\ \medskip
\dot{\bfsigma}+\omega^2_{\sf n}\bfxi^2+\varpi\Int{\partial\Omega}{}\mathbb T(\bfw,{\phi})\cdot\bfn=\bfF \ \ \mbox{in $[0,2\pi]$}\,,\\
\dot{\bfxi}-\bfsigma=\bfG\,,\ \ \mbox{in $[0,2\pi]$}\,.
\ea
\eeq{ciapa}
where $\bfF:=\bff|_{\Omega_0}$.
Thus, 
from \lemmref{3.2} we infer the following important result.
\Bl For any given $\lambda_0,\omega_{\sf n}^2, \zeta_0,\varpi\in (0,\infty)$, the operator $\mathscr Q_0$ is a homeomorphism.
\EL{Hom}
\par
Next, define
$$
\mathscr Q:\bsfW\in {\W}\times W^2_\sharp\mapsto \mathscr Q(\bsfW)=\zeta_0\partial_t\bsfW+\mathscr L_2(\bsfW)\in \M\times W^1_\sharp\,,
$$
with $\mathscr L_2$ given in \eqref{L2}.
\lemmref{Hom} allows us to prove the following theorem, representing the main result of this section.
\Bt For any given $\lambda_0,\omega_{\sf n}^2, \zeta_0,\varpi\in (0,\infty)$,
the operator $\mathscr Q$ is Fredholm of index 0.
\ET{3.1}
{\em Proof.} Since $\mathscr Q=\mathscr Q_0+\mathscr P\mathscr K$,  by \lemmref{Hom},  the stated property will follow provided we show that the map
$$
\mathscr K :\bsfW\in\W\mapsto \mathscr K(\bsfW)\in \M
$$
is compact. Let $\{\bsfW_k\}$ be a bounded sequence in $\W$. This means  that there is $M>0$ independent of $k$ such that 
\be
\|\bfw_k\|_{\calw_\sharp^2}+\|\bfsigma_k\|_{W_\sharp^{1}}\le M\,,
\eeq{gatto1}
with $\bfsigma_k:=\bfw_k|_{\Omega_0}$. 
We may then select sequences (again denoted by $\{\bfw_k,\bfsigma_k\}$) and find $(\bfw_*,\bfsigma_*)\in \mathcal W^{2}_\sharp\times W^{1,2}_\sharp$ such that
\be
\bfw_k\to{\bfw_*} \ \, \mbox{weakly in $\W$\,;}\ \ \bfsigma_k\to{\bfsigma_*}\,  \  \mbox{strongly in $L^{\infty}(0,2\pi)$.}
\eeq{2.18}
Due to the linearity of $\mathscr K$, without loss of generality  we may take $\bfw_*\equiv\bfsigma_*\equiv\0$, so that we must  show that
\be
\lim_{k\to\infty}\int_0^T\|\mathscr K(\bsfW_k)\|_{2,\Omega}^2=0\,.
\eeq{gatto2}
From \eqref{2.18},  the compact embeddings $W^{2,2}(\Omega)\subset W^{1,4}(\Omega_R)$ for all $R>R_*$, and Lions-Aubin lemma we infer
\be
\int_{0}^{2\pi}\left(\|\bfw_k(t)\|_{4,\Omega_R}^2+\|\nabla\bfw_k(t)\|_{4,\Omega_R}^2\right){\rm d}t\to 0\ \ \mbox{as $k\to\infty$, for all $R>R_*$\,.}
\eeq{2.19}
Further,   
$$
\int_{0}^{2\pi} \|\bsfu_0\cdot\nabla\bfw_k(t)\|_{2,\Omega}^2
\le \|\bsfu_0\|_4^2\int_{0}^{2\pi} \|\nabla\bfw_k(t)\|_{4,\Omega_R}^2+ \|\bsfu_0\|_{4,\Omega^R}^2\int_{-\pi}^{\pi}\|\nabla\bfw_k(t)\|_{4,\Omega}^2\,,
$$
which, by  \propref{3.1}, the embedding $W^{2,2}(\Omega)\subset W^{1,4}(\Omega)$,  \eqref{gatto1}, \eqref{2.19}  and the arbitrariness of $R$ furnishes
\be
\lim_{k\to\infty}\int_{0}^{2\pi} \|\bsfu_0\cdot\nabla\bfw_k(\tau)\|_{2}^2=0\,.
\eeq{2.21}
Likewise, 
$$
\Int{0}{2\pi} \|\bfw_k(t)\cdot\nabla\bsfu_0\|_{2,\Omega}^2
\le \|\nabla\bsfu_0\|_{4}^2\Int{0}{2\pi} \|\bfw_k(t)\|_{4,\Omega_R}^2 + 
\|\nabla\bsfu_0\|_{4,\Omega^R}^2\Int{0}{\pi}\|\bfw_k(t)\|_{4,\Omega}^2\,,
$$
so that, by  \lemmref{Xemb},  \eqref{gatto1}, and \eqref{2.19} we deduce, as before,
\be
\lim_{k\to\infty}\int_{0}^{2\pi} \|\bfw_k(t)\cdot\nabla\bsfu_0\|_{2,\Omega}^2=0\,.
\eeq{2.22}
Finally,
$$
\int_{0}^{2\pi}\|\bfsigma_k\cdot\nabla\bsfu_0\|_{2,\Omega}^2\le {2\pi}\, \|\bfsigma_k\|_{L^\infty(0,2\pi)}^2\|\nabla\bsfu_0\|_2^2\,,
$$
which, by \lemmref{Xemb} and \eqref{2.18}$_2$ furnishes
\be
\lim_{k\to\infty}\int_{0}^{2\pi}\|\bfsigma_k\cdot\nabla\bsfu_0\|_{2,\Omega}^2=0
\eeq{3.38}
Combining \eqref{2.21}--\eqref{3.38} we thus arrive at \eqref{gatto2},
which completes the proof of the theorem.\par
\hfill$\square$
\setcounter{equation}{0}
\section{Reformulation of the Problem in Banach Spaces}\label{sec:rif}
In this section we shall write the bifurcation problem, stated in Section \ref{sec:Form}, as a system of operator equations in an appropriate functional setting that, thanks to the results proved in previous sections,  will enable us to apply \theoref{3.1_ar}. The first step is to split \eqref{04} into its averaged and oscillatory components. Set
$$
\bfw:=\bar{\bfw}+(\bfw-\bar{\bfw}):=\bsfu+\bsfw\,,\ \ \bfeta=\bar{\bfeta}+(\bfeta-\bar{\bfeta}):=\bar{\bfeta}+\bfxi\,,\ \ {\sf p}=\bar{\sf p}+({\sf p}-\bar{\sf p}):=\bar{\sf p}+{\sf q}\,,\ \ \mu:=\lambda-\lambda_0\,,
$$
and 
$$\bsfu_0:=\bfu_0(x;\lambda_0)\,, \ \ \tilde{\bfu}_0:=\bfu_0(x;\mu+\lambda_0)-\bsfu_0.
$$ 
From \eqref{04} we thus get
\be\ba{cc}\medskip\left.\ba{ll}\medskip
-\lambda_0(\partial_1\bsfu-\bsfu_0\cdot\nabla\bsfu-\bsfu\cdot\nabla\bsfu_0)-\Delta\bsfu+\nabla\bar{\sf p}=\bfN_1(\bsfu,\bsfw,\mu)\\
\Div\bsfu=0\ea\right\}\ \ \mbox{in $\Omega$\,,}\\ \medskip
\bsfu=\0\ \ \mbox{at $\partial\Omega$}\,,\\
\omega_{\sf n}^2\bar{\bfeta}+\varpi\Int{\partial\Omega}{}\mathbb T(\bsfu,\bar{\sf p})\cdot\bfn=\0\,,
\ea
\eeq{7.1}
where
\be
\bfN_1:=\left\{\ba{ll}\ms-\mu(\partial_1\bsfu-\bsfu_0\cdot\nabla\bsfu-\bsfu\cdot\nabla\bsfu_0)-(\mu+\lambda_0)\left[\tilde{\bfu}_0\cdot\nabla\bsfu+\bsfu\cdot\nabla\tilde{\bfu}_0+\bsfu\cdot\nabla\bsfu+\bar{(\bfsigma-\bsfw)\cdot\nabla\bsfw}\right]\ \mbox{in $\Omega$}\\
\0\ \mbox{in $\Omega_0$}
\,,
\ea\right.
\eeq{7.2}
and, with the time-scaling $\tau:=\zeta\,t$,
\be\ba{cc}\medskip\left.\ba{ll}\medskip
\zeta\partial_\tau\bsfw-\lambda_0(\partial_1\bsfw-\bsfu_0\cdot\nabla\bsfw+(\bfsigma-\bsfw)\cdot\nabla\bsfu_0)-\Delta\bsfw+\nabla{\sf q}=\bfN_2(\bsfu,\bsfw,\mu)\\
\Div\bsfw=0\ea\right\}\ \ \mbox{in $\Omega\times [0,2\pi]$}\\ \medskip
\bsfw=\bfsigma\ \ \mbox{at $\partial\Omega\times[0,2\pi]$}\\
\zeta\dot{\bfsigma}+\omega_{\sf n}^2\bfxi+\varpi\Int{\partial\Omega}{}\mathbb T(\bsfw,{\sf q})\cdot\bfn=\0\,,\ \ \zeta\dot{\bfxi}-\bfsigma=\0\ \ \mbox{in $[0,2\pi]$}\,,
\ea
\eeq{7.3}
where
\be
\bfN_2:=\left\{\ba{ll}\medskip\ba{ll}\medskip\mu(\partial_1\bsfw-\bsfu_0\cdot\nabla\bsfw+(\bfsigma-\bsfw)\cdot\nabla\bsfu_0)
+(\mu+\lambda_0)\Big[\tilde{\bfu}_0\cdot\nabla\bsfw+(\bsfw-\bfsigma)\cdot\nabla\tilde{\bfu}_0+\bar{(\bfsigma-\bsfw)\cdot\nabla\bsfw}\\
\hspace*{5.5cm}+{(\bfsigma-\bsfw)\cdot\nabla\bsfw}+\bsfu\cdot\nabla\bsfw+(\bfsigma-\bsfw)\cdot\nabla\bsfu\Big]\ \mbox{in $\Omega$}\ea\\
\0 \ \mbox{in $\Omega_0$}
\ea\right.
\eeq{7.4}
\Bl
Let $\bfu_0=\bfu_0(\lambda)$, be the velocity field of the solution determined in \theoref{exi},  corresponding to $\lambda>0$, and let   $(\bsfu,\bsfw,\mu)\in X^2(\Omega)\times {\sf W}^2_\sharp\times\real$. Then, the following properties hold:
\be
\bfN_1(\bsfu,\bsfw,\mu)\in \mathcal Y(\Omega)\,;\ \ \bfN_2(\bsfu,\bsfw,\mu)\in {\sf L}^2_\sharp\,.
\eeq{7.5}
\EL{7.1}
{\em Proof.} Taking into account   \eqref{7.2} and the proof of \lemmref{L1}, it is easy to verify that to prove the first property in \eqref{7.5} it suffices to show
\be
\bar{(\bfsigma-\bsfw)\cdot\nabla\bsfw}\in Y(\Omega)\equiv \cald_0^{-1,2}(\Omega)\cap L^2(\Omega)\,.
\eeq{7.6}
By Schwarz inequality, integration by parts and elementary embedding theorems we get, for arbitrary $\bfphi\in\cald_0^{1,2}(\Omega)$, 
$$
|\big(\bar{(\bfsigma-\bsfw)\cdot\nabla\bsfw},\bfphi)|\le (\|\bfsigma\|_{L^4_\sharp}+\|\bfw\|_{L^4(L^4)})\|\bfw\|_{L^4(L^4)}\|\nabla\bfphi\|_2\le c\, \|\bfw\|_{\mbox{\scriptsize $\W$}}^2\|\nabla\bfphi\|_2\,.
$$
By a similar argument,
$$
\|\bar{(\bfsigma-\bsfw)\cdot\nabla\bsfw}\|_2\le (\|\bfsigma\|_{L^4_\sharp}+\|\bfw\|_{L^4(L^4)})\|\nabla\bfw\|_{L^4(L^4)} \le c\, \|\bfw\|_{\mbox{\scriptsize $\W$}}^2\,,
$$
which completes the proof of \eqref{7.6}. We next observe that, by \theoref{exi} we have $\bsfu_0,\tilde{\bfu_0}\in W^{1,\infty}(\Omega)$, whereas $\bsfu\in L^\infty(\Omega)$, by \lemmref{Xemb}. Thus, from \eqref{7.4} and also bearing in mind    the proof of the first property in \eqref{7.5} just given, it is easy to convince oneself that in order to show the second property in \eqref{7.5}, we can limit ourselves to prove $(\bfsigma-\bsfw)\cdot\nabla\bsfu\in L^2(L^2(\Omega))$. However, the latter follows at once  because  $\bfsigma\in W_\sharp^1\subset L^\infty(0,2\pi)$, $\bsfw\in L^2(W^{2,2})\subset L^2(L^\infty)$ and $\bsfu\in X^2(\Omega)\subset D^{1,2}(\Omega)$. 
\par\hfill$\square$\par
Setting $\mathscr N_i:=\mathscr P\bfN_i$, $i=1,2$, and recalling the definition of the operators $\mathscr L_1$ and $\mathscr L_2$ given in \eqref{2.2} and  \eqref{L2} respectively, we deduce that equations \eqref{7.1} and \eqref{7.3}
can equivalently be rewritten as follows
\be\ba{ll}\medskip
\mathscr L_1(\bsfU)=\mathscr N_1(\bsfU,\bsfW,\mu)\ \ \mbox{in $\mathcal Y$}\\
\zeta\partial_\tau\bsfW+\mathscr L_2(\bsfW)=\mathscr N_2(\bsfU,\bsfW,\mu)\ \ \mbox{in ${\sf L}^2_\sharp$} \,,
\ea
\eeq{7.7}
where $\bsfU:=(\bsfu,\bar{\bfeta})\in X^2(\Omega)\times\real^3$ and $\bsfW:=(\bsfw,\bfxi)\in \W\times W^1_\sharp$.
\setcounter{equation}{0}
\section{A Time-Periodic Bifurcation Theorem}\label{sec:Bif}  
Our goal is to apply \theoref{3.1_ar} to \eqref{7.7} and derive corresponding sufficient conditions for the occurrence of time-periodic bifurcation. To this end, we shall   consider separately the assumptions (H1)--(H4) made in that theorem and  formulate them  appropriately for the case at hand. 
\subsection*{\sf Formulation of  Assumption ({H1})} By \lemmref{fredo}, the operator $\mathscr L_1$ is Fredholm of index 0. Therefore, the assumption (H1) on $\mathscr L_1$ being a homeomorphism is satisfied if ${\sf N}[\mathscr L_1]=\{\0\}$, that is
\tag{H1$^\prime$}
\be
\mathscr L_1(\bsfU)=\0\ \ \Longrightarrow\ \ \bsfU=\0\,.
\eeq{H1'}
Recalling the definition of $\mathscr L_1$, we find that \eqref{H1'} is equivalent to the following request: If $(\bsfu,{\sf r},\bfchi)\in X^2(\Omega)\times D^{1,2}(\Omega)\times\real^3$ is a solution to the problem
$$\ba{cc}\medskip\left.\ba{ll}\medskip
-\lambda_0(\partial_1\bsfu-\bsfu_0\cdot\nabla\bsfu-\bsfu\cdot\nabla\bsfu_0)-\Delta\bsfu+\nabla {\sf r}=\0\\
\Div\bsfu=0\ea\right\}\ \ \mbox{in $\Omega$\,,}\\ \medskip
\bsfu=\0\ \ \mbox{at $\partial\Omega$}\,,\\
\omega_{\sf n}^2\bar{\bfchi}+\varpi\Int{\partial\Omega}{}\mathbb T(\bsfu,\bar{\sf p})\cdot\bfn=\0\,,
\ea
$$
then, necessarily, $\bsfu\equiv\nabla {\sf r}\equiv\bfchi\equiv\0$. According to \theoref{2.1}, this  means that the steady-state solution ${\sf s}(\lambda)$ in \eqref{sfaco} is unique for all $\lambda\in U(\lambda_0)$. 
\subsection*{\sf Formulation of  Assumption ({H2})} According to \theoref{5.1} the operator $\mathscr L_2-{\rm i}\,\zeta\,\cali$ is Fredholm of index 0 and, moreover $\Sigma:=\sigma(\mathscr L_2)\cap \{{\rm i}\real\}$ is constituted only by eigenvalues of finite algebraic multiplicity (a.m.). Therefore, the assumption (H2) can be formulated as follows:
\tag{H2$^\prime$}
\be
\mbox{There is $\nu_0:={\rm i}\zeta_0\in\Sigma$ with a.m.=1, and $k\nu_0\not\in \Sigma$, for all $|k|>1$.} 
\eeq{H2'}
\setcounter{equation}{0}\renewcommand{\theequation}{\arabic{section}.\arabic{equation}}Taking into account the definition of $\mathscr L_2$ given in \eqref{L2}, we show that  assumption \eqref{H2'} implies, in particular, that the eigenvalue problem
\be\ba{cc}\medskip\left.\ba{ll}\medskip
-\i\zeta_0\bfw-\lambda_0\,[\partial_1\bfw-\bfu_0\cdot\nabla\bfw+(-i\zeta_0{\bfxi}-\bfw)\cdot\nabla\bfu_0]=\Delta\bfw-\nabla {\sf p}\\
\Div\bfw=0\ea\right\}\ \ \mbox{in $\Omega$}\,,\\ \medskip
\bfw(x)=-\i\zeta_0{{\bfxi}}\,, \ \mbox{ $x\in\partial\Omega$}\,,
\\
(-\zeta_0^2+\omega_{\sf n}^2)\bfxi+\varpi\Int{\partial\Omega}{} \mathbb T(\bfw,{\sf p})\cdot\bfn=\0\,,
\ea
\eeq{8.1}
has a corresponding one-dimensional eigenspace $(\bfw,\bfxi)\in  Z^{2,2}\times\real^3$.
\subsection*{\sf Formulation of  Assumption ({H3})} In view of \theoref{3.1} this assumption is automatically satisfied in our case.
\subsection*{\sf Formulation of  Assumption ({H4})} The assumption \eqref{H1'} combined with \theoref{2.1} entails that the map
$$
\mu\in U(0)\subset\real\mapsto \bfu_0(\mu+\lambda_0)\in X^2(\Omega)
$$
with $\bfu_0(\mu+\lambda_0)$ velocity field of the solution determined in \theoref{exi} and corresponding to $\mu+\lambda_0$, is analytic. In addition,
the nonlinear operators $\mathscr N_i$, $i = 1, 2,$ are
(at most) quadratic in $(\bsfu,\bsfw)$  and then, by \lemmref{7.1}, analytic in those variables. Therefore, we conclude that also assumption (H4) is satisfied in the case at hand.
\smallskip\par
Our final comment, before stating the bifurcation theorem, regards the assumption \eqref{nupr} and its formulation in the context of our problem. To this end, it is easy to show from \eqref{7.4} that, in the case at hand, we have
$$
S_{011}=\partial_1\bsfw-\bsfu_0\cdot\nabla\bsfw+(\bfsigma-\bsfw)\cdot\nabla\bsfu_0
+\lambda_0\big[\tilde{\bfu}^\prime_0(0)\cdot\nabla\bsfw+(\bsfw-\bfsigma)\cdot\nabla\tilde{\bfu}_0^\prime(0)\big]
$$
where the prime denotes differentiation with respect to $\mu$. Therefore, denoting by $\nu = \nu(\mu)$ the
eigenvalues of $\mathscr L_2 + \mu S_{011}$, we may apply \cite[Proposition 79.15 and Corollary 79.16]{Z1} to show  that  the map $$\mu\in U(0)\mapsto \nu(\mu)\in\comp$$ is well defined and of class $C^\infty$. 
We are now in a position to state our bifurcation result.
\Bt Let $\omega_{n}^2$, $\varpi>0$. Suppose there exists $(\lambda_0,\bsfu_0,\zeta_0)$ such that  assumptions \eqref{H1'} and \eqref{H2'} hold and, moreover,
$$ 
\Re[\nu^\prime(0)]\neq 0\,.
$$
Then, the following properties are valid. \smallskip\\
{\rm (a)} {\rm Existence.} There is an analytic family
\be
\big({\bsfu}(\varepsilon),\bsfw(\varepsilon),\bfxi(\varepsilon),\zeta(\varepsilon),\mu(\varepsilon)\big)\in X^2(\Omega)\times {\sf W}^2_\sharp \times W^1_\sharp\times \real_+\times\real
\eeq{fam1}
of solutions to \eqref{7.1}, \eqref{7.3}, for all  $\varepsilon$ in a neighborhood $\mathcal I(0)$ of\, $0\in\real$. Moreover, let $(\bsfw_0,\bfxi_0)\in Z^{2,2}_{\mbox{\tiny $\comp$}}\times \comp^3$ be a normalized eigenfunction of the operator $\mathscr L_2$ corresponding to the eigenvalue $\i\zeta_0$, and set $(\bsfw_1,\bfxi_1):=\Re[(\bsfw_0,\bfxi_0)\,{\rm e}^{-\i\tau}]$. Then  
\be
\big({\bsfu}(\varepsilon),\bsfw(\varepsilon)-\varepsilon\,\bsfw_1,\bfxi(\varepsilon)-\varepsilon\,\bfxi_1,\zeta(\varepsilon),\mu(\varepsilon)\big)\to (0,0,\zeta_0,0)\ \ \mbox{as $\varepsilon\to 0$}\,.
\eeq{Ar.101}
\par\noindent
{\rm (a)} {\rm Uniqueness.}
There is a neighborhood  $$\calu(\0,\0,\0,\zeta_0,0)\subset X^2(\Omega)\times {\sf W}^2_\sharp\times W^1_\sharp\times \real_+\times \real$$ such that every (nontrivial) $2\pi$-periodic solution to \eqref{7.1}, \eqref{7.2},  lying in $\calu$ must coincide, up to a phase shift, with a member of the family \eqref{fam1}.
\smallskip\par\noindent
{\rm (a)} {\rm Parity.}  The functions $\zeta(\varepsilon)$ and $\mu(\varepsilon)$ are even:
$$
\zeta(\varepsilon)=\zeta(-\varepsilon)\,,\ \ \mu(\varepsilon)=\mu(-\varepsilon)\,,\ \ \mbox{for all $\varepsilon\in\cali(0)$\,.} 
$$
Consequently, the bifurcation due to these solutions is either subcritical or supercritical, a two-sided bifurcation being excluded.\footnote{Unless $\mu\equiv 0$.}
\ET{8.1}
\setcounter{equation}{0}\renewcommand{\theequation}{A.\arabic{equation}}
\section*{Appendix} 
\addcontentsline{toc}{section}{Appendix}
We shall give a proof of \theoref{3.1_ar}. We begin with some preparatory results.
Let $v_0$ be a normalized eigenvector of $L_2$ corresponding to the eigenvalue $\nu_0$, and set
\be
v_1:=\Re[v_0\,{\rm e}^{-{\rm i}\,\tau}]\,,\ \ v_2:=\Im[v_0\,{\rm e}^{-{\rm i}\,\tau}]\,. 
\eeq{vee}    
\begin{lemmaA} {\sl Under the assumption {\rm (H2)}, we have ${\rm dim}\,{\sf N}\,[\mathscr Q]=2$, and $\{v_1,v_2\}$ is a basis in ${\sf N}\,[\mathscr Q]$.}
\label{3.1_00}
\end{lemmaA}
{\em Proof.} Clearly, $\mathcal S:= {\rm span}\,\{v_1,v_2\}\subseteq {\sf N}[\mathscr Q]$. Conversely, take $w\in {\sf N}[\mathscr Q]$, and expand it in Fourier series
$$
w=\sum_{\ell=-\infty}^\infty w_\ell\,{\rm e}^{-{\rm i}\,\ell\,\tau}\,;\ w_\ell:=\frac1{2\pi}\int_{-\pi}^\pi w(\tau)\,{\rm e}^{{\rm i}\,\ell\,\tau}\,d\tau\,,\ \ w_0\equiv\bar w=0.
$$
Obviously, $w_\ell\in \calw_{\comps}\equiv {\sf D}_{\comps}[L_2]$. From $\mathscr Q(w)=0$ we deduce 
$$
-\ell\,\nu_0\,w_\ell+L_2(w_\ell)=0\,,\ \ w_\ell\in{\sf D}_{\comps}[L_2]\,,\ \ \ell\in\mathbb Z, 
$$ 
which, by (H2) and the fact that $w_0=0$, implies $w_\ell=0$ for all $\ell\in \mathbb Z\backslash\{\pm 1\}$. Thus, recalling that $\nu_0$ is simple, we infer  $w\in\cals$ and the lemma follows.\QED
\smallskip\par
Denote by $\langle\cdot,\cdot\rangle$ the scalar product in $\calz$ and set
$$
(w_1|w_2):=\int_{-\pi}^{\pi}\langle w_1(s),w_2(s)\rangle\,{\rm d}s\,,\ \ w_1,w_2\in \calz_{2\pi,0}\,.
$$
Let $L^\dagger_2$ be the adjoint of  $L_2$. Since $\nu_0$ is simple and $L_2-\nu_0I$ is Fredholm of index 0 (by (H2)), from classical results (e.g. \cite[Section 8.4]{Z}), it follows that there exists at least one  element  $v_0^\dagger\in{\sf N}_{\comps}[L_2^\dagger-\nu_0\,I]$ such that $\langle v_0^\dagger,v_0\rangle\neq 0$. Without loss, we may take 
\be
\langle v_0^\dagger,v_0\rangle=\pi^{-1}\,.
\eeq{3.6_00}
We then define
$$
v_1^\dagger:=\Re[v_0^\dagger\,{\rm e}^{{\rm i}\,\tau}]\,,\ \ v_2^\dagger:=\Im[v_0^\dagger\,{\rm e}^{{\rm i}\,\tau}]\,, 
$$
and observe that,
by \eqref{vee} and \eqref{3.6_00}, 
\be\ba{ll}\medskip 
(v_1|v_1^\dagger)=(v_2|v_2^\dagger)=1\,,\ \ (v_2|v_1^\dagger)=(v_1|v_2^\dagger)=0\,,\\
( (v_1)_\tau|v_1^\dagger)=0\,,\ \ ((v_1)_\tau|v_2^\dagger)=-1\,. 
\ea
\eeq{3.7_00}
Set
$$
\hat{\calz}_{2\pi,0}=\big\{w\in {\calz}_{2\pi,0}: \ (w|v_1^\dagger)=(w|v_2^\dagger)=0\big\}\,,\ \  
\hat{\calw}_{2\pi,0}={\calw}_{2\pi,0}\cap \hat{\calh}_{2\pi,0}\,.
$$

\begin{lemmaA} {\sl Let {\rm (H2)} and {\rm (H3)} hold. Then, the operator $\mathscr Q$ maps  $\hat{\calw}_{2\pi,0}$ onto $\hat{\calh}_{2\pi,0}$ homeomorphically.}
\label{3.2_00}
\end{lemmaA}
{\em Proof.} By (H3), $\mathscr Q$ is Fredholm of index 0, whereas by Lemma A.\ref{3.1_00} ${\rm dim}\,{\sf N}\,[\mathscr Q]=2$. From  classical theory of Fredholm operators (e.g. \cite[Proposition 8.14(4)]{Z}) it then follows that ${\rm dim}\,{\sf N}\,[\mathscr Q^\dagger]=2$ 
where
$$
\mathscr Q^\dagger=\zeta_0(\cdot)_\tau+L_2^\dagger
$$
is the adjoint of $\mathscr Q$. In view of the stated properties of $v_0^\dagger$, we infer that ${\rm span}\,\{v_1^\dagger,v_2^\dagger\}={\sf N}\,[\mathscr Q^\dagger]$, and the lemma follows from another classical result on Fredholm operators (e.g. \cite[Proposition 8.14(2)]{Z}).\QED
\smallskip\par 
Let
$$
L_2(\mu):=L_2+\mu\,S_{011}\,.
$$
since, by (H2), $\nu_0$ is a simple eigenvalue of $L_2(0)\equiv L_2$,  denoting by $\nu(\mu)$ the eigenvalues of $L_2(\mu)$, it follows (e.g.  \cite[Proposition 79.15 and Corollary 79.16]{Z1}) that in a neighborhood of $\mu=0$ the map $\mu\mapsto\nu(\mu)$ is well defined and of class $C^\infty$ and, further,
$$
\nu'(0)=\langle v_0^\dagger, S_{011}(v_0)\rangle
.
$$
Using the latter, by  direct inspection we show
\be   
\Re[\nu'(0)]=\pi^{-1}(S_{011}(v_1)|v_1^\dagger)\,.
\eeq{nu0}
\par
{\bf Proof of Theorem \ref{3.1_ar}(a)}. 
In order to ensure that the solutions we are looking for are non-trivial, we endow \eqref{Ar.5} with the side condition 
\be  
(w|v_1^\dagger)=\varepsilon\,,\ \ (w|v_1^\dagger)=0\,, 
\eeq{Ar.8}
where $\varepsilon$ is a real parameter ranging in a neighborhood, $\cali(0)$, of $0$.
We next scale $v$ and $w$  by setting
$v=\varepsilon\,{\sf v}$, ${w}=\varepsilon\, {\sf w}$, so that problem \eqref{Ar.5},
 \eqref{Ar.8} becomes
\be\ba{ll}\medskip  
L_1({\sf v})={\mathcal N}_1(\varepsilon, {\sf v},{\sf w},\mu)\,,\ \mbox{in $\calv$}\,;\\ 
\zeta_0\, {\sf w}_\tau +L_2({\sf w)={\mathcal N}_2(\varepsilon, \zeta, {\sf v},{\sf w},\mu)\,,\ \mbox{in $\calz_{2\pi,0}$}}\,,
\ \
({\sf w}|v_1^\dagger)=1\,,\ \ ({\sf w}|v_1^\dagger)=0\,,
\ea
\eeq{3.11}
where 
$$\ba{ll}\medskip
\mathcal N_1(\varepsilon,\sfv,\sfw,\mu):=(1/\varepsilon)\,N_1(\varepsilon{\sf v}, \varepsilon\sfw,\mu)\,, \\
\mathcal N_2(\varepsilon,\zeta,\sfv,\sfw,\mu):=(1/\varepsilon)\,N_2(\varepsilon{\sf v}, \varepsilon\sfw,\mu)+(\zeta_0-\omega){\sf w}_\tau\,.\ea
$$ Set ${\sf U}:=(\mu, \zeta,{\sf v},{\sf w})$, and consider the map
$$\ba{cc}\medskip
F: (\varepsilon, {\sf U}):=(\varepsilon,{\sf U})\in \cali(0)\times \left(U(0)\times V(\zeta_0)\times \calu\times \calw_{2\pi,0}\right)\\ \smallskip
\mapsto
\Big( L_1(\sfv)-\mathcal N_1(\varepsilon,\sfv,\sfw,\mu),\ 
\mathscr Q(\sfw)-\mathcal N_2(\varepsilon,\zeta,\sfv,\sfw,\mu),\
({\sf w}|v_1^\dagger)-1,\ ({\sf w}|v_2^\dagger)\Big)
\in \calv\times \calz_{2\pi,0}\times \real^2\,,
\ea
$$
with $U(0)$ and $V(\zeta_0)$  neighborhoods of 0 and  $\zeta_0$, respectively. Since, by (H4), we have in particular $\mathcal N_1(0,0,v_1,0)=
\mathcal N_2(0,\zeta_0,v_1,0)=0$, using  
\eqref{3.7_00}$_1$ and Lemma A.\ref{3.1_00} we deduce that, at $\varepsilon=0$,  
the equation $F(\varepsilon,{\sf U})=0$ has the solution ${\sf U}_0=(0,\zeta_0,0,v_1)$. Therefore, since by (H4) we have that $F$ is analytic at $(0,{\sf U_0})$, by the  analytic version of the Implicit Function Theorem (e.g. \cite[Proposition 8.11]{Z}), to show the existence statement -including the validity of \eqref{Ar.10}- it suffices to show that the Fr\'echet derivative, $D_{{\mbox{\tiny {\sf U}}}}F(0,{\sf U}_0)$ is a bijection. 
Now,   in view of the assumption (H4), it easy to see that the Fr\'echet derivative of $\mathcal N_1$ at $(\varepsilon=0, \sfv=0,\sfw=v_1,\mu=0)$ is equal to 0, while that of $\mathcal N_2$ at $(\varepsilon=0, {\sf U}={\sf U}_0)$ is equal to $ -\zeta\,(v_{1})_\tau+\mu\,S_{011}(v_1)$\,.
Therefore, $D_{{\mbox{\tiny {\sf U}}}}F(0,{\sf U}_0)$ is a bijection  if  we  prove that for any
$({\sf f}_1,{\sf f}_2,{\sf f}_3,{\sf f}_4)\in \calv\times\calz_{2\pi,0}\times\real\times\real$, the following set of equations has one and only one solution $(\mu,\zeta,\sfv,\sfw)\in \real\times\real\times \calu\times \calw_{2\pi,0}$:
\be\ba{rl}\medskip
{L}_1(\sfv)=&\!\!\!\!{\sf f}_1\ \ \mbox{in $\calv$}\\ \medskip
\mathscr Q(\sfw)= &\!\!\!\!-\zeta\, (v_{1})_\tau+\mu\,S_{011}(v_1)+{\sf f}_2\ \ \mbox{in $\calz_{2\pi,0}$}\,,\\
({\sf w}|\bfv_1^\dagger)=&\!\!\!\!{\sf f}_3\,,\ \ ({\sf w}|\bfv_2^\dagger)={\sf f}_4 \ \ \mbox{in $\real$}\,,
\ea
\eeq{3.12}
In view of (H1), for any given ${\sf f}_1\in \calv$, equation \eqref{3.12}$_1$ has one and only one solution $\sfv\in \calu$. Therefore, it remains to prove the existence and uniqueness property only for the system of equations \eqref{3.12}$_{2-4}$
To this aim, we observe that, by Lemma A.\ref{3.2_00}, for a given ${\sf f}_2\in \calz_{2\pi,0}$,  equation \eqref{3.12}$_2$ possesses a unique solution $\sfw_1\in\hat{\calw}_{2\pi,0}$ if and only if its right-hand side is in $\hat{\calz}_{2\pi,0}$, namely,
$$
\big(-\zeta\, (v_{1})_\tau+\mu\,S_{011}(v_1)+{\sf f}_2|v_1^\dagger\big)=\big(-\zeta\,( v_{1})_\tau+\mu\,S_{011}(v_1)+{\sf f}_2|v_2^\dagger\big)=0\,. 
$$
Taking into account \eqref{3.7_00}$_{2}$ the above conditions will be satisfied provided we can find $\mu$ and $\zeta$  satisfying the following algebraic system
\be\ba{rl}\medskip
\mu(\,S_{011}(v_1)|v_1^\dagger)&\!\!\!\!=-({\sf f}_2|v_1^\dagger)\\
\zeta+\mu\,(\,S_{011}(v_1)|v_2^\dagger)&\!\!\!\!=-({\sf f}_2|v_2^\dagger)\,.
\ea
\eeq{3.13}
However, by virtue of \eqref{nu0}, \eqref{nupr} this system possesses a uniquely determined solution $(\mu,\zeta)$, which ensures the existence of a unique solution ${\sf w}_1\in \hat{\calw}_{2\pi,0}$ to \eqref{3.12}$_2$ corresponding to the  selected values of $\mu$ and $\zeta$. We now set
$$ 
\sfw:={\sf w}_1+\alpha\,v_1+\beta\,v_2\,,\ \ \alpha\,,\, \beta\in\real\,.
$$
Clearly, by Lemma A.\ref{3.1_00}, ${\sfw}$ is also a solution to \eqref{3.12}$_2$. We then 
choose $\alpha$ and $\beta$ in such a way that $\sfw$ satisfies both conditions \eqref{3.12}$_{3,4}$  for any given ${\sf f}_i\in\real$, $i=1,2$. This choice is made possible by virtue of \eqref{3.7_00}$_1$.
We have thus shown that $D_{\mbox{\tiny {\sf U}}}F(0,{\sf U}_0)$ is surjective. To show that it is also injective,
set ${\sf f}_i=0$ in \eqref{3.12}$_{2-4}$. From \eqref{3.13} and \eqref{nu0}, \eqref{nupr} it then follows $\mu=\zeta=0$ which in turn implies, by \eqref{3.12}$_2$ and  Lemma A.\ref{3.1_00}, $\sfw=\gamma_1\,v_1+\gamma_2\,v_2$, for some $\gamma_i\in\real$, $i=1,2$. Replacing this information back in \eqref{3.12}$_{3,4}$ with ${\sf f}_3={\sf f}_4=0$, and using \eqref{3.7_00}$_1$ we conclude $\gamma_1=\gamma_2=0$, which proves  the claimed injectivity property.
Thus, $D_{\mbox{\tiny {\sf U}}}F(0,{\sf U}_0)$ is a bijection, and the proof of the existence statement in (a) is completed. 
\par
{\bf Proof of Theorem \ref{3.1_ar}(b)}.
Let $(z,s)\in \calu\times\calw_{2\pi,0}$  be a  $2\pi$-periodic solution to \eqref{Ar.5} with $\zeta\equiv\tilde{\zeta}$ and $\mu\equiv\tilde\mu$.
By the uniqueness property associated with  the implicit function theorem, the proof of the claimed uniqueness
amounts to show that we can find a sufficiently small $\rho>0$ such that if
\be
\|z\|_{\calu}+\|s\|_{\calw_{2\pi,0}}+|\tilde\zeta-\zeta_0|+|\tilde\mu|<\rho\,,
\eeq{3.14}
then there exists a neighborhood of $0$, $\cali(0)\subset\real$, such that
\be\ba{cc}\medskip
s=\eta\, v_1+\eta\,{\sf s}\,,\, \ z=\eta\,{\sf z}\,, \ \mbox{for all $\eta\in\cali(0)$},\, \\
|\tilde\zeta-\zeta_0|+|\tilde\mu|+\|{\sf z}\|_{\calu}+\|{\sf s}\|_{\calw_{2\pi,0}}\to 0\ \ \mbox{as $\eta\to 0$}\,.\ea
\eeq{3.15}
To this end, we notice that, by \eqref{3.7}$_1$, we may write
\be
s={\sigma} +\ts
\eeq{3.16}
where ${\sigma}=(s|v^\dagger_1)\,v_1+(s|v^\dagger_2)\,v_2$ and 
\be
(\tilde{\sf s}|v^\dagger_i)=0\,,\ \ i=1,2\,.
\eeq{3.17}
We next make the simple but important observation that if we modify $s$ by a  constant phase shift in time, $\delta$, namely, $s(\tau)\to s(\tau+\delta)$,  the shifted function is still  a $2\pi$-periodic solution to \eqref{Ar.5}$_2$ and, moreover, by an appropriate choice of $\delta$, 
\be   
{\sigma}=\eta\, v_1\,,
\eeq{3.18}
with $\eta=\eta(\delta)\in\real$. (The proof of \eqref{3.18} is straightforward, once we take into account the definition of $v_1$ and $v_2$.)
Notice that from \eqref{3.14}, \eqref{3.16}--\eqref{3.18} it follows that
\be
|\eta| +
\|{\ts}\|_{\calw_{2\pi,0}}
\to 0 \ \ \mbox{as $\rho\to 0$}\,.
\eeq{3.19} 
From \eqref{3.5} we thus get 
\be 
{L}_1(z)= N_1(z,\eta\, v_1+\ts,\tilde\mu) 
\eeq{3.20}
and, recalling Lemma A.\ref{3.1_00}, 
\be
\mathscr Q(\ts)=\eta(\zeta_0-\zeta)(v_1)_\tau+(\zeta_0-\zeta)\ts_\tau+N_2(z,\eta\,v_1+\ts,\tilde\mu)\,.
\eeq{3.21}
In view of (H4) and \eqref{3.14}, we easily deduce deduce
$$
N_1(z,\eta\, v_1+\ts,\tilde\mu)=R_{110}z(\eta\,v_1+\ts)+R_{101}z\tilde\mu+R_{020}(\eta\,v_1+\ts)^2+n_1(z,\eta,\ts,\tilde\mu)\,,
$$
where 
$$
\|n_1(z,\eta,\ts,\tilde\mu)\|_{\calv}\le \epsilon(\rho) \,\left(\|z\|_{\calu}+\|\ts\|_{\calw_{2\pi,0}}+\eta^2\right)\,,\ \ \
\epsilon (\rho)\to 0\ \mbox{as}\ \rho\to 0\,,$$
so that, by \eqref{3.20} and (H1) we obtain by taking $\rho$ sufficiently small
\be
\|z\|_{\calu}\le c_1\,\big(|\eta|^2+\|\ts\|_{\calw_{2\pi,0}}^2+\epsilon(\rho)\|\ts\|_{\calw_{2\pi,0}}\big)\,.
\eeq{3.22}
Likewise,
\be\ba{rl}\medskip
N_2(z,\eta\,v_1+\ts,\tilde\mu)
=&\!\!\!\!S_{011}(\eta\,v_1+\ts)\tilde\mu+S_{110}z(\eta\,v_1+\ts)+S_{101}z\tilde\mu\\&+S_{200}z^2+S_{020}(\eta\,v_1+\ts)^2+n_2(z,\eta,\ts,\tilde\mu)\,,
\ea
\eeq{3.23}
where $n_2$ enjoys the same property as $n_1$.
From \eqref{3.21},  \eqref{3.23} and \eqref{3.7_00}$_1$ we infer, according to Lemma A.\ref{3.2_00}, that the
following (compatibility) conditions must be satisfied
$$\ba{ll}\medskip
-\eta\,\tilde\mu\,(S_{011}(v_1)|v_1^\dagger)=\big((\zeta_0-\zeta)\ts_\tau+ S_{011}\ts\tilde\mu+S_{110}z(\eta\,v_1+\ts)|v_1^\dagger\big)\\ \medskip
\hspace*{3.5cm}+\big(S_{200}z^2+S_{020}(\eta\,v_1+\ts)^2|v_1^\dagger\big)+(n_2|v_1^\dagger)\\ \medskip
\eta\,(\zeta-\zeta_0)=\big((\zeta_0-\zeta)\ts_\tau+ S_{011}\ts\tilde\mu+S_{110}z(\eta\,v_1+\ts)|v_2^\dagger\big)\\ \medskip
\hspace*{3.5cm}+\big(S_{200}z^2+S_{020}(\eta\,v_2+\ts)^2|v_2^\dagger\big)++(n_2|v_2^\dagger)\,,
\ea
$$
so that, from \eqref{3.9} and the property of $n_2$  we show
\be \ba{ll}\medskip
|\eta|\,\big(|\tilde\mu|+|\zeta-\zeta_0|\big)\le c_2 \big(|\zeta-\zeta_0|+|\tilde\mu|\big)\,\|\ts\|_{\calw_{2\pi,0}}
+|\eta|\,\|z\|_{\calu}+\|z\|_{\calu}^2\\
\hspace*{3.5cm}+\|\ts\|_{\calw_{2\pi,0}}^2+\eta^2\big)+\epsilon(\rho)\big(\|z\|_{\calh}+\|\ts\|_{\calw_{2\pi,0}}\big)\,.\ea
\eeq{3.24}
Also, applying Lemma A.\ref{3.2_00} to \eqref{3.21} and using \eqref{3.23}, \eqref{3.14} with $\rho$ sufficiently small we get
\be
\|\ts\|_{\calw_{2\pi,0}}\le c_3\,\big(|\eta|\,(|\tilde\mu|+|\zeta-\zeta_0|)+(|\eta|+|\tilde\mu|+\epsilon(\rho))\,\|z\|_{\calu}+\|z\|_{\calu}^2+\eta^2\big)\,.
\eeq{3.25}
Summing side by side \eqref{3.22}, \eqref{3.24} and $(1/(2c_3))\times$\eqref{3.25}, and taking again $\rho$ small enough, we thus arrive at
$$
|\eta|\,\big(|\tilde\mu|+|\zeta-\zeta_0|\big)+\|z\|_{\calu}+\|\ts\|_{\calw_{2\pi,0}}\le c_4\,\eta^2\,,
$$
from which we establish the validity of \eqref{3.15}$_2$, thus concluding the proof of the uniqueness property (b).
\par 
{\bf Proof of Theorem \ref{3.1_ar}(c)}
We notice that if $\big(v(-\varepsilon),w(-\varepsilon;\tau)\big)$ is the solution corresponding to $-\varepsilon$, we have $\big(w(-\varepsilon;\tau+\pi)|v_1^\dagger\big)=\varepsilon \,v_1$, which, by part (b), implies that, up to a phase shift, $\big(v(-\varepsilon),w(-\varepsilon;\tau)\big)=\big(v(\varepsilon),w(\varepsilon;\tau)\big)$. This, in turn, furnishes $\zeta(-\varepsilon)=\zeta(\varepsilon)$ and $\mu(-\varepsilon)=\mu(\varepsilon)$. From the latter and the analyticity of $\mu$ we then obtain
that either $\mu\equiv0$ or else there is an integer $k\ge 1$ such that 
$$
\mu(\varepsilon)=\varepsilon^{2k}\mu_k+O(\varepsilon^{2k+2})
\ \ \mu_k
\in\real\backslash\{0\}\,.
$$
Thus, $\mu(\varepsilon)<0$ or $\mu(\varepsilon)>0$, according to whether $\mu_k$ is negative or positive. The theorem is completely proved.     
\QED

\ed